\documentclass[11pt]{combine}
\usepackage{amsmath,amsfonts,amssymb,amsthm,epsfig, graphicx, hyperref}
\usepackage[dvipsnames,table,dvipsnames*, svgnames*, hyperref]{xcolor}
\usepackage{natbib,amsmath,amssymb,amsthm,graphicx,setspace,paralist,booktabs,rotating,subcaption,float,color}
\usepackage{multirow}
\usepackage{fancyhdr}
\usepackage{bibunits}
\usepackage[small]{titlesec}

\newtheorem{theorem}{Theorem}
\newtheorem{corollary}{Corollary}
\newtheorem{proposition}{Proposition}

\newtheorem{lemma}{Lemma}
{
\theoremstyle{definition}
\newtheorem{definition}{Definition}
\newtheorem{example}{Example}
\newtheorem{remark}{Remark}
}
\newcommand{\beq}{\begin{equation}}
\newcommand{\eeq}{\end{equation}}
\newcommand{\beas}{\begin{align*}}
\newcommand{\eeas}{\end{align*}}
\newcommand{\bea}{\begin{align}}
\newcommand{\eea}{\end{align}}
\newcommand{\bei}{\begin{itemize}}
\newcommand{\eei}{\end{itemize}}
\newcommand{\ben}{\begin{enumerate}}
\newcommand{\een}{\end{enumerate}}
\newcommand{\bet}{\begin{theorem}}
\newcommand{\eet}{\end{theorem}}
\newcommand{\bel}{\begin{lemma}}
\newcommand{\eel}{\end{lemma}}
\newcommand{\bep}{\begin{proposition}}
\newcommand{\eep}{\end{proposition}}

\newcommand{\bed}{\begin{definition}}
\newcommand{\eed}{\end{definition}}
\newcommand{\bec}{\begin{corollary}}
\newcommand{\eec}{\end{corollary}}
\newcommand{\bex}{\begin{example}}
\newcommand{\eex}{\end{example}}

\newcommand{\R}{\mathbb{R}}
\newcommand{\E}{\mathbb{E}}

\newcommand{\argmax}{\mathop{\rm arg\max}}

\def\xx{\bold{x}}

\numberwithin{equation}{section}

\begin{document}


\title{\scshape Optimal Permutation Recovery in Permuted Monotone Matrix Model}
\author{Rong Ma$^1$, T. Tony Cai$^2$ and Hongzhe Li$^1$ \\
Department of Biostatistics, Epidemiology and Informatics$^1$\\
Department of Statistics$^2$\\
University of Pennsylvania\\
Philadelphia, PA 19104}
\date{}
\maketitle
\thispagestyle{empty}

	\begin{abstract}
Motivated by  recent research on quantifying  bacterial growth dynamics based on genome assemblies, we consider a permuted monotone matrix model $Y=\Theta\Pi+Z$, where the rows represent different  samples, the columns represent contigs in genome assemblies and the elements represent log-read counts after preprocessing steps and Guanine-Cytosine (GC) adjustment.  In this model, $\Theta$ is an unknown mean matrix with monotone entries for each row, $\Pi$ is a permutation matrix that permutes the columns of $\Theta$, and $Z$ is a noise matrix. This paper studies the problem of  estimation/recovery of $\Pi$ given the observed noisy matrix $Y$. We propose an estimator based on the best linear projection, which  is shown to be minimax rate-optimal for both exact recovery, as measured by the 0-1 loss, and partial recovery, as quantified by the normalized Kendall's tau distance. Simulation studies demonstrate the superior empirical performance of the  proposed estimator over alternative methods. We demonstrate the methods  using  a synthetic metagenomics dataset of 45 closely related bacterial species  and a real metagenomic dataset to compare the bacterial growth dynamics between the responders and the non-responders of the IBD patients after 8 weeks of treatment.
	\bigskip
	
	\noindent\emph{KEY WORDS}: Kendall's tau; Microbiome growth dynamics;  Minimax lower bound;  Sorting.
\end{abstract}

\section{INTRODUCTION}

\subsection{A Motivation Example from Microbiome Studies}

The statistical problem considered in this paper is motivated  by the problem of estimating the  bacterial growth dynamics based on shotgun metagenomics data \citep{myhrvold2015distributed,abel2015sequence,korem2015growth, brown2016measurement}.  The growth dynamics of microbial populations reflects their physiological states and drives variation of microbial  compositions, which provide important feature summary of the microbes in a given community.  One way of studying such communities is through shotgun metagenomic sequencing, which  involve direct DNA sequencing of all the microbiome genomes in a given microbial community.  \cite{korem2015growth} presented the first paper on  quantifying the bacterial growth dynamics based on shotgun metagenomics data, where the  uneven sequencing read coverage resulting from the bidirectional DNA replications provides information on the rates of microbial DNA replications.  For  bacterial species with known complete genome sequences,  \cite{korem2015growth} proposed to use the peak-to-trough ratio (PTR) of read coverages to quantify the bacterial growth dynamics after aligning the sequencing reads to the complete genome sequences. 

However,  in many applications, it is of importance to quantify the bacterial growth dynamics based on genome assemblies for the bacterial species with unknown  genomes.  These genome assemblies may represent new bacterial species that we have seen or sequenced before. The genome assembly of a bacterium species consists of a collection of contigs (called bin) constructed based on the overlapping of the sequencing reads \citep{megahit, maxbin}. Compared to the complete genome, the genome assembled bins are more fragmented and often contained errors or contaminations.   The noisy read coverage data due to intraspecific variations, interspecific/intraspecific repeated sequences, limited sequencing depths and the inability of binning algorithms to correctly cluster all the contigs  further complicate the estimation of growth dynamics based on read coverages of the contigs.   Besides these noisy count data, one key difficulty in estimating the growth dynamic based on contig counts is that
the accurate locations of the contigs on the original genome are  unknown. It is therefore not feasible to measure the microbial  growth rate directly using peak-to-trough coverage ratio  for the assembled genomes \citep{brown2016measurement,gao2018quantifying}.

\cite{brown2016measurement} presented the first method (called iRep) of estimating the bacterial growth dynamics based on genome assemblies, where the contigs are ordered based on the GC-adjusted counts for each sample separately. However, due to noise in the count data, such an ordering method often leads to wrong ordering of the contigs and therefore inaccurate estimates of the growth dynamics.  \cite{gao2018quantifying} developed  a computational algorithm, DEMIC, to accurately compare growth dynamics of a given assembled species  existing in multiple samples by taking advantage of highly fragmented contigs assembled in typical metagenomics studies. One key step of DEMIC is to apply   a  principal components analysis (PCA)-based method  to recover the true ordering of the contigs along the underlying unknown bacterial complete genomes.   \cite{gao2018quantifying} reported  excellent  empirical performance of DEMIC over existing methods. The goal of this paper is to provide a rigorous statistical framework  to study the problem of optimal permutation recovery in a permuted monotone matrix model. 

\subsection{A Permuted Monotone Matrix Model}
For a given genome assembly with  $p$ contigs, DEMIC first obtains the read coverage for each of the sliding window of size 5000 bps, denoted by $X_{ijl}$ for the $i$th sample, $j$th contig and $k$th window. In order to account for the GC-content of the $k$th window, \cite{gao2018quantifying} considered the following mixed-effects model, 
\[
\log_2 X_{ijk}=\alpha+GC_{jk}\beta+W_{ij}+e_{ijk},
\]
where  $GC_{jk}$ is the centred GC count of the $k$th window of the $j$th contig, $W_{ij}$ is the sample- and contig- specific random intercept, $\alpha$ is the intercept, $\beta$ is the regression coefficient, and $e_{ijk}$ is the random error. This model is fitted for each contig to obtain the best linear unbiased predictor of $W_{ij}$, which is used as the GC-adjusted log-read count $Y_{ij}$ for the $i$th sample and $j$th contig.  Here $Y_{ij}$ can be regarded as average read coverage over non-overlapping windows of a contig and is approximately normally distributed.

Let $Y$ be the GC-adjusted log-contig count matrix of $n$ samples and $p$ contigs of a genome assembly with $Y_{ij}$ as its entries. Given this, we consider the following permuted monotone matrix model:
\beq \label{m2}
Y = \Theta \Pi+Z,
\eeq
where $\Theta\in \R^{n\times p}$ is an unknown nonnegative signal matrix with nondecreasing rows, $Z\in \R^{n\times p}$ is a zero-mean  noise matrix, and $\Pi \in \R^{p\times p}$ is a permutation matrix corresponding to some permutation $\pi$ from the symmetric group $\mathcal{S}_p$. That is, after a suitable permutation of the columns of $Y$, all the rows of the mean matrix are nondecreasing sequences.   In microbiome applications, 
$\Theta$ is the matrix of true log-coverage of $n$ samples over $p$ contigs along the circular genome of the bacterium, which is generally hypothesized to have non-decreasing rows.  $\Pi$  represents a permutation due to unknown locations of the contigs relative to the replication origin. 
Throughout this paper,  we denote the parameter space
\[
(\Theta,\pi)\in\mathcal{D} = \bigg\{ \Theta = (\theta_{ij})\in \R^{n\times p}, \pi\in \mathcal{S}_p: \, 0\le \theta_{i,j-1}\le \theta_{i,j} <\infty \text{ for all $1\le i\le n, 2\le j\le p$}
\bigg\}.
\]
The focus of this paper is to optimally estimate the permutation $\pi$ from the noisy observation $Y$.

\subsection{Related Problems and Other Applications}

The permutation recovery problem under permuted monotone matrix model bears some similarity to other problems studied in machine learning literature, including 
the  feature matching between two sets of observations  \citep{collier2016minimax}  and 
linear regression model with permuted data,  where the correspondences between the response and the predictors are unknown  \citep{pananjady2016linear, slawski2017linear,pananjady2017denoising}. More recently, \cite{flammarion2019optimal} considered the problem of statistical seriation, which has a close affinity to our model (\ref{m2}). However, the focus of \cite{flammarion2019optimal} is to optimally estimate the signal matrix $\Theta$ rather than the underlying permutation.

Model (\ref{m2}) can be thought as a natural extension of the shape constrained matrix denoising model studied in the isotonic regression literature. Specifically, under Model (\ref{m2}) with known $\Pi={\bf I}_p$, risk bounds and the minimax rate-optimal estimator for $\Theta$ under the Frobenius norm was obtained in \cite{chatterjee2015risk} for $n=1$ and later in \cite{chatterjee2018matrix} for general $n>1$. 
Using the idea of optimal transport, a minimax optimal estimator of the underlying signals was obtained by \cite{rigollet2018uncoupled}.  However, their goal is not to recover the underlying permutation.

Besides the microbiome applications, the  permuted monotone matrix model is generic and has other applications.   For instance, the problem of permutation recovery is usually equivalent to statistical ranking/sorting from noisy observations, which arises commonly in finance \citep{currie2011finance}, sport analytics \citep{deshpande2016estimating}, and recommendation systems \citep{rendle2009learning}. Specifically, in the latter case, the task of tag recommendation is to provide a user with a personalized ranked list of tags for a specific item. Under the permuted monotone matrix model, we can treat the entries of $Y$, say $Y_{ij}$, as an indicator of the $j$th tag being related to the $i$th item by a given customer,  and $\Theta$ as  a probability matrix characterizing the customer's tagging preferences across multiple items. As a result, recovering the underlying permutation provides a solution of a tag recommender.

\subsection{Main Contributions and Organization}

\begin{figure}[h!]
	\centering
	\includegraphics[angle=0,width=13cm]{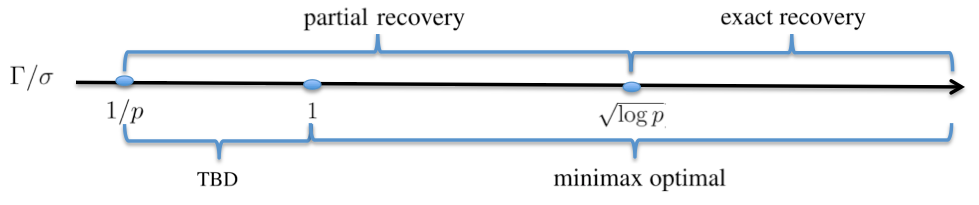}
	\caption{A graphical illustration of the main result obtained in this paper about the regions of the signal-to-noise ratio $\Gamma/\sigma$ that correspond to exact/partial recovery, and the region with proved minimax optimality.} 
	\label{rec}
\end{figure}

In this paper,  we investigate the problem  of  permutation recovery in the permuted monotone matrix model (\ref{m2}), which relies on certain invariance property of the singular subspace of the monotone matrices. The properties of the proposed method in terms of both the exact and partial recovery 
are studied in detail. In particular, we obtained regions of the signal-to-noise ratio (defined later as $\Gamma/\sigma$) that are subject to exact/partial recovery (Figure \ref{rec}). 
For both exact and partial permutation recovery, we obtained the matching minimax lower bounds and established the minimax rate-optimality of the proposed method over a wide range of parameter space (Figure \ref{rec}). For partial recovery, the proof of the lower bound relies on a version of Fano's lemma and the sphere packing of the symmetric group equipped with the Kendall's tau metric.

The rest of this paper is organized as follows. After a brief introduction of notation and definitions, we present in Section 2 the proposed permutation estimator. The theoretical properties of the proposed method are studied, first under a more illustrative linear growth model in Section 3 and then under a general growth model in Section 4. Section 5 provides results on minimax lower bounds and the optimality of the proposed estimator.  We evaluate the methods using both simulated data, synthetic and real  microbiome  datasets and compare with other methods in Section 6. In Section 7, we discuss some implications and extensions of the methods. Finally, the proofs of our main results are given in Section 8.

\subsection{Notation and Definitions}

Throughout, we define the permutation $\pi$ as a bijection from the set $\{1,2,...,p\}$ onto itself. For simplicity, we denote $\pi=(\pi(1),\pi(2),...,\pi(p))$. All permutations of the set $\{1,2,...,p\}$ form a symmetric group, equipped with the function composition  operation $\circ$, denoted as $\mathcal{S}_p$. For any $\pi\in \mathcal{S}_p$, we denote $\pi^{-1}\in \mathcal{S}_p$ as its group inverse, so that $\pi\circ\pi^{-1}=\pi^{-1}\circ\pi=id$, and denote $\text{rev}(\pi)=(\pi(p),\pi(p-1),...\pi(1))$. In particular, we may use $\pi$ and its corresponding permutation matrix $\Pi\in \R^{p\times p}$ interchangeably, depending on the context. For a vector $\bold{a} = (a_1,...,a_n)^\top \in \mathbb{R}^{n}$, we define the $\ell_p$ norm $\| \bold{a} \|_p = \big(\sum_{i=1}^n a_i^p\big)^{1/p}$, and the $\ell_\infty$ norm $\| \bold{a}\|_{\infty} = \max_{1\le j\le n}  |a_{i}|$.  For a matrix $\Theta\in \R^{p_1\times p_2}$, we denote $\Theta_{.i}\in \R^{p_1}$ as its $i$-th column, $\Theta_{i.}\in \R^{p_2}$ as its $i$-th row, and denote its (ordered) singular values as $\lambda_1(\Theta)\ge \lambda_2(\Theta)\ge...\ge \lambda_{\min\{p_1, p_2\}}(\Theta)$. Furthermore, for sequences $\{a_n\}$ and $\{b_n\}$, we write $a_n = o(b_n)$ if $\lim_{n} a_n/b_n =0$, and write $a_n = O(b_n)$, $a_n\lesssim b_n$ or $b_n \gtrsim a_n$ if there exists a constant $C$ such that $a_n \le Cb_n$ for all $n$. We write $a_n\asymp b_n$ if $a_n \lesssim b_n$ and $a_n\gtrsim b_n$. For a finite set $A$, we denote $|A|$ as its cardinality. We use the logical symbols $\land$ and $\lor$ to represent ``and" and ``or," respectively. 
Lastly, $C, C_0, C_1,...$ are constants that may vary from place to place.

\section{PERMUTATION RECOVERY VIA BEST LINEAR PROJECTION}

In the following, we first make some key observations about the connection between the underlying permutation $\pi$ and the column linear projections of the observed matrix $Y$, which motivate our construction of the proposed estimator.

\subsection{Linear Projection}
Given the observed noisy matrix $Y$, we consider the class of the linear projection statistics of the form $w^\top Y\in \R^{p}$ where $w\in \R^n$ and $\|w\|_2=1$. Intuitively, by projecting each column of $Y$ onto the subspace generated by $w$, the components of $w^\top Y$ (hereafter referred as ``projection scores") would quantify the relative position of the columns of $Y$, so that their order statistics can be used to recover the original orders of the columns of $\Theta$. To fix ideas, we define the following ranking operator.
\bed [Ranking Operator]
The ranking operator $\frak{r}: \R^p \to \mathcal{S}_p$ is defined such that for any vector $\xx\in \R^p$, $\frak{r}(\xx)$ is the vector of ranks for components of $\xx$ in increasing order. Whenever there are ties, increasing orders are assigned from left to right.
\eed
For example, given a vector $\xx = (2,5,1,6,2)^\top$, we have $\frak{r}(\xx) = (2,4,1,5,3)$. The following proposition concerning the invariance property of the column spacing of $\Theta$ is the key to our construction of the minimax optimal estimator.
\bep \label{prop.inv}
Suppose $(\Theta,\pi)\in \mathcal{D}$. For any nonnegative unit vector $w\in \R^n$, we have
\beq
\frak{r}(w^\top \Theta\Pi)=\pi^{-1}.
\eeq
\eep
Apparently, under the noiseless setting, any nonnegative unit vector would lead to the exact recovery of the underlying permutation as in this case the relative orders of the columns are exactly coded by the relative magnitudes of the projection scores $w^\top Y=w^\top \Theta\Pi$. However, with the noisy observations, $w^\top Y= w^\top \Theta\Pi+w^\top Z$ so that the relative orders of the columns are only partially preserved by the noisy projection scores $w^\top Y$, up to some random perturbations. 

Consequently, the best linear projection vector $w_0$ would correspond to the case where $w_0^\top \Theta\Pi$ has the most separated components such that their relative orders are most immune to the random noises. Specifically, since for any given $w\in \R^n$, the $i$-th component of $w^\top\Theta\Pi$ has the expression $w^\top \Theta\Pi e_i$ where $\{ e_i\}_{i=1}^p$ is the canonical basis of the Euclidean space $\R^p$, we define 
\[
w_0 = \argmax_{\substack{w\in \R^n\\ \|w\|_2=1}}\sum_{\substack{1\le i, j\le p\\i\ne j}} (w^\top \Theta\Pi e_i-w^\top \Theta\Pi e_j)^2= \argmax_{\substack{w\in \R^n\\ \|w\|_2=1}}\sum_{i=1}^p \bigg(w^\top \Theta\Pi e_i-\frac{1}{p}\sum_{j=1}^pw^\top \Theta\Pi e_j \bigg)^2,
\]
which maximizes the pairwise distances of the components under the squared distance. 
Now since $w_0$ relies on the unknown $\Theta\Pi$ and is not computable from the data, we substitute $\Theta\Pi$ by its sample/noisy counterpart $Y$ and define our data-driven best linear projection vector as
\beq \label{w.hat}
\hat{w}=\argmax_{\substack{w\in \R^n \\ \|w\|_2=1}}\sum_{i=1}^p \bigg(w^\top Y e_i-\frac{1}{p}\sum_{i=1}^pw^\top Y e_i\bigg)^2,
\eeq
which is actually the first eigenvector of the symmetric matrix 
\beq \label{A}
A = Y\sum_{i=1}^p\bigg(e_i-\frac{1}{p}\sum_{i=1}^pe_i\bigg)\bigg(e_i-\frac{1}{p}\sum_{i=1}^pe_i\bigg)^\top Y^\top,
\eeq
and can be immediately solved by performing an eigen-decomposition on $A$. 
Once $\hat{w}$ is obtained, we define our proposed permutation estimator as
\beq \label{p.est}
\hat{\pi}=(\frak{r}(\hat{w}^\top Y))^{-1}.
\eeq 
Intuitively, the projection vector $\hat{w}$ assigns different weights to the rows of $Y$ so that more weight is given to the rows whose elements are better separated and therefore more informative in  distinguishing the columns of $Y$ or $\Theta$.

\subsection{Evaluation Criteria}

The main focus of this paper is to investigate the theoretical properties of our proposed estimator (\ref{p.est}) under various loss measures and parameter spaces. For any given estimator $\check{\pi}$, we  first consider the 0-1 loss
\[
\ell(\check{\pi},\pi) = 1\{ \check{\pi}\ne \pi\},
\]
with the corresponding risk
$\E \ell(\check{\pi},\pi) =P(\check{\pi}\ne \pi)$.
The 0-1 loss is used to evaluate the exact recovery, which can be  a  strong requirement  for practical applications. As an alternative, we also consider the more flexible partial recovery, where the loss function is given by the normalized Kendall's tau distance \citep{kendall1938new} defined as
\beq \label{tau.loss}
\tau_K(\pi_1,\pi_2) = \frac{\{ \text{$\#$ of discordant pairs between $\pi_1$ and $\pi_2$}\}}{{n \choose 2}}.
\eeq
Technically, for two permutations $\pi_1$ and $\pi_2$, the set of discordant pairs is defined as
\[
\mathcal{G}(\pi_1,\pi_2)=\{ (i,j): i<j, [\pi_1(i)<\pi_1(j) \land \pi_2(i)>\pi_2(j)] \lor [\pi_1(i)>\pi_1(j) \land \pi_2(i)<\pi_2(j)]  \}
\]
so that the numerator in (\ref{tau.loss}) is equal to the cardinality $|\mathcal{G}(\pi_1,\pi_2)|$, which, in fact, is also the minimum number of pairwise adjacent transpositions converting $\pi_1^{-1}$ into $\pi_2^{-1}$ \citep{diaconis1988group}. The denominator ${n \choose 2}$ ensures that $\tau_K (\pi_1,\pi_2) \in [0,1]$ where $\tau_K (\pi_1,\pi_2)=0$ corresponds to $\pi_1=\pi_2$.

\section{A LINEAR GROWTH MODEL}

We start with a simpler case where the pair $(\Theta,\pi)$ is from the subspace
\begin{align} \label{linear.growth}
\mathcal{D}_L=\bigg\{ (\Theta,\pi)\in\mathcal{D}: \qquad \begin{aligned}
&\theta_{ij}=a_i\eta_j+b_i,
\text{ where $a_i,b_i\ge 0$ for $1\le i\le n$, }\\ 
&\text{$0\le \eta_j\le \eta_{j+1}$ for $1\le j\le p-1$}
\end{aligned} \bigg\}.
\end{align}
In other words, each row of $\Theta$ has a linear growth pattern with possibly different intercepts and different slopes.
In the context of bacterial growth dynamics, this model is sometimes referred as the Cooper-Helmstetter model \citep{cooper1968chromosome,bremer1977examination} that  associates the copy number of genes with their relative distances to the replication origin.  Specifically, $a_i$ is the ratio of genome replication time and doubling time,  which can be used to quantify the bacterial growth dynamics for the $i$th sample, $\eta_j$ is related to distance from the replication origin for the $j$th contig, and $b_i$ is related to the read counts at the replication origin and the sequencing depth. If the bacterium is non-dividing in sample $i$, $a_i$ is zero.

For the linear growth model \eqref{linear.growth},  there are two key quantities that are relevant to  permutation recovery. 
\bed
For any $\Theta\in \mathcal{D}_L$, we define 
\beq
\Gamma= \bigg(\sum_{i=1}^n a_i^2\bigg)^{1/2}\cdot \min_{1\le i< j\le p} | \eta_{i}-\eta_{j}  |
\eeq
as the \textit{local minimal signal gap} of $\Theta$, and define
\beq
\Lambda = \bigg(\sum_{i=1}^n a_i^2\bigg) \cdot\frac{1}{p}\sum_{1\le i<j\le p} (\eta_i-\eta_j)^2
= \bigg(\sum_{i=1}^n a_i^2\bigg) \cdot\sum_{j=1}^p (\eta_j-\bar{\eta})^2\eeq
as the \textit{global signal strength} of $\Theta$, where $\bar{\eta}=\sum_{j=1}^p {\eta_j}/p$.
\eed 

Intuitively, both quantities involve the set $\{|\eta_j-\eta_i| \}_{1\le i<j\le p}$ and the $\ell_2$ norm of the vector $a=(a_1,...,a_n)^\top$, which characterize the column spacings and the growth rates (slopes) of $\Theta$, respectively. Throughout this paper, we assume

\indent ({\bf A1}) the additive noise matrix $Z\in \R^{n\times p}$ has i.i.d. entries $z_{ij} \sim N(0,\sigma^2)$.\\
The Gaussian assumption simplifies our theoretical analysis. But this is not essential because all the theoretical results remain true if $Z$ has independent sub-Gaussian entries with parameters bounded by $\sigma^2$. 
The following theorem provides conditions on $\Gamma$ and $\Lambda$ such that exact recovery of $\pi$ can be obtained by $\hat{\pi}$ in (\ref{p.est}).

\bet[Exact Recovery, Linear] \label{linear.thm.1}
Suppose (A1) hold, $(\Theta,\pi)\in \mathcal{D}_L$ and $\Theta$ satisfies
\beq \label{cond.thm}
\Gamma >C_0 \sigma\sqrt{\log p}, \quad\Lambda >  C_1 \sigma^2(n\max\{ \sigma^2n/\Gamma^2,1\}+\sqrt{np\max\{ \sigma^2n/\Gamma^2,1\}})
\eeq
for some $C_0,C_1> 0$. Then with probability at least $1-O(p^{-c})$ for some constant $c>0$, up to a permutation reversion, we have $\hat{\pi}=\pi$.
\eet
\begin{remark}
	Due to non-identifiability between $\hat{w}$ and $-\hat{w}$ defined in (\ref{w.hat}), in Theorem \ref{linear.thm.1}, as well as all the other theoretical results concerning $\hat{\pi}$, the statement is up to a possible reversion of $\hat{\pi}$. For example,  for  permutation $\pi=(2,4,1,5,3)$,
	its reversion would be $\text{rev}(\pi)=(4,2,5,1,3)$. In fact, such indeterminacy can be avoided by noting that $a_i\ge 0$ for all $i$'s, but we will not pursue such a direction in this study as the practical interest only concerns relative orders of the permuted elements.
\end{remark}

Since $\Gamma$ depends on $a$ only through its $\ell_2$ norm $\|a\|_2$, the local minimal signal gap (MSG) condition $\Gamma \ge C\sigma\sqrt{\log p}$ allows for the presence of non-informative signals in the sense that some components of $a$ can be 0. In contrast, the condition on $\Lambda$ (GSS) depends on a trade-off between $\Gamma$ and $\sigma\sqrt{n}$. One the one hand, when $\Gamma > \sigma\sqrt{n}$, the condition on $\Lambda$ becomes $\Lambda \ge \sigma^2(C_0n+C_1\sqrt{np})$, which is independent of $\Gamma$, and is minimax optimal for left singular subspace estimation \citep{cai2018rate}. On the other hand, when $\Gamma < \sigma\sqrt{n}$, a stronger condition on $\Lambda$ is posed, as a compensation for small $\Gamma$. 

In some cases, the GSS condition in (\ref{cond.thm}) can be implied by the MSG condition. We summarize our results in the following proposition.

\bep \label{prop.1}
Suppose $\Gamma/\sigma>1/p$ and the MSG condition hold. Then the GSS condition can be implied by either one of the following conditions 
\begin{itemize}
	\item[(i)] $\Gamma \gtrsim \sigma\sqrt{n}$;
	\item[(ii)] $\Gamma \lesssim \sigma\sqrt{n}$, and either $(\sigma^4n^2/\Gamma^4)^{1/3}\lesssim p \lesssim \sigma^2n^2/\Gamma^2$ or  $p\gtrsim \sigma^2n^2/\Gamma^2+(\sigma^3n/\Gamma^3)^{2/5}$.
\end{itemize}
\eep

We next  turn to the partial recovery and study the rate of convergence of $\hat{\pi}$ measured by the normalized Kendall's tau distance under the linear growth model. In particular, we will assume an approximate uniform assignment of $\{\eta_j\}_{j=1}^p$ over some subinterval of $[0,\infty)$. In other words, the minimal element and maximal element of the set $\{|\eta_j-\eta_{j+1}| \}_{j=1}^{p-1}$ should have roughly the same magnitude, so that $\Gamma=\|a\|_2\cdot\min_{1\le j\le p-1}|\eta_{j}-\eta_{j+1}|\asymp\|a\|_2\cdot\max_{1\le j\le p-1}|\eta_{j}-\eta_{j+1}|$. This is equivalent to assuming that the contigs in genome assemblies are approximately uniformly spaced  along the circular genome.

\bet[Partial Recovery, Linear] \label{linear.thm.2}
Suppose (A1) hold, $(\Theta,\pi)\in \mathcal{D}_L$, and $\Theta$ satisfies
\begin{itemize}
	\item[(i)] there exist some $C_0>0$ such that $\max_{1\le j\le p-1}|\eta_{j}-\eta_{j+1}|< C_0\min_{1\le j\le p-1}|\eta_{j}-\eta_{j+1}|$
	for all $p>0$, and
	\item[(ii)] $\Lambda > C_1\sigma^2\big(\max\big\{ \frac{\sigma^2(n+\log p)^2}{\Gamma^2},n\}+\sqrt{p}\max\big\{ \frac{\sigma(n+\log p)}{\Gamma},\sqrt{n}\big\}\big)$ for some $C_1>0$.
\end{itemize}
Then, up to a permutation reversion,
\begin{align*}
\E [\tau_K(\hat{\pi},\pi) ] &\le 1\land \bigg(\frac{c_0\sigma}{p\Gamma}\min\bigg\{1,e^{-\Gamma^2/2\sigma^2}\log\bigg(1+\frac{2\sigma^2}{\Gamma^2} \bigg)  \bigg\}+\frac{c_1e^{-\Gamma^2/2\sigma^2}}{p(\Gamma/\sigma+\sqrt{8/\pi})}+\frac{c_2}{p^{c+2}}\bigg)
\end{align*}
for some $c,c_0,c_1,c_2>0$. 
\eet

\begin{figure}[h!]
	\centering
	\includegraphics[angle=0,width=13cm]{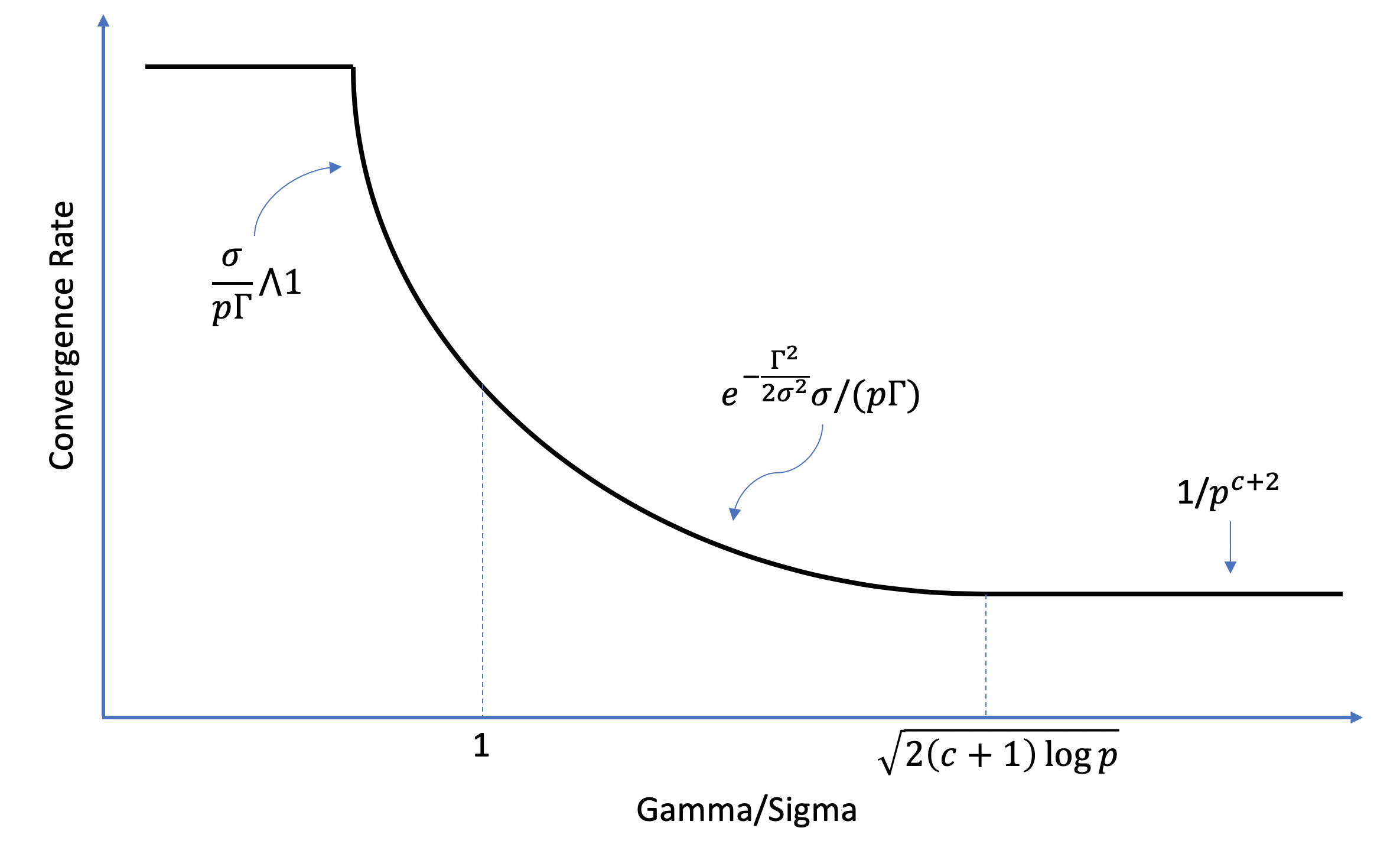}
	\caption{A graphical illustration of the risk upper bound for $\E [\tau_K(\hat{\pi},\pi) ]$, as a function of signal-to-noise ratio $\Gamma/\sigma$.} 
	\label{risk.curve}
\end{figure}

\begin{remark}
	The risk upper bound derived in the above theorem can be simplified as
	\[
	\E [\tau_K(\hat{\pi},\pi) ]\lesssim \left\{ \begin{array}{ll}
	\frac{\sigma}{p\Gamma} \land 1& \textrm{if $\Gamma/\sigma \to 0$}\\
	\frac{\sigma}{p\Gamma}e^{-\Gamma^2/2\sigma^2}+{1}/{p^{c+2}} & \textrm{otherwise}
	\end{array} \right.
	\]
	for some $c>0$. In the case of $\Gamma/\sigma \to \infty$, simple calculation yields $e^{-\Gamma^2/2\sigma^2}{\sigma}/{(p\Gamma)}+{1}/{p^{c+2}} \asymp e^{-\Gamma^2/2\sigma^2}{\sigma}/{\Gamma}$ when $\Gamma < \sigma\sqrt{2(c+1)\log p},$ whereas $e^{-\Gamma^2/2\sigma^2}{\sigma}/{(p\Gamma)}+{1}/{p^{c+2}} \asymp {1}/{p^{c+2}}$ when $\Gamma \ge \sigma\sqrt{2(c+1)\log p}.$ As a result, we also have
	\beq \label{rate}
	\E [\tau_K(\hat{\pi},\pi) ]\lesssim \left\{ \begin{array}{ll}
		{1}/{p^{c+2}} & \textrm{if $\Gamma/\sigma \ge \sqrt{2(c+1)\log p}$}\\
		\frac{\sigma}{p\Gamma}e^{-\Gamma^2/2\sigma^2}& \textrm{if $1\lesssim \Gamma/\sigma < \sqrt{2(c+1)\log p}$}\\
		\frac{\sigma}{p\Gamma} \land 1 & \textrm{if $\Gamma/\sigma \lesssim 1$}\\
	\end{array} \right..
	\eeq
	See Figure \ref{risk.curve} for an illustration. 
\end{remark}

In general, Theorem \ref{linear.thm.2} shows that, even with a weaker condition on $\Gamma$ that  is below the requirement for the exact recovery, our proposed estimator $\hat{\pi}$ is still able to obtain a partial recovery of $\pi$ with an exponential rate of convergence if $\Gamma/\sigma \gtrsim 1$ and a polynomial rate of convergence if $1/p<\Gamma/\sigma\lesssim 1$. As for $\Lambda$, the requirement is essentially the same as the exact recovery, except for an additional $\log p$ term, which is negligible in the exact recovery scenario.

Some implications about the practically preferable settings of $n$ and $p$ should be clarified. Firstly, although Theorem 1 implies that the difficulty for exact recovery increases as $p$ grows (see also Table 1 from our simulations), our theory suggests a wide range of feasible choices for $p$. For example, if the underlying signals $\theta_{ij}$ and the noise level $\sigma^2$ are of constant order, then we have $\Gamma\asymp  \sqrt{n}$ and $\Lambda \asymp np^3$, so the conditions of Theorem 1 imply that the exact recovery can be guaranteed as long as $\log p \lesssim n$. In other words, $p$ is allowed to grow exponentially with $n$, which is in line with the modern high-dimensional setting. Secondly, our Theorem 2 implies that, even if some conditions (such as MSG) for the exact recovery are not satisfied, one can still hope to partially recover the underlying permutation. In accordance to our theoretical result (\ref{rate}), our numerical results (Figure 4) show that, for the partial recovery, increasing $p$ indeed reduces the overall risk of the proposed estimator. Finally, as to the sample size $n$, we argue that, without assuming additional structural assumptions such as row-sparsity, it is very unlikely that including more samples will result in a worse estimate (see Table 1 and Figure \ref{simu.k.1} for numerical evidences).

\section{A GENERAL GROWTH MODEL}

In this section we study the permutation recovery over the general parameter space $\mathcal{D}$ where the growth pattern is not necessarily linear and therefore is more realistic inasmuch as the noisy nature of the shotgun metagenomic datasets \citep{boulund2018computational, gao2018quantifying}. The analysis relies on a deeper understanding of the relationship between the row-monotonic matrices and its leading singular vectors. 

Specifically, for any $\Theta\in \mathcal{D}$, we define the row-centered matrix
\beq \label{Theta'}
\Theta'  = \Theta(I-p^{-1}ee^\top) \in \R^{n\times p}
\eeq
whose singular value decomposition (SVD) is given by $\Theta'=\sum_{i=1}^r \lambda_i(\Theta') u'_i{v'}_i^\top$, with $r\le \min\{n,p\}.$
The following proposition is essential to our analysis of the general growth model.

\bep \label{rsv.prop}
Let $\Theta'$ be defined as above, then its first right singular vector $v'_1$ is a monotone vector, i.e., either $v'_{11}\le v'_{12}\le...\le v'_{1p}$ or  $v'_{11}\ge v'_{12}\ge...\ge v'_{1p}.$
\eep

Together with Proposition \ref{prop.inv}, the above proposition justifies our construction of the permutation estimator $\hat{\pi}$ using a PCA based approach. To overcome the identifiability issue, we further assume $\lambda_1(\Theta')$ has multiplicity one. We first   introduce the several quantities that play the key roles in permutation recovery over $\mathcal{D}$.

\bed 
For any $\Theta\in \mathcal{D}$ and the corresponding $\Theta'$ defined as above, we define 
\[
\Gamma=\min_{1\le i<j\le p}|{u'}_1^\top(\Theta'_{.i}-\Theta'_{.j})|=\lambda_1(\Theta')\min_{1\le i<j\le p}|v'_{1i}-v'_{1j}|,
\]
as the \textit{{local minimal signal gap}}, define 
\[
\Xi =\max_{1\le i\le p-1}\|\Theta'_{.i}-\Theta'_{.i+1}\|_2=\max_{1\le i\le p-1}\bigg(\sum_{j=1}^r \lambda^2_j(\Theta')|v'_{ji}-v'_{j,i+1}|^2\bigg)^{1/2},
\]
as the {\textit{local maximal signal gap}}, and define 
\[
\Lambda = \lambda_1^2(\Theta')-\lambda_2^2(\Theta')
\]
as the {\textit{global signal strength}} of $\Theta$.
\eed

In particular, the above definitions of $\Gamma$ and $\Lambda$ generalize the ones given earlier in the linear growth model as these quantities coincide for $\Theta\in\mathcal{D}_L$. The following theorem concerns the exact permutation recovery with $\hat{\pi}$ over $\mathcal{D}$.

\bet[Exact Recovery, General] \label{general.thm.1}
Suppose (A1) hold, $n\lesssim p$, $(\Theta,\pi)\in \mathcal{D}$, and $\Theta$ satisfies $\Gamma> C_0\sigma\sqrt{\log p}$ and
\[
\Lambda >  C_1\sigma^2\bigg[\bigg(n+\frac{\Xi^2}{\sigma^2}\bigg)\max\bigg\{\frac{(n+\log p)\sigma^2}{\Gamma^2},1\bigg\}+\sqrt{p}\bigg(\sqrt{n}+\frac{\Xi}{\sigma}\bigg)\max\bigg\{\frac{\sigma\sqrt{n+\log p}}{\Gamma},1  \bigg\}\bigg]
\]
for some $C_0 ,C_1> 0$. Then with probability at least $1-O(p^{-c})$ for some constant $c>0$, up to a permutation reversion, we have $\hat{\pi}=\pi$.
\eet

As in the case of linear growth model (Theorem \ref{linear.thm.1}), in Theorem \ref{general.thm.1}, to guarantee exact recovery, we need the MSG condition $\Gamma>C_0\sigma\sqrt{\log p}$. Unlike the linear growth model, here $\Gamma$ only implicitly depends on the elements of $\Theta$ through its spectral quantities, which makes its interpretation less clear. To address this issue, we make the following observation that links the minimal singular vector gap $\min_{1\le i<j\le p}|v'_{1i}-v'_{1j}|$ in the definition of $\Gamma$ to the elements of $\Theta$.

\bep \label{prop.suff}
Let $\Theta'$ in (\ref{Theta'}) be such that there exists a $\delta>0$ being the lower bound of the normalized minimum gap between any two entries in the same row, i.e. 
\[
\min_{1\le k\le n} \frac{|\theta'_{k,i}-\theta'_{k,j}|}{\|\Theta'_{k.}\|_2} \ge \delta \quad\text{for some $i\ne j$.}
\]
Then the first singular vector $v'_1\in \R^p$ of $\Theta'$ satisfies $|v'_{1,i}-v'_{1,j}| \ge \delta$.
\eep

Consequently, the implicit requirement that $\min_{1\le i<j\le p}|v'_{1i}-v'_{1j}|$ is large can be guaranteed when the normalized minimum distance $\min_{1\le i<j\le p}\min_{1\le k\le n}{|\theta'_{k,i}-\theta'_{k,j}|}/{\|\Theta'_{k.}\|_2}$ is large.
Our next theorem concerns the partial recovery over the general parameter space $\mathcal{D}$.

\bet[Partial Recovery, General] \label{general.thm.2}
Suppose (A1) hold, $n\lesssim p$, $(\Theta,\pi)\in \mathcal{D}$, and $\Theta$ satisfies
\begin{itemize}
	\item[(i)] there exits some $C_0>0$ such that $\max_{1\le j\le p-1}|v'_{1j}-v'_{1,j+1}|<C_0\min_{1\le j\le p-1}|v'_{1j}-v'_{1,j+1}|$ for all $p>0$, and
	\item[(ii)] $\Lambda >  C_1\sigma^2\big[\max\big\{\frac{(n+\log p)^2\sigma^2}{\Gamma^2},n+\frac{\Xi^2}{\sigma^2}\big\}+\sqrt{p}\max\big\{\frac{\sigma(n+\log p)}{\Gamma},\sqrt{n}+\frac{\Xi}{\sigma}  \big\}\big]$ for some $C_1>0$.
\end{itemize}
Then, up to a permutation reversion,
\begin{align*}
\E [\tau_K(\hat{\pi},\pi) ] &\le 1\land \bigg(\frac{c_0\sigma}{p\Gamma}\min\bigg\{1,e^{-\Gamma^2/2\sigma^2}\log\bigg(1+\frac{2\sigma^2}{\Gamma^2} \bigg)  \bigg\}+\frac{c_1e^{-\Gamma^2/2\sigma^2}}{p(\Gamma/\sigma+\sqrt{8/\pi})}+\frac{c_2}{p^{c+2}}\bigg)
\end{align*}
for some $c,c_0,c_1,c_2>0$. 
\eet
Condition (i) of Theorem \ref{general.thm.2} parallels the one given in Theorem \ref{linear.thm.2}. It essentially requires an even distancing of the elements (the projected columns of $\Theta$) whose ordering is to be tracked by $\hat{\pi}$. In contrast, in both Theorem \ref{general.thm.1} and \ref{general.thm.2}, the conditions on $\Lambda$ are slightly more complicated than those in Theorem \ref{linear.thm.1} and \ref{linear.thm.2}, as it further depends on the relative magnitude between $\Xi/\sigma$ and $\sqrt{n}$. In particular, if $\Xi/\sigma \lesssim \sqrt{n}$, the conditions reduce to the ones required in the linear growth models. Interestingly, the risk upper bound obtained in Theorem \ref{general.thm.2} remains the same as in the linear growth model, which only depends on $p$ and the signal-to-noise ratio $\Gamma/\sigma$.

\section{MINIMAX LOWER BOUNDS AND OPTIMALITY}\label{minimax}

In this section, we establish the minimax lower bounds for both exact and partial recovery considered in previous sections, in relation to different levels of the signal-to-noise ratio $\Gamma/\sigma$. In the following theorem, we show the MSG condition for exact recovery is asymptotically sharp.

\begin{figure}[h!]
	\centering
	\includegraphics[angle=0,width=13cm]{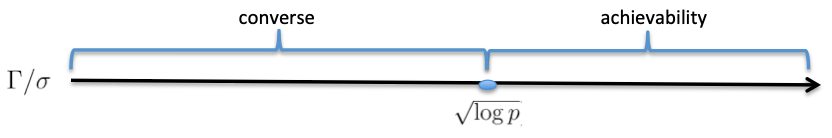}
	\caption{\footnotesize A graphical illustration of the achievability/converse result for exact recovery.} 
	\label{achiv}
\end{figure}

\bet \label{linear.low.bnd.thm.1}
Suppose (A1) hold. Let $ \mathcal{D}_1=\mathcal{D}_L \cap\{(\Theta,\pi):\Gamma \le \frac{\sigma}{4} \sqrt{\log p}\}$ and $ \mathcal{D}'_1=\mathcal{D} \cap\{(\Theta,\pi):\Gamma \le \frac{\sigma}{4} \sqrt{\log p}\}$. Then for any $p\ge 10$, we have
\[
\inf_{\hat{\pi}}\sup_{(\Theta,\pi)\in \mathcal{D}'_1}P(\hat{\pi}\ne \pi) \ge \inf_{\hat{\pi}}\sup_{(\Theta,\pi)\in \mathcal{D}_1}P(\hat{\pi}\ne \pi) \ge 0.3,
\]
where the infimum is over all the permutation estimators $\hat{\pi}$. 
\eet
This theorem along with Theorem \ref{linear.thm.1} and Theorem \ref{general.thm.1} indicates that our proposed estimator is minimax rate-optimal over $\mathcal{D}_L$ and $\mathcal{D}$ in terms of the MSG condition on $\Gamma$.
In light of Proposition \ref{prop.1}, in some situations the MSG condition can be both necessary and sufficient for the exact recovery, which includes practically important cases such as $n\asymp p$, $n< \log p$, etc. Using the information-theoretic language, we have therefore obtained both the achievability result, i.e., the existence of an algorithm or estimator that exactly recovers signal with high probability, and the converse result, namely, an upper
bound on the probability of exact recovery that applies to
any estimators \citep{cullina2016improved}. See Figure \ref{achiv} for an illustration.

Our next theorem establishes a minimax lower bound for the expected rate of convergence for the partial recovery.
\bet \label{linear.low.bnd.thm.2}
Suppose (A1) hold, $ \mathcal{D}_2(t)=\mathcal{D}_L \cap\{(\Theta,\pi):\Gamma \asymp t\}$ and $ \mathcal{D}'_2(t)=\mathcal{D} \cap\{(\Theta,\pi):\Gamma \asymp t\}$ and $t/\sigma \ge 2$. Then there exist constants $C_1,C_2>0$ such that
\[
\inf_{\hat{\pi}}\sup_{(\Theta,\pi)\in \mathcal{D}'_2(t)}\E [\tau_K(\hat{\pi},\pi) ] \ge \inf_{\hat{\pi}}\sup_{(\Theta,\pi)\in \mathcal{D}_2(t)}\E [\tau_K(\hat{\pi},\pi) ] \ge \frac{C_1\sigma }{pt}e^{-t^2/2\sigma^2}+\frac{C_2}{p^2}.
\]
\eet

Comparing the above minimax lower bound to the risk upper bounds obtained in Theorem \ref{linear.thm.2} and \ref{general.thm.2}, we conclude that our proposed estimator $\hat{\pi}$ is minimax rate-optimal in terms of the partial recovery for both the linear growth model and the general growth model over the range whenever $\Gamma/\sigma$ does not  diminish (Figure \ref{rec}). In particular, in Theorem 5 and 6, since the minimax lower bounds only concern the worst-case scenarios, the same lower bounds should hold for any parameter spaces whenever the same worst cases are included. Similarly, the assumption (A1) does not pose a restriction to the general applicability of such results.

\section{NUMERICAL STUDIES}

\subsection{Simulation with Model-Generated Data}

To demonstrate our theoretical results and compare with alternative methods, we generate data from model (\ref{m2}) with various configurations of the signal matrix $\Theta$. We compare the empirical performance of our proposed estimator $\hat{\pi}$ with the following alternatives:
\begin{itemize}
	\item $\pi_{\text{mean}}$ : Order the columns of $Y$ by the magnitude of its column means;
	\item $\pi_{\max}$: Order the columns of $Y$ by the magnitude of its column maximums.
\end{itemize}
We use both the 0-1 loss and the normalized Kendall's tau distance in comparing these methods. Due to the identifiability issue, the performance of each estimator is evaluated up to a complete reversion of the permutation. For example, we use $\min\{ \tau_K(\hat{\pi},\pi),\tau_K(\hat{\pi},\text{rev}(\pi)) \}$ as the empirical Kendall's tau distance. 
By symmetry, we set the underlying permutation $\pi=id$. 
The signal matrix $\Theta=(\theta_{ij})\in \R^{n\times p}$ is generated under the following four regimes:
\begin{itemize}
	\item[(i)] $S_1(\alpha,n,p)$: For any $1\le j\le p$, $\theta_{ij}=\log(1+j\alpha_i +\beta_i)$ where $\alpha_i\sim \text{Unif}(\alpha/2,\alpha)$ for $1\le i\le n/2$, $\alpha_i\sim \text{Unif}(0,0.01)$ for $n/2<i\le n$, and $\beta_i\sim \text{Unif}(1,3)$ for all $1\le i\le n$;
	\item[(ii)] $S_2(\alpha,n,p)$: For any $1\le j\le p$, $\theta_{ij}=j\alpha_i +\beta_i$ where $\alpha_i\sim \text{Unif}(\alpha/2,\alpha)$ for $1\le i\le n/2$, $\alpha_i\sim \text{Unif}(0,\alpha/10)$ for $n/2<i\le n$, and $\beta_i\sim \text{Unif}(1,3)$ for all $1\le i\le n$;
	\item[(iii)] $S_3(\alpha,n,p)$: For any $1\le j\le p$, $\theta_{ij}=\log(1+j\alpha_i +\beta_i)$ where $\alpha_i\sim \text{Unif}(\alpha/2,\alpha)$ for $1\le i\le 3$, $\alpha_i\sim \text{Unif}(0,0.01)$ for $4<i\le n$, and $\beta_i\sim \text{Unif}(1,3)$ for all $1\le i\le n$;
	\item[(iv)] $S_4(\alpha,n,p)$: For any $1\le j\le p$, $\theta_{ij}=j\alpha_i +\beta_i$ where $\alpha_i\sim \text{Unif}(\alpha/2,\alpha)$ for $1\le i\le 3$, $\alpha_i\sim \text{Unif}(0,\alpha/10)$ for $4<i\le n$, and $\beta_i\sim \text{Unif}(1,3)$ for all $1\le i\le n$.
\end{itemize}
Specifically, under each regime, the sample-specific ``growth rate" parameter $\alpha_i$ is randomly and uniformly generated either from the interval $[\alpha/2,\alpha]$ or an interval with much smaller values, namely, $[0,\alpha/10]$ in $S_2$ and $S_4$ and $[0,0.01]$ in $S_1$ and $S_3$. By construction, the four regimes consist of the nonlinear growth model where the signals spread out over many samples ($S_1$) or concentrate at a few rows ($S_3$) and the linear growth model where the signals spread out over many samples $(S_2)$ or concentrate at a few rows $(S_4)$. In particular, in accordance to our theory, for the supposedly ``non-informative" samples, we allow the corresponding growth rates to be small but non-zero, which shows the flexibility of our proposed method.
The entries of $Z$ are drawn from i.i.d. centred normal distributions whose variance $\sigma^2$ will be given explicitly.
In each setting, we evaluate the empirical performance of each method over a range of $n$, $p$ or $\alpha$. Each setting is repeated for 200 times. 

For  the exact recovery, in Table \ref{table:t2}, we reported the empirical risks of the estimators under the 0-1 loss for various regimes and parameter combinations. The noise level $\sigma^2$ is chosen for each regime to better illustrate the differences in the empirical risks among the estimators. From our simulation results, in consistent to our theory, our proposed estimator  has the smallest empirical risk over all the settings, and the estimation risk decreases as we increase $\alpha$, $n$ or decrease $p$.

\begin{table}[h!]
	\centering
	\caption{The empirical risks of the estimators under the 0-1 loss based on 200 simulations for various combinations of the parameters $(p,n,\alpha)$. 
		$\hat{\pi}$: proposed method; $\pi_{mean}$: mean-based method; $\pi_{max}$: max-based method.	}
	\begin{tabular}{lcccccccc}
		\hline
		\multirow{2}{2.95em}{$p=75$\\$n=40$}&\multicolumn{2}{c}{$S_1(\sigma^2=0.025)$} &\multicolumn{2}{c}{$S_2(\sigma^2=0.1)$}&\multicolumn{2}{c}{$S_3(\sigma^2=0.0075)$}&\multicolumn{2}{c}{$S_4(\sigma^2=0.025)$}  \\
		&   $\alpha=0.1$  &0.2 & 0.1 & 0.2 & 0.1  & 0.2 & 0.1 & 0.2  \\  
		\hline  
		$\hat{\pi}$ &    0.775  & 0.575 &0.415 & 0.000 & 0.025&0.020&0.025&0.000  \\
		$\pi_{mean}$ & 0.925 & 0.815 & 0.955 & 0.015   &0.155&0.135&0.880&0.005   \\
		$\pi_{max}$ &  1.000 & 1.000 & 1.000  & 0.995    &0.995&0.970&0.840&0.430 \\
		\hline
		
		\multirow{2}{3.3em}{$n=40$\\$\alpha=0.1$}&\multicolumn{2}{c}{$S_1(\sigma^2=0.025)$} &\multicolumn{2}{c}{$S_2(\sigma^2=0.1)$}&\multicolumn{2}{c}{$S_3(\sigma^2=0.0075)$}&\multicolumn{2}{c}{$S_4(\sigma^2=0.025)$}  \\
		&   $p=60$  &90 & 60 & 90 & 60  & 90 & 60 & 90  \\  
		\hline  
		$\hat{\pi}$ &    0.410  & 0.930 &0.340 & 0.470 & 0.010&0.115&0.000&0.010  \\
		$\pi_{mean}$ & 0.720 & 0.985 & 0.910 & 0.980   &0.070&0.245&0.775&0.900   \\
		$\pi_{max}$ &  1.000 & 1.000 & 1.000  & 1.000    &0.975&1.000&0.815&0.875 \\
		\hline
		
		\multirow{2}{3.3em}{$p=75$\\$\alpha=0.1$}&\multicolumn{2}{c}{$S_1(\sigma^2=0.025)$} &\multicolumn{2}{c}{$S_2(\sigma^2=0.1)$}&\multicolumn{2}{c}{$S_3(\sigma^2=0.0075)$}&\multicolumn{2}{c}{$S_4(\sigma^2=0.025)$}  \\
		&   $n=40$  &60 & 40 & 60 & 40  & 60 & 40 & 60  \\  
		\hline  
		$\hat{\pi}$ &    0.765  & 0.440 &0.475 & 0.095 & 0.050&0.020&0.010&0.005  \\
		$\pi_{mean}$ & 0.920 & 0.645 & 0.940 & 0.700   &0.175&0.045&0.900&0.905   \\
		$\pi_{max}$ &  1.000 & 1.000 & 1.000  & 1.000    &0.995&0.995&0.855&0.820 \\
		\hline
	\end{tabular}
	\label{table:t2}
\end{table}

\begin{figure}[h!]
	\centering
	\includegraphics[height=0.15\textheight,angle=0,width=13cm]{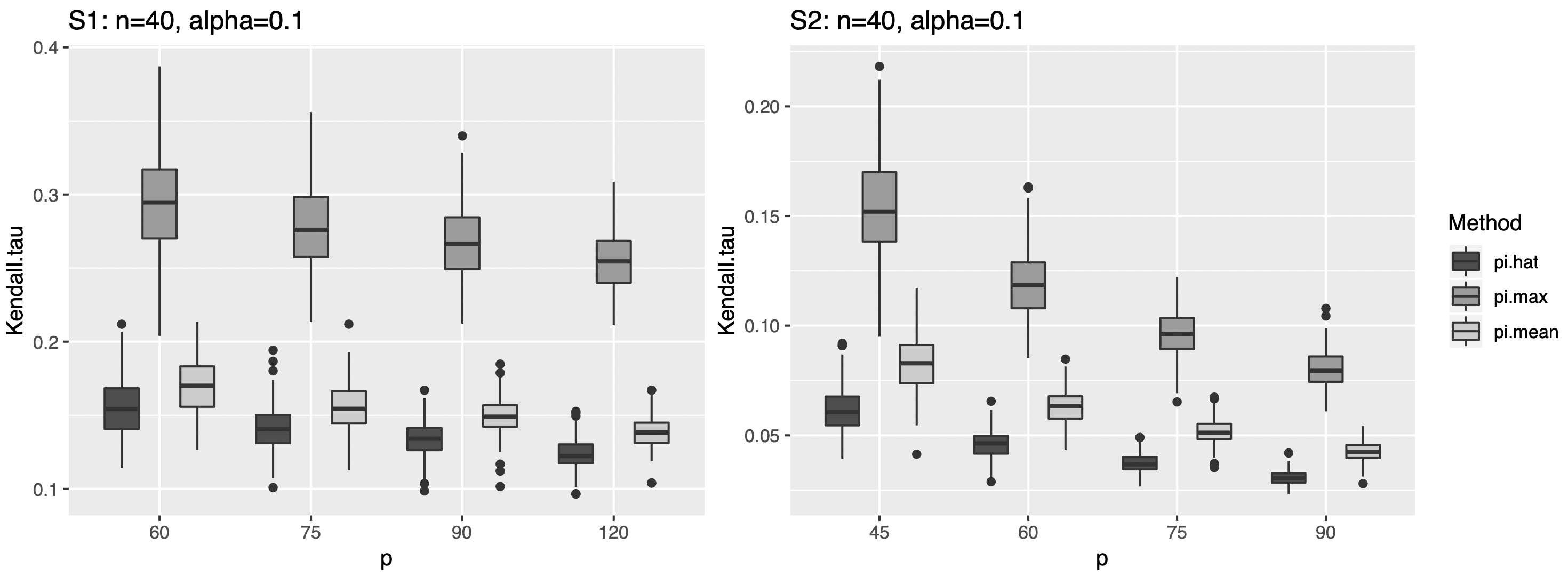}
	\includegraphics[height=0.15\textheight,angle=0,width=13cm]{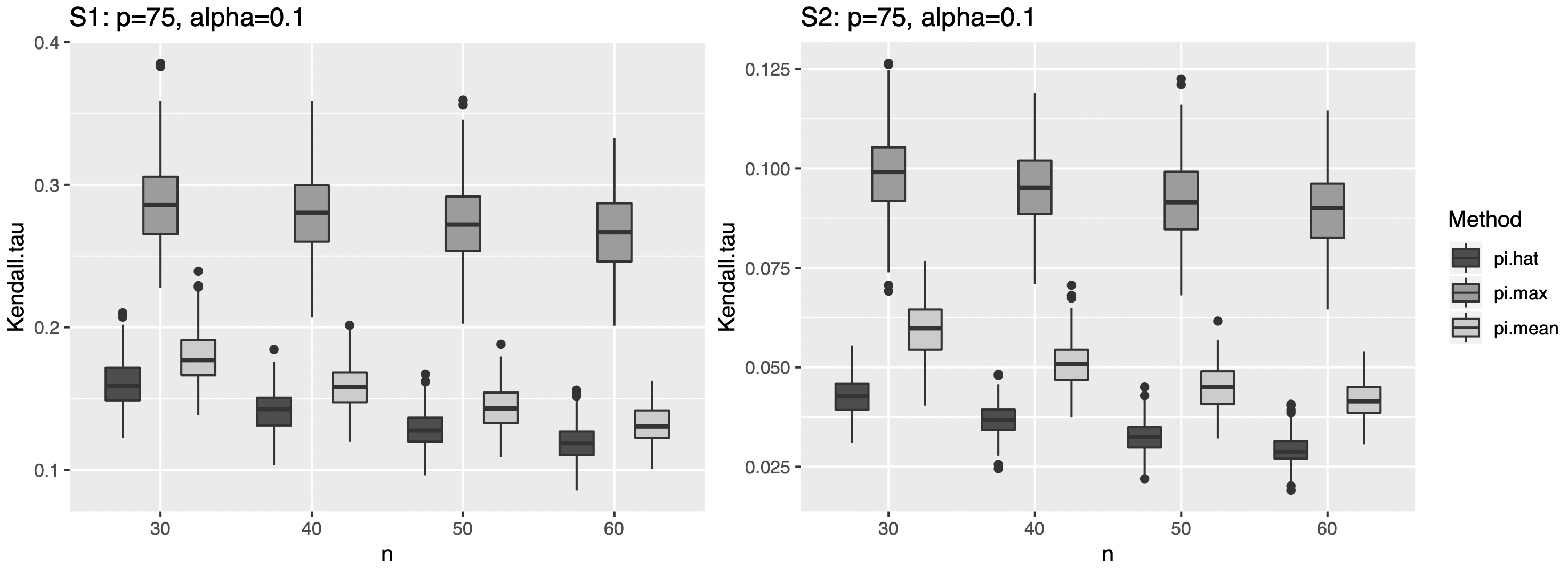}
	\includegraphics[height=0.15\textheight,angle=0,width=13cm]{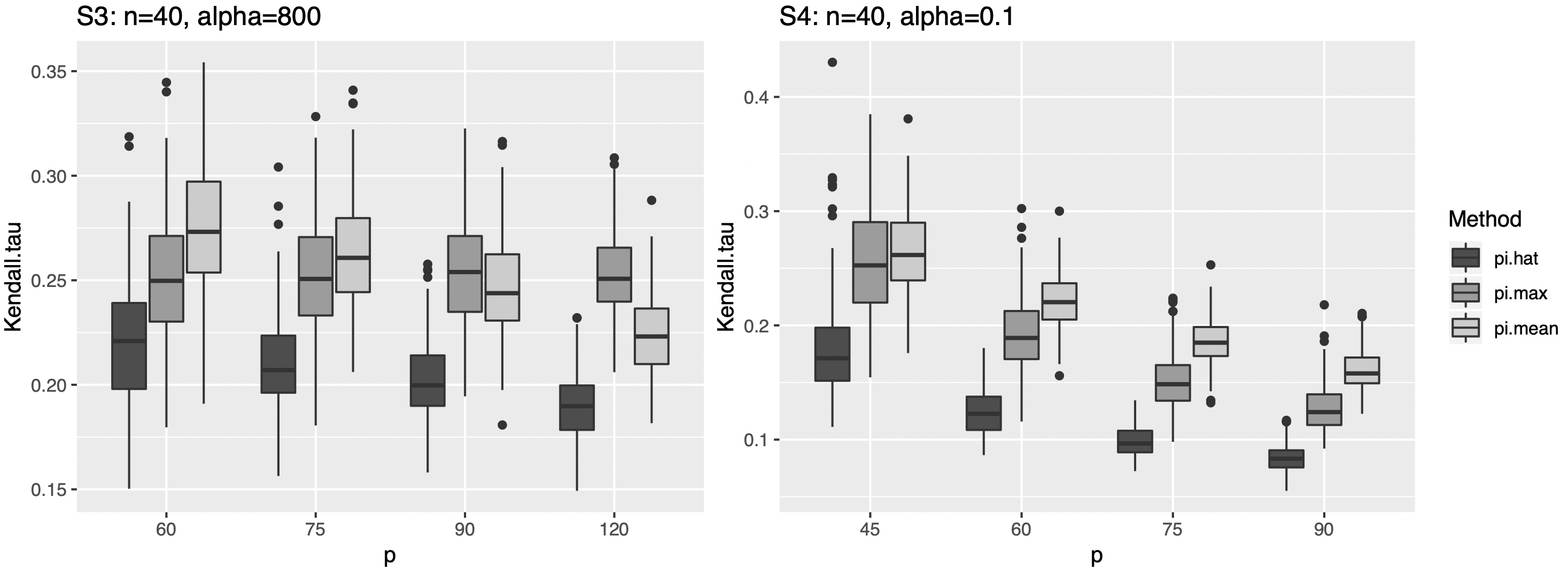}
	\includegraphics[height=0.15\textheight,angle=0,width=13cm]{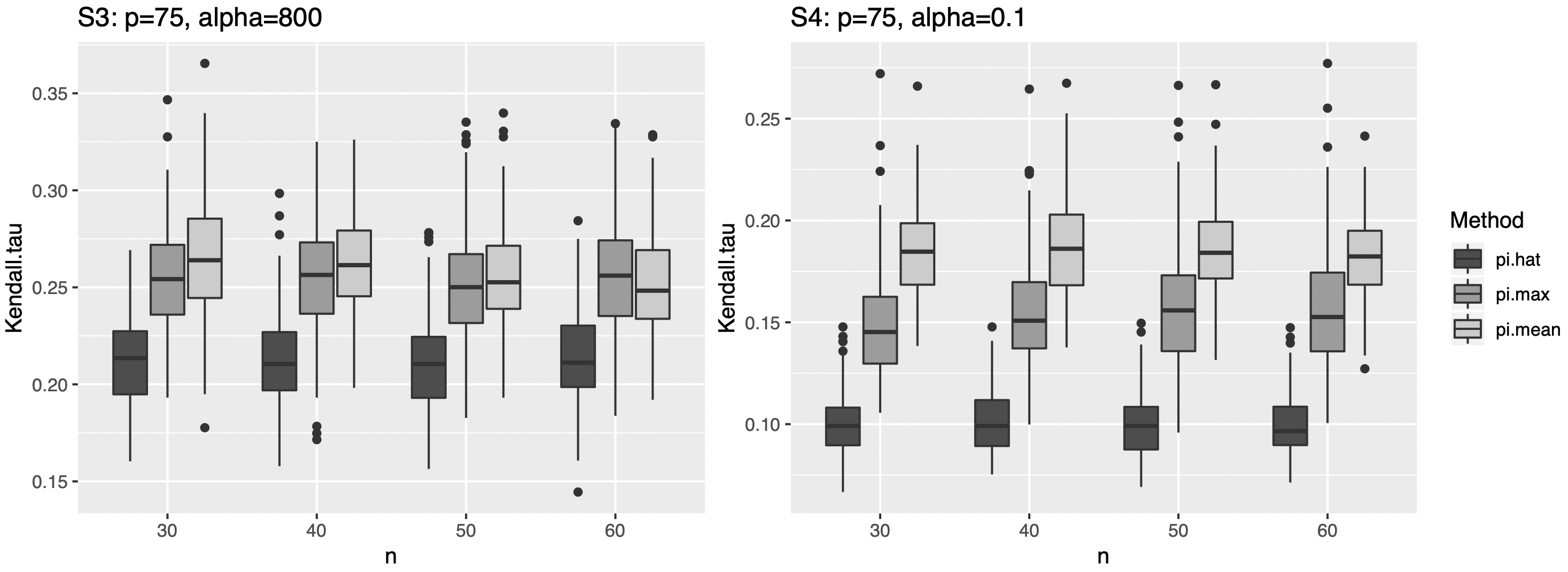}
	\caption{\footnotesize Boxplots of the empirical normalized Kendall's distance between the estimated and true permutations under different models. $\hat{\pi}$: proposed estimator; $\pi_{\text{mean}}$: mean-based estimator; 
		$\pi_{\max}$: max-based estimatior.} 
	\label{simu.k.1}
\end{figure}

For partial recovery, in Figure \ref{simu.k.1}, we show boxplots of the empirical normalized Kendall's tau between each estimator and the true permutation $\pi$. Again, our proposed method outperforms the alternatives in all the cases. As expected from our theory, under all the four regimes, increasing $p$ while keeping other parameters fixed results to smaller estimation risk. As for the dependence on $n$, under $S_1$ and $S_2$, increasing $n$ leads to smaller risk as it is equivalent to increasing $\Gamma$, whereas under $S_3$ and $S_4$, the risk roughly remains the same across different $n$'s as in these case $\Gamma$ doesn't change much.

To offer more intuitive interpretation of why $\hat{\pi}$ performs better than the alternative methods, we assessed the weight vectors $\hat{w}$ of our proposed estimator $\hat{\pi}$ under each regime after 200 rounds of simulations (Figure 3 in Supplemented Material). In comparison, the weight vector for $\pi_{mean}$ is simply $(1/\sqrt{n},...,1/\sqrt{n})$, which assigns equal weight to all the samples. On the other hand, since $\pi_{\max}$ cannot be written in the form of $(\frak{r}(w^\top Y))^{-1}$ for some weight vector $w$ and therefore does not belong to the class of linear projection estimators, we reported instead the pseudo-weight vector $\tilde{w}\in \R^n$ where the $i$-th component is the proportion that the $i$-th sample is used among the $p$ coordinates. In general, we found that $\tilde{w}\in \R^n$ assigns larger weights to only a few samples among those with higher signal strength, and the weight vector for $\pi_{mean}$ fails to distinguish the informative samples from the non-informative ones. In contrast, the weight vectors $\hat{w}$ for our proposed estimator $\hat{\pi}$ would automatically adapt to the varying signal strengths across the samples and assign larger weights to the samples with more significant signal changes. This also explains the interesting phenomenon in Figure  \ref{simu.k.1} that, under the regime $S_1$ and $S_2$, $\hat{\pi}$ and $\pi_{\text{mean}}$ perform better than $\pi_{\max}$, whereas under $S_3$ and $S_4$, $\hat{\pi}$ and $\pi_{\max}$ perform better. In summary, methods that are able to detect and assign larger weight to the more informative samples would perform better than methods that are not. Observably, $\hat{\pi}$ combines the advantages of $\pi_{\text{mean}}$ and $\pi_{\max}$ in that it finds the best weights (projection scores) in a data-driven manner.

\subsection{Evaluation Using  Synthetic Metagenomic Data}

We evaluate the empirical performance of our proposed method using a synthetic metagenomic sequencing dataset  used  in \cite{gao2018quantifying} by generating sequencing reads based on 45 bacterial genomes.  Instead of estimating the PTRs, which was the focus of \cite{gao2018quantifying}, our goal is to recover the unknown relative orders of the contigs assembled in typical metagenomics studies. In addition to assisting the estimation of PTRs, such ordering of the contigs  could be of independent interest for other applications, including genome assemblies based on shotgun metagenomics data.

\cite{gao2018quantifying} presented a synthetic shotgun metagenomic sequencing dataset  of a community of 45  phylogenetically related species from 15 genera of five different phyla with known RefSeq ID, taxonomy and replication origin \citep{synth} (see Figure 2 in our Supplementary Material). To generate metagenomics reads, reference genome sequences of randomly selected three species in each genus were downloaded from NCBI.  Read coverages were generated along the genome based on an exponential distribution with a specified peak-to-trough ratio and a function of accumulative distribution of read coverages along the genome was calculated. Sequencing reads were next generated using the above accumulative distribution function and a random  location of each read on the genome, until the total read number achieved a randomly assigned average coverage between 0.5 and 10 folds for the species in a sample. Sequencing errors including substitution, insertion and deletion were simulated in a position- and nucleotide-specific pattern according to a recent study on metagenomic sequencing error profiles of Illumina. 

For the final dataset, the average nucleotide identities  (ANI) between species within each genus ranged from 66.6\% to 91.2\%  The probability of one species existing in each of the 50 simulated samples was set as 0.6, and a total of 1,336 average coverages and the corresponding PTRs were randomly and independently assigned.  After the same processing and filtering  steps and CG-adjustment step as  in \cite{gao2018quantifying}, the final dataset included genome assemblies of 41 species. For each species, we obtained 
the permuted matrix of log-contig counts, with the number of samples ranging from 29 to 46, and the number of contigs ranging from 47 to 482. 


Our proposed method ($\hat{\pi}$) was used to estimate the unknown orders of the contigs for each species and each sample. As a comparison, we also considered the iRep estimator proposed in \cite{brown2016measurement}, where the contigs  of a given species were ordered for each sample separately based on the read counts observed. 
We evaluate these methods by comparing the estimated contig orders to their true orders as measured by the normalized Kendall's tau distance. To generalize our evaluation to diverse metagenomic datasets, we also evaluate the effect of sample size as well as contig numbers by randomly selecting subsets of samples or contigs from each dataset. The selection was made with replacement. 

\begin{figure}[h!]
	\centering
	\includegraphics[angle=0,width=16cm]{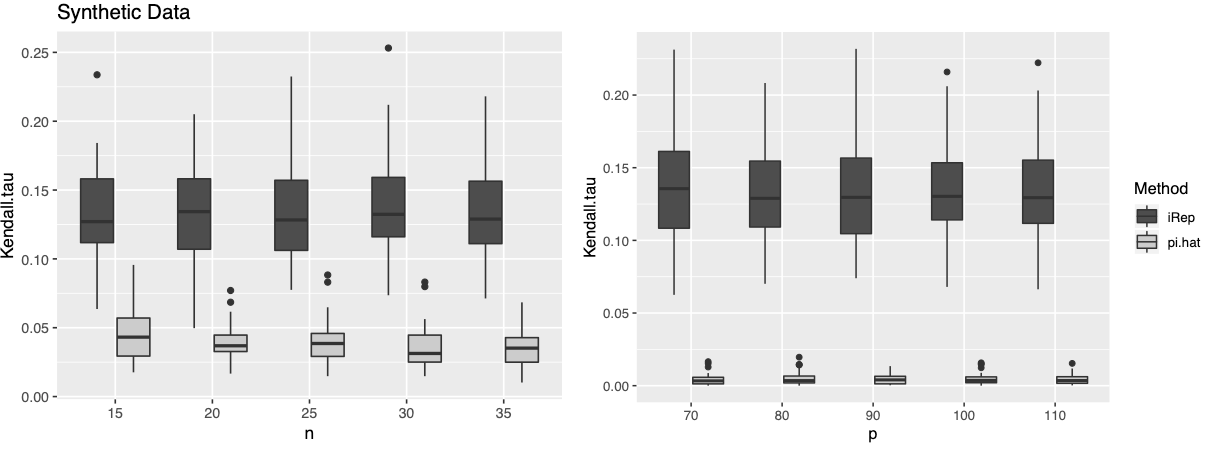}
	\caption{Boxplots of the normalized Kendall's distance between the estimated contig orders and the true orders for different sample sizes $n$ and different numbers of contigs $p$. The lighter ones correspond to our proposed method and the darker ones correspond to the iRep estimation method.}
	\label{synth.data}
\end{figure}

The results are summarized in Figure \ref{synth.data} by comparing the normalized Kendall's tau distances.  As $n$ or $p$ varies, our proposed estimator performs consistently better than iRep in recovering the true contig orders, which explains partially  why the DEMIC algorithm worked better in estimating the bacterial growth dynamics. The results of our methods are not sensitive to the sample size and the number of contigs from the genome assemblies.  Our estimator also shows smaller variability.

\subsection{Analysis of a Real Microbiome Metagenomic Data Set}

Finally, we complete our numerical studies by analyzing a real metagenomic dataset from the Pediatric Longitudinal Study of Elemental Diet and Stool Microbiome Composition (PLEASE) study, a prospective cohort study
to investigate the treatment effects on the gut microbiome and reduction of inflammation  in pediatric Crohn's disease patients \citep{lewis2015inflammation}. In particular, sequencing data from the fecal samples of 86 Crohn's disease children were obtained at baseline, 1 week and 8 weeks after antiTNF or enteral diet treatment. In our analysis, the sequencing data at the 8th week after treatment was used to compare the bacterial growth dynamics for non-responders ($n=34$) and responders ($n=47$). The reads were downloaded from NCBI short read archive (SRP057027) with the corresponding metadata.
After the same coassembly, alignment and binning steps as in \cite{gao2018quantifying}, the DEMIC algorithm was applied to estimate the bacterial growth rate of a given species represented by a contig cluster (bin) for each sample. In particular, DEMIC applied our proposed method to the GC-adjusted contig coverage data to recover the original order of the contigs. 
After obtaining the ordered contigs, a simple linear regression was fitted to obtain estimates of the PTRs (ePTRs). 

In order to compare the baterial growth rates between responders and non-responders,   our analysis focused on ePTRs of 8 contig clusters  over subsets of the non-responders ($n_1$) and the responders ($n_2$) after 8 weeks of treatment with $\min\{n_1,n_2\}>5$. Other  contig clusters were rare and only appeared in a few samples.  For each contig cluster, we compare the ePTRs of the responders and non-responders  Wilcoxon rank sum test (Table S.1 in Supplementary Material). The taxonomic annotations of these eight  contig clusters were  obtain by applying  the BAT algorithm \citep{von2019robust} that compares the metagenomic assembled bins to a taxonomy database. In Table S.1, we show the final taxonomic annotations for each bin to the finest possible resolution, with the lineage scores indicating the quality of each taxonomic classification.

Among the 8 contig clusters, bin.026  showed a significant difference in ePTRs between responders and non-responders after either antiTNF or enteral diet treatment for 8 weeks ($p$=0.0418), where the growth rate was higher in    Crohn's disease patients who did not respond to the treatment. The  taxonomic classification (Table S.1) shows that this  contig cluster belongs to the phylum {\it Firmicutes} and the order {\it Clostridiales}.   Since BAT algorithm was not able to classify the  order {\it Clostridiales} to finer taxonomic level of known species, this contig cluster may represent a new species that is  important to the treatment outcome of  Crohn's disease patients.

\section{DISCUSSION}

In this paper, partial recovery was studied under the normalized Kendall's tau distance. Another commonly used metric is the normalized Spearman's footrule distance defined by
\[
\rho(\pi_1,\pi_2)= \frac{2}{p(p-1)}\sum_{i=1}^p |\pi_1(i)-\pi_2(i)|,\quad\quad \pi_1,\pi_2\in\mathcal{S}_p.
\]
A celebrated result by \cite{diaconis1977spearman} shows that $\tau_K(\pi_1,\pi_2)\le \rho(\pi_1,\pi_2)\le 2\tau_K(\pi_1,\pi_2),$
which means  the two distances are equivalent.
As a consequence, all the theoretical results presented in this paper concerning the Kendall's tau distance  also hold for the Spearman's footrule distance without any change.

The minimax optimality of the proposed estimator $\hat{\pi}$ was investigated in Section \ref{minimax} by examining the asymptotic sharpness of the MSG condition for exact recovery, and by obtaining the matching minimax risk lower bound for partial recovery. There are a few issues that deserve further investigation.  For both exact and partial recovery, it is unclear to what extent the GSS condition is necessary. In our risk analysis, the perturbation bound for the left singular subspace \citep{cai2018rate} was used. In fact, similar results can be obtained using the concentration bound for the linear functionals of singular vectors \citep{koltchinskii2016perturbation}. Nevertheless, it remains to show whether the GSS condition is also asymptotically sharp.
In addition, in Theorem \ref{linear.low.bnd.thm.2}, the matching minimax lower bound was obtained only for nonvanishing $\Gamma/\sigma$. It remains to show whether the rate $\sigma/(p\Gamma)$ is minimax optimal when $\Gamma/\sigma\to 0$. The difficulty lies in finding  a $p^{1+\delta}$-sphere packing of the group $\mathcal{S}_p$ equipped with the Kendall's tau distance for any $0<\delta<1$, while the pairwise $\ell_2$ distances of the packing elements are also well controlled. Some initial steps have been made in the so-called rank modulation theory \citep{barg2010codes,mazumdar2013constructions}. 

There are several related problems that are also of significant theoretical and practical interest. Firstly, although we used the Kendall's tau distance or the equivalent Spearman's footrule distance as the metric for partial recovery, other distances such the Hamming distance, Spearman's rank correlation distance, and Ulam's distance have also been used as the performance metrics for partial recovery in other permutation estimation problems \citep{golouglu2015new,mukherjee2016estimation}. It is therefore of interest to see how $\hat{\pi}$ performs under these losses. Secondly, our proposed estimator $\hat{\pi}$ implicitly performs a (linear) dimension reduction technique and only uses the information contained in the first eigenvector of $A$ in (\ref{A}). A natural extension is to consider the eigen-subspace spanned by the first $k$ eigenvectors and to estimate the permutation in a sequential manner.

The present paper focuses on the estimation of the permutation matrix $\Pi$. It is also of interest to estimate the underlying signal matrix $\Theta$ or some functionals of it. For example, in microbiome growth dynamics studies, it is of significant interest to estimate the peak-to-trough ratio $\exp(\theta_{kp}-\theta_{k1})$ for $k=1,...,n$, which measures the microbial growth rate for the $k$th sample, and to identify the samples with peak-to-trough ratio of 1.  It is also interesting to identify the bacteria that show differential growth dynamics between disease and normal individuals.  Finally,  robust  permutation recovery methods that can relax the Gaussian or sub-Gaussian assumption of the noise in the permuted monotone matrix model are needed. For example, in some applications, the columns of the noise matrix are not independent, or the variance levels across the noise matrix are not identical. In these cases, we argue that, as long as the marginal distributions of the noise matrix entries remain sub-Gaussian, the analytical framework of the current paper can still be applied, but with more efforts to control the underlying heteroskedasticity. Toward this end, results from the recent work of \cite{zhang2018heteroskedastic} can be very useful, in terms of the new technical tools that parallel the ones used in the current paper to analyse the homoskedastic PCA (cf. Lemma 2 and 3).  Finally, to account for non-informative samples, sparse PCA \citep{CMW2013SPCA,YZ2013} can be considered. These are interesting problems left for future research.

\section{PROOFS OF THE MAIN THEOREMS}

In this section, we prove Theorems \ref{linear.thm.1} and \ref{linear.thm.2} in detail and briefly sketch the proofs of Theorems \ref{general.thm.1} and \ref{general.thm.2}. We also prove the minimax lower bounds in Theorems \ref{linear.low.bnd.thm.1} and  \ref{linear.low.bnd.thm.2}. Proofs of other results including the technical lemmas can be found in the online Supplementary Materials.

\paragraph{Proof of Theorem \ref{linear.thm.1}.} Let $X=\Theta+Z$. It follows that $Y=X\Pi$. By right invariance of the 0-1 loss with respect to permutation composition, we have
\[
\ell((\frak{r}(\hat{w}^\top Y))^{-1},\pi) = \ell((\frak{r}(\hat{w}^\top X\Pi))^{-1},\pi) = \ell((\frak{r}(\hat{w}^\top X))^{-1}\circ \pi,\pi) = \ell((\frak{r}(\hat{w}^\top X))^{-1},id).
\]
Thus it suffices to study the risk $\E\ell((\frak{r}(\hat{w}^\top X))^{-1},id)=P((\frak{r}(\hat{w}^\top X))^{-1}\ne id)$. In fact,
\begin{align} \label{risk.upper.bnd}
&P((\frak{r}(\hat{w}^\top X))^{-1}\ne id) \le P\bigg(\bigcup_{i=1}^{p-1}\bigg\{ \sum_{k=1}^n\hat{w}_kX_{ki}\ge  \sum_{k=1}^n\hat{w}_kX_{k,i+1} \bigg\}\bigg) \nonumber  \\
&\le \sum_{i=1}^{p-1}P\bigg(  \sum_{k=1}^n\hat{w}_kX_{ki}\ge  \sum_{k=1}^n\hat{w}_kX_{k,i+1}  \bigg)=\sum_{i=1}^{p-1}P\bigg(  \sum_{k=1}^n\hat{w}_k(X_{ki}-X_{k,i+1})\ge  0 \bigg),
\end{align}
which further reduces to obtaining an upper bound for $P_i=P\big(  \sum_{k=1}^n\hat{w}_k(X_{ki}-X_{k,i+1})\ge  0 \big)$.
By definition, $\hat{w}$ is the first eigenvector of $A = Y\big(I-\frac{1}{p}ee^\top\big)\big(I-\frac{1}{p}ee^\top\big)  Y^\top$.
Simple calculation yields $\Pi\big(I-\frac{1}{p}ee^\top\big)\big(I-\frac{1}{p}ee^\top\big) \Pi^\top = \big(I-\frac{1}{p}ee^\top\big)\big(I-\frac{1}{p}ee^\top\big)$ for any $\Pi\in \mathcal{S}_p$. So $\hat{w}$ is also the first eigenvector of $A=X\big(I-\frac{1}{p}ee^\top\big)\big(I-\frac{1}{p}ee^\top\big)  X^\top \equiv TT^\top$, where $T\in \R^{n\times p}$. Note that $T$ admits the decomposition $T = \Theta'+E \in \R^{n\times p}$ where $E_{ij}\sim N(0, (p-1)\sigma^2/p)$ and $\Theta'=a\eta'^\top$,$\eta'_j=\eta_j-\frac{1}{p}\sum_{i=1}^p\eta_i$.
In particular, $T_{.i}=X_{.i}-\bar{X}_{row}$ where $\bar{X}_{row}=p^{-1}\sum_{i=1}^pX_{.i} \in \R^n$ is the vector of row means of $X$. We denote $\phi_{ij}=T_{.i}-T_{.j}=X_{.i}-X_{.j}$ and denote $w = a/\|a\|_2 \in \R^n$ as the first eigenvector of the rank-one matrix $\Theta'\Theta'^\top$. 
Now following (\ref{risk.upper.bnd}), we have
\begin{align*}
&P_i=P\bigg( \sum_{k=1}^n \hat{w}_k (X_{ki}-X_{k,i+1})\ge 0 \bigg) = P\bigg( w^\top \phi_{ij}+(\hat{w}-w)^\top \phi_{ij}\ge 0 \bigg) \\
&= P\bigg( w^\top \phi_{i,i+1}+(\hat{w}-w)^\top \phi_{ij}\ge 0, |1-(\hat{w}^\top w)^2|\le\delta \bigg) \\
&\quad+P\bigg( w^\top \phi_{i,i+1}+(\hat{w}-w)^\top \phi_{ij}\ge 0,  |1-(\hat{w}^\top w)^2|>\delta  \bigg)
\end{align*}
for some $\delta>0$.
By definition, up to a change of sign for $\hat{w}$, we have $0\le \hat{w}^\top w\le 1$. Then $ |1-(\hat{w}^\top w)^2|\le\delta$ implies $|(\hat{w}-w)^\top \phi_{i,i+1}| \le \|\hat{w}-w\|_2\|\phi_{i,i+1}\|_2\le \sqrt{2\delta}\|\phi_{i,i+1}\|_2$,
where the first inequality follows from Cauchy-Schwartz and the second inequality used $\|\hat{w}-w\|_2 = \sqrt{2(1-\hat{w}^\top w)}\le \sqrt{2(1-(\hat{w}^\top w)^2)}.$
Thus
\begin{align} \label{eq5.2}
P_i \le P\bigg( w^\top \phi_{i,i+1}\ge- \sqrt{2\delta}\|\phi_{i,i+1}\|_2 \bigg) +P\bigg( |1-(\hat{w}^\top w)^2|>\delta  \bigg).
\end{align}
The following lemmas provide upper bounds for the two probability events in the last expression.
\bel \label{lem.term1}
Under the conditions of Theorem \ref{linear.thm.1}, denote $\Gamma_i=\|a\|_2(\eta_i-\eta_{i+1})$, then for any $\delta>0$, we have
\begin{align}\label{term1.eq}
&P\bigg( w^\top \phi_{i,i+1}\ge- \sqrt{2\delta}\|\phi_{i,i+1}\| \bigg) \le\Phi\bigg( C_1\sqrt{\delta}\Psi_i^{1/2}+\frac{\Gamma_i}{\sigma} \bigg)+\frac{C_2}{p^c}
\end{align}
for $\Psi_i=(\sqrt{n}+ \sqrt{\log p})^2+\frac{\Gamma_i^2}{\sigma^2} +\frac{|\Gamma_i|}{\sigma}\sqrt{\log p}$ and some constants $C_1,C_2,c>0$.
\eel

\bel \label{conc.lem}
Suppose $\lambda_1^2(\Theta') \ge C\sigma^2(n+\sqrt{np})$ for some $C>0$, it follows that
\[
P\bigg( |1-(\hat{w}^\top w)^2|\le C_1\frac{\sigma^2(\lambda_1^2(\Theta')+\sigma^2p)(n+\log p)}{\lambda_1^4(\Theta')}\bigg)\ge 1-\frac{C_2}{p^c}
\]
for some $C_1,C_2,c>0.$
\eel

Now since $\frac{1}{p}\sum_{1\le i<j\le p}(\eta_i-\eta_j)^2=\sum_{j=1}^p\big(\eta_j-\frac{1}{p}\sum_{i=1}^p \eta_i \big)^2=\sum_{j=1}^p {\eta'_j}^2,$
we have $\lambda_1^2(\Theta')=\Lambda>C_0\sigma^2(n+\sqrt{np})$ for some $C_0>0$.
Set $\delta = C_0\sigma^2\frac{(\lambda_1^2(\Theta')+\sigma^2p)(n+\log p)}{\lambda_1^4(\Theta')}.$
It follows that $\delta=o(1)$.
Combining Lemma \ref{lem.term1} and Lemma \ref{conc.lem}, we have
\begin{align} \label{bd.eq1}
&P_i\le \Phi\bigg( {C\sqrt{\delta}}\bigg[ (\sqrt{n}+\sqrt{\log p})^2+\frac{\Gamma_i^2}{\sigma^2} +\frac{|\Gamma_i|}{\sigma}\sqrt{\log p} \bigg]^{1/2}+\frac{\Gamma_i}{\sigma} \bigg)+\frac{C}{p^c}
\end{align}
for some $C,c>0$.
The rest of the analysis is divided into several cases.
\paragraph{Case 1. $\log p \lesssim n$.} In this case, we have $P_i \le \Phi\big( {C\sqrt{\delta}}\big[ n+\frac{\Gamma_i^2}{\sigma^2} +\frac{|\Gamma_i|}{\sigma}\sqrt{\log p} \big]^{1/2}+\frac{\Gamma_i}{\sigma} \big)+\frac{C}{p^c}$.
In addition, if $|\Gamma_i|/\sigma \lesssim \sqrt{n}$, we have $P_i \le \Phi\big( C\sqrt{\delta n}+\frac{\Gamma_i}{\sigma} \big)+\frac{C}{p^c}\le \frac{C'}{p^c},$
where the last inequality follows from $\sqrt{\log p}\lesssim \Gamma/\sigma \le |\Gamma_i|/\sigma \lesssim \sqrt{n}$ and $\Lambda \gtrsim \sigma^2n\big(\frac{\sigma^2n}{\Gamma^2}+\frac{\sigma\sqrt{p}}{\Gamma} \big)$. If instead $|\Gamma_i|/\sigma \gtrsim \sqrt{n}$, we have $P_i \le \Phi\big( {C\sqrt{\delta}}\frac{|\Gamma_i|}{\sigma} +\frac{\Gamma_i}{\sigma} \big)+\frac{C}{p^c}\le \frac{C'}{p^c}$,
where the last inequality follows from $|\Gamma_i|/\sigma \gtrsim \sqrt{n}\gtrsim \sqrt{\log p}$ and $\delta =o(1)$.
Hence, in Case 1, (\ref{bd.eq1}) can be bounded by $O(p^{-c})$.

\paragraph{Case 2. $\log p \gtrsim n$.} In this case, we have $P_i \le \Phi\big( {C\sqrt{\delta}}\big[ \log p+\frac{\Gamma_i^2}{\sigma^2} +\frac{|\Gamma_i|}{\sigma}\sqrt{\log p} \big]^{1/2}+\frac{\Gamma_i}{\sigma} \big)+\frac{C}{p^c}$.
In addition, since $|\Gamma_i| \ge \Gamma \gtrsim \sigma\sqrt{\log p}$ and $\delta=o(1)$, we have $P\big( \sum_{k=1}^n \hat{u}_k (X_{ki}-X_{k,i+1})\ge 0 \big) \le \Phi\big( \frac{C\sqrt{\delta}}{\sigma}|\Gamma_i|+\frac{\Gamma_i}{\sigma} \big)+\frac{C}{p^c}\le \frac{C'}{p^c}$.
This shows that, in Case 2, (\ref{bd.eq1}) can also be bounded by $O(p^{-c})$. 

As a result, it follows that, up to a change of sign for $\hat{w}$, $P((\frak{r}(\hat{w}^\top X))^{-1}\ne id) =O(p^{-c})$ for some constant $c>0$.
\qed

\paragraph{Proof of Theorem \ref{linear.thm.2}.} Firstly, by invariance property of Kendall's tau distance, $\E [\tau_K(\hat{\pi},\pi) ]=\E [\tau_K((\frak{r}(\hat{w}^\top X))^{-1},id) ]=\E [\tau_K((\frak{r}(\hat{w}^\top X)),id) ].$
It then follows
\begin{align*}
\E [\tau_K(\hat{\pi},\pi) ]&= \frac{2}{p(p-1)}\sum_{i< j} P( [\frak{r}(\hat{w}^\top X)]_i \ge  [\frak{r}(\hat{w}^\top X)]_j )\\
&= \frac{2}{p(p-1)}\sum_{i< j} P\bigg( \sum_{k=1}^n\hat{w}_k(X_{ki}-X_{kj})\ge 0 \bigg).
\end{align*}
The summation in the last expression can be divided into two parts, namely, the consecutive differences and non-consecutive differences, i.e.,
\begin{align*}
\sum_{i< j} P\bigg( \sum_{k=1}^n\hat{w}_k(X_{ki}-X_{kj})\ge 0 \bigg)
&= \sum_{(i,j): j=i+1} P\bigg(\sum_{k=1}^n\hat{w}_k(X_{ki}-X_{kj})\ge 0 \bigg)\\
&\quad+\sum_{(i,j): j> i+1}P\bigg( \sum_{k=1}^n\hat{w}_k(X_{ki}-X_{kj})\ge 0 \bigg).
\end{align*}
In the following, we first show
\beq \label{eq.1}
P_i=P\bigg(\sum_{k=1}^n\hat{w}_k(X_{ki}-X_{k,i+1})\ge 0 \bigg) \le \frac{ce^{-\Gamma^2/2\sigma^2}}{\Gamma/\sigma+\sqrt{\Gamma^2/\sigma^2+8/\pi}}+\frac{C}{p^c}
\eeq
so that
\beq \label{eq17}
\sum_{(i,j): j=i+1} P\bigg(\sum_{k=1}^n\hat{w}_k(X_{ki}-X_{kj})\ge 0 \bigg) \le \frac{cpe^{-\Gamma^2/2\sigma^2}}{\Gamma/\sigma+\sqrt{\Gamma^2/\sigma^2+8/\pi}}+\frac{C}{p^c}.
\eeq
for some $C,c>0$. 
Then we show that 
\beq \label{eq.2}
\sum_{(i,j): j> i+1}P\bigg( \sum_{k=1}^n\hat{w}_k(X_{ki}-X_{kj})\ge 0 \bigg) \le C\frac{p\sigma}{\Gamma}\min\bigg\{1,e^{-\Gamma^2/2\sigma^2}\log\bigg(1+\frac{2\sigma^2}{\Gamma^2} \bigg) \bigg\}+\frac{C}{p^c}.
\eeq
Combining (\ref{eq17}) and (\ref{eq.2}), we conclude that
\begin{align*}
\E [\tau_K(\hat{w}^\top Y,\pi) ] &\le \frac{C\sigma}{p\Gamma}\min\bigg\{1,e^{-\Gamma^2/2\sigma^2}\log\bigg(1+\frac{2\sigma^2}{\Gamma^2} \bigg)  \bigg\}+\frac{Ce^{-\Gamma^2/2\sigma^2}}{p(\Gamma/\sigma+\sqrt{8/\pi})}+\frac{C}{p^{c+2}},
\end{align*}
which completes the proof, as the bound $\E [\tau_K(\hat{w}^\top Y,\pi) ]\le 1$ is trivial.

{\bf{Proof of (\ref{eq.1}).}}
Following the same argument as the proof of Theorem \ref{linear.thm.1}, we have for $1\le i\le p-1$ and $\delta = \frac{\sigma^2(n+\log p)(\lambda_1^2+\sigma^2p)}{\lambda_1^4}$, $P_i \le P\big( w^\top \phi_{ij}\ge- \sqrt{2\delta}\|\phi_{i,i+1}\|_2 \big) +P\big( |1-(\hat{w}^\top w)^2|>\delta  \big)$,
where the second term can be bounded using Lemma \ref{conc.lem}. For the first term, by Lemma \ref{lem.term1}, for $\lambda_1^4(\Theta')\ge \sigma^2(p+\sigma^2\lambda_1^2(\Theta'))(\log p+n)$, we have
\begin{align*}
&P_i\le \Phi\bigg( {C\sqrt{\delta}}\bigg[ (\sqrt{n}+\sqrt{\log p})^2+\frac{\Gamma_i^2}{\sigma^2} +\frac{|\Gamma_i|}{\sigma}\sqrt{\log p} \bigg]^{1/2}+\frac{\Gamma_i}{\sigma} \bigg)+\frac{C}{p^c}.
\end{align*}
Using same argument as the proof of Theorem \ref{linear.thm.1}, it holds that $P_i \le \Phi\big(\frac{\Gamma_i}{\sigma} \big)+\frac{C}{p^c}.$
Equation (\ref{eq.1}) then follows by using formula 7.1.13 of \cite{abramowitz1965handbook} that	$\Phi(-t) < \frac{2}{t+\sqrt{t^2+8/\pi}}\phi(t)$ for $t\ge 0$.

{\bf{Proof of (\ref{eq.2})}.} For the set of indices $S=\{ (i,j): 1\le i<j\le p, j>i+1\}$, we further divide it into two subsets $S_1 = \{(i,j):1\le i<j\le p, j> i+\lfloor\sigma\sqrt{C\log p}/\Gamma \rfloor \}$ and $S_2 = \{(i,j):1\le i<j\le p, i+1<j\le  i+\lfloor\sigma\sqrt{C\log p}/\Gamma \rfloor \}$ for some constant $C>0$. Apparently we have the decomposition
\begin{align} \label{eq.2.5}
\sum_{(i,j): j> i+1}P\bigg( \sum_{k=1}^n\hat{w}_k(X_{ki}-X_{kj})\ge 0 \bigg) &= \sum_{(i,j)\in S_1}P\bigg( \sum_{k=1}^n\hat{w}_k(X_{ki}-X_{kj})\ge 0 \bigg)\nonumber\\
&\quad+\sum_{(i,j)\in S_2}P\bigg( \sum_{k=1}^n\hat{w}_k(X_{ki}-X_{kj})\ge 0 \bigg)
\end{align}
For the first term, by construction, it can be shown using the same argument (see supplementary materials) in Theorem \ref{linear.thm.1} that 
\beq \label{eq.3}
\sum_{(i,j)\in S_1}P\bigg( \sum_{k=1}^n\hat{w}_k(X_{ki}-X_{kj})\ge 0 \bigg) \le \frac{C|S_1|}{p^c} \le \frac{C}{p^{c_0}}.
\eeq
Now for the second term in (\ref{eq.2.5}), similar argument yields,  for $(i,j)\in S_2$,
\begin{align*}
P\bigg( \sum_{k=1}^n \hat{w}_k (X_{ki}-X_{kj})\ge 0 \bigg) & \le \frac{c\exp(-|i-j|^2\Gamma^2/(2\sigma^2))}{|i-j|\Gamma/\sigma+\sqrt{|i-j|^2\Gamma^2/\sigma^2+8/\pi}}+\frac{C}{p^c}.
\end{align*}
Note that, on the one hand,
\[
\frac{e^{-|i-j|^2\Gamma^2/(2\sigma^2)}}{|i-j|\Gamma/\sigma+\sqrt{|i-j|^2\Gamma^2/\sigma^2+8/\pi}}\le \frac{\sigma e^{-|i-j|^2\Gamma^2/2\sigma^2}}{|i-j|\Gamma}.
\]
We have
\[
\sum_{(i,j)\in S_2}\frac{\sigma e^{-|i-j|^2\Gamma^2/2\sigma^2}}{|i-j|\Gamma} = \frac{\sigma}{\Gamma}\sum_{k=2}^{p\land \lfloor\sqrt{\log p}/\Gamma\rfloor} \bigg(\frac{e^{-k^2\Gamma^2/2\sigma^2}p}{k}-e^{-k^2\Gamma^2/2\sigma^2} \bigg)=T_1-T_2.
\]
For the rest of the proof, we assume $\Gamma \le \sigma\sqrt{\log p}/2$, otherwise the set $S_2$ will vanish.
Then
\begin{align*}
T_1 &= \frac{\sigma p}{\Gamma}\sum_{k=2}^{p\land \lfloor\sigma\sqrt{\log p}/\Gamma\rfloor} \frac{e^{-k^2\Gamma^2/2\sigma^2}}{k}\le  \frac{\sigma p}{\Gamma}\int_{1}^{p\land \sigma\sqrt{\log p}/\Gamma} \frac{e^{-x^2\Gamma^2/2\sigma^2}}{x}dx
\end{align*}
where the last inequality used monotonicity of the integrand.
The integral in the last inequality, after change of variable, can be bounded by an exponential integral $\text{Ei}(\Gamma^2/2\sigma^2)$, which has an upper bound
\begin{align*}
\int_{1}^{p\land \sigma\sqrt{\log p}/\Gamma} \frac{e^{-x^2\Gamma^2/2\sigma^2}}{x}dx&=\frac{1}{2}\int_{\Gamma^2/2\sigma^2}^{(p^2\Gamma^2/2\sigma^2)\land(\log p/2)} \frac{e^{-t}}{t}dt \le \frac{1}{2}\int_{\Gamma^2/2\sigma^2}^{\infty} \frac{e^{-t}}{t}dx\le e^{-\Gamma^2/2\sigma^2}\log\bigg(1+\frac{2\sigma^2}{\Gamma^2} \bigg)
\end{align*}
so that $T_1 \le \frac{\sigma p}{\Gamma}e^{-\Gamma^2/2\sigma^2}\log\big(1+\frac{2\sigma^2}{\Gamma^2} \big).$
For $T_2$, we have $T_2 = \frac{\sigma}{\Gamma}\sum_{k=2}^{p\land\lfloor\sigma\sqrt{\log p}/\Gamma\rfloor} {e^{-k^2\Gamma^2/2\sigma^2}} \ge \frac{\sigma}{\Gamma}e^{-2\Gamma^2/\sigma^2}.$
Therefore,
\beq \label{eq.4}
\sum_{(i,j)\in S_2}P\bigg( \sum_{k=1}^n\hat{w}_k(X_{ki}-X_{kj})\ge 0 \bigg) \le \frac{C\sigma p}{\Gamma}e^{-\Gamma^2/2\sigma^2}\log\bigg(1+\frac{2\sigma^2}{\Gamma^2} \bigg)+\frac{C}{p^c}.
\eeq
On the other hand, note that
\[
\frac{e^{-|i-j|^2\Gamma^2/(2\sigma^2)}}{|i-j|\Gamma/\sigma+\sqrt{|i-j|^2\Gamma^2/\sigma^2+8/\pi}}\le c{e^{-|i-j|^2\Gamma^2/(2\sigma^2)}}.
\]
We have
\begin{align*}
&\sum_{(i,j)\in S_2}{e^{-|i-j|^2\Gamma^2/(2\sigma^2)}}=\sum_{k=2}^{p\land\lfloor\sigma\sqrt{\log p}/\Gamma\rfloor} {pe^{-k^2\Gamma^2/(2\sigma^2)}}-\sum_{k=2}^{p\land\lfloor\sigma\sqrt{\log p}/\Gamma\rfloor} {ke^{-k^2\Gamma^2/(2\sigma^2)}}\\
&\le p\int_1^{\infty}e^{-k^2\Gamma^2/2\sigma^2}dk-2e^{-2\Gamma^2/\sigma^2}\le Cp{\sigma}/{\Gamma}.
\end{align*}
Thus
\beq \label{eq.5}
\sum_{(i,j)\in S_2} P\bigg( \sum_{k=1}^n\hat{w}_k(X_{ki}-X_{kj})\ge 0 \bigg) \le Cp{\sigma}/{\Gamma}+\frac{C'}{p^c}.
\eeq
Combining (\ref{eq.4}) and (\ref{eq.5}), we have
\beq \label{eq.6}
\sum_{(i,j)\in S_2} P\bigg( \sum_{k=1}^n\hat{w}_k(X_{ki}-X_{kj})\ge 0 \bigg) \le C\frac{p\sigma}{\Gamma}\min\bigg\{1,e^{-\Gamma^2/2\sigma^2}\log\bigg(1+\frac{2\sigma^2}{\Gamma^2} \bigg) \bigg\}+\frac{C}{p^c}
\eeq
Combining (\ref{eq.3}) and (\ref{eq.6}), we have (\ref{eq.2}). 
\qed

\paragraph{Proof of Theorem \ref{general.thm.1} and Theorem \ref{general.thm.2}.} Here we only provide a sketch of the proofs. We refer the readers to our Supplementary Material for detailed proofs. The proofs follow essentially from the same argument as the proofs of Theorem \ref{linear.thm.1} and Theorem \ref{linear.thm.2}, respectively. However, in place of Lemma \ref{conc.lem} used therein, we need the following lemma that provides a perturbation bound for the leading eigenvector of approximate rank-one matrices, which could be of independent interest.

\bel \label{conc.lem3}
Suppose $p\gtrsim n$ and $\lambda_1^2(\Theta') \ge \lambda_2^2(\Theta')+C\sigma^2(n+\sqrt{np})$ for some $C>0$. Let $w={u'}_1$ be the first left singular vector of $\Theta'$, it follows that,
\[
P\bigg( |1-(\hat{w}^\top w)^2|\le \frac{C\sigma^2(\lambda_1^2(\Theta')+\sigma^2p)(n+\log p)}{(\lambda_1^2(\Theta')-\lambda_2^2(\Theta'))^2}\bigg)\ge 1-\frac{C}{p^c}.
\]
\eel
The proof of Lemma \ref{conc.lem3} is nontrivial, which depends on a combination of the generic perturbation bound obtained by \cite{cai2018rate} and new concentration inequalities of approximate rank-one matrices (see Supplementary Materials).

\paragraph{Proof of Theorem \ref{linear.low.bnd.thm.1}.} The proof relies on the following lemma adapted from \citep{tsybakov2009introduction}.

\bel \label{low.bnd.lem}
Assume that for some integer $M\ge 2$ there exist distinct parameters $\theta_0,...,\theta_M$ from the parameter space $\Theta$ and mutually absolutely continuous probability measures $P_0,...,P_M$ with $P_j=P_{\theta_j}$ for $j=0,1,...,M$, defined on a common probability space $(\Omega, \mathcal{F})$ such that the averaged K-L divergence $\frac{1}{M}\sum_{j=1}^M D(P_j,P_0) \le \frac{1}{8}\log M.$ Then, for every measurable mapping $\hat{\theta}: \Omega\to \Theta$,
\[
\max_{j=0,...,M}P_j(\hat{\theta}\ne \theta_j)\ge \frac{\sqrt{M}}{\sqrt{M}+1}\bigg( \frac{3}{4}-\frac{1}{2\sqrt{\log M}}\bigg).
\]
\eel

We construct the $(M+1)=p$ points parameter space as follows. We define $p$ permutations from $\mathcal{S}_p$ as an identity plus $(p-1)$ consecutive swaps, i.e., $\pi_0=id$, $\pi_k = (k, k+1)$ for $k=1,...,p-1.$
The signal matrix $\Theta_0=a\eta^\top$ where $a=(1,...,1)^\top \in \R^n$ and $\eta = (0,\delta,...,(p-1)\delta)^\top \in \R^p$, $\delta = \frac{\sigma}{4}\sqrt{\log p/n}$. In this way, we have $\Gamma=\|a\|_2\cdot \min_{1\le i\le p-1}|\eta_j-\eta_{j+1}| = \frac{\sigma}{4}\sqrt{\log p}$. Let $P_k$ corresponds to the joint probability measure of $Y$ under $(\Theta_0,\pi_k)$ for $k=0,1,...,p-1$, and let $p_k$ be the pdf of $P_k$, we have $p_0(x) = \prod_{i=1}^n\prod_{j=1}^p\phi_{\eta_j}(x_{ij})$, $p_k(x)=\prod_{i=1}^n\prod_{j=1}^p\phi_{\eta_{\pi_k(j)}}(x_{ij})$ for $k=1,...,p-1$,
where $\phi_{\mu}$ is the pdf of Gaussian distribution $N(\mu,\sigma^2)$. Now we calculate the KL-divergence
\begin{align*}
&D(P_k,P_0) = \int \log\bigg(\frac{p_k(x)}{p_0(x)} \bigg)p_0 (x)dx=\int \frac{n}{2\sigma^2}\sum_{i=1}^p[(x_{1j}-\eta_{\pi_k(j)})^2-(x_{1j}-\eta_j)^2]p_0(x)dx\\
&=\frac{n\delta^2}{\sigma^2}=\frac{\log p}{16}.
\end{align*}
Then, we have for $p\ge 10$, $\frac{1}{p-1}\sum_{k=1}^{p-1}D(P_k,P_0)=\frac{\log p}{16} \le \frac{1}{8}\log (p-1).$
It follows from Lemma \ref{low.bnd.lem} that, $ \inf_{\hat{\pi}}\sup_{(\pi,\Theta)\in \mathcal{D}_1}P(\hat{\pi}\ne \pi) \ge \inf_{\hat{\pi}}\max_{j=0,...,p-1}P_j(\hat{\pi}\ne \pi_j) \ge 0.3$ as long as $p\ge 10$. In addition, $\inf_{\hat{\pi}}\sup_{(\pi,\Theta)\in \mathcal{D}'_1}P(\hat{\pi}\ne \pi) \ge \inf_{\hat{\pi}}\sup_{(\pi,\Theta)\in \mathcal{D}_1}P(\hat{\pi}\ne \pi)$ as $\mathcal{D}_1\subset\mathcal{D}'_1$.
\qed

\paragraph{Proof of Theorem \ref{linear.low.bnd.thm.2}.} The proof relies on the following lemma from \cite{tsybakov2009introduction}.

\bel \label{lower.lem}
Assume that $M\ge 2$ and suppose that $\Theta$ contains elements $\theta_0,\theta_1,...,\theta_M$ such that: (i) $d(\theta_j,\theta_k)\ge 2s>0$ for any $0\le j<k\le M$; (ii) for any $j=1,...,M$, $\frac{1}{M}\sum_{j=1}^M D(P_j,P_0)\le \alpha \log M$ with $0<\alpha<1/8$ and $P_j=P_{\theta_j}$ for $j=0,1,...,M$. Then
\[
\inf_{\hat{\theta}}\sup_{\theta\in\Theta}P_\theta(d(\hat{\theta},\theta)\ge s) \ge \frac{\sqrt{M}}{1+\sqrt{M}}\bigg( 1-2\alpha-\sqrt{\frac{2\alpha}{\log M}}\bigg) >0.
\]
\eel

We also need the following sphere packing lemma proved by \cite{mao2017minimax}, which is a direct consequence of the well-celebrated Varshamov-Gilbert bound.

\bel \label{packing.lem}
For any $r< p/2$, there exists a subset $\mathcal{Q}_r$ of $\mathcal{S}_p$ such that (i) $\log|\mathcal{Q}_r|\ge \frac{r}{5}\log(p/r)$, (ii) for any elements $\pi_1,\pi_2\in \mathcal{Q}_r$, we have ${p\choose 2}\cdot\tau_K(\pi_1,\pi_2)\ge r$, and (iii) for any $\pi\in\mathcal{Q}_r$, we have  $\|\pi-id\|_2^2\le 2r$.
\eel

For $t/\sigma \ge 2$, we set $r=\frac{p\sigma}{t}e^{-t^2/2\sigma^2}<p/2$. Let $\pi_0=id$ and $\pi_1,...,\pi_{|\mathcal{Q}_r|}$ be the elements of $\mathcal{Q}_r$. The signal matrix $\Theta_0=a\eta^\top$ where $a=(1/\sqrt{320n},...,1/\sqrt{320n})^\top\in \R^n$ and $\eta=(t,...,pt)^\top\in \R^p$. Let $P_k$ be the joint probability measure of $Y$ under $(\Theta_0,\pi_k)$ for $k=0,1,...,|\mathcal{Q}_r|$, and let $p_k$ be the pdf of $P_k$. By Lemma \ref{packing.lem}, the KL-divergence
\begin{align*}
&D(P_k,P_0) = \int \log\bigg(\frac{p_k(x)}{p_0(x)} \bigg)p_0 (x)dx=\frac{t^2}{320\sigma^2}\|\pi_k-id\|_2^2\le \frac{pt}{160\sigma} e^{-t^2/2\sigma^2}.
\end{align*}
and therefore
\[
\frac{1}{p-1}\sum_{k=1}^{p-1}D(P_k,P_0) \le \frac{pt}{160\sigma} e^{-t^2/2\sigma^2} \le \frac{p\sigma}{80t}e^{-t^2/2\sigma^2}\log\bigg( \frac{t}{\sigma}e^{t^2/2\sigma^2} \bigg)\le\frac{1}{16}\log |\mathcal{Q}_r|.
\]
Without loss of generality, we assume $|\mathcal{Q}_r| \ge 2$. By Lemma \ref{lower.lem}, it then follows that,
\[
\inf_{\hat{\pi}}\sup_{(\Theta,\pi)\in \mathcal{D}_2(t)} P\bigg(\tau_K(\hat{\pi},\pi) \ge \frac{\sigma}{2pt}e^{-t^2/2\sigma^2}\bigg) \ge C_1,
\]
for some absolute constant $C_1>0$. By Markov's inequality, we have
\[
\inf_{\hat{\pi}}\sup_{(\Theta,\pi)\in \mathcal{D}_2(t)} \E [\tau_K(\hat{\pi},\pi)] \ge \frac{\sigma}{2pt}e^{-t^2/2\sigma^2} \inf_{\hat{\pi}}\sup_{(\Theta,\pi)\in \mathcal{D}_2(t)} P\bigg(\tau_K(\hat{\pi},\pi) \ge \frac{\sigma}{2pt}e^{-t^2/2\sigma^2}\bigg) \ge \frac{C_1\sigma}{pt}e^{-t^2/2\sigma^2}.
\]
The relationship $\inf_{\hat{\pi}}\sup_{(\Theta,\pi)\in \mathcal{D}'_2(t)} \E [\tau_K(\hat{\pi},\pi)] \ge \inf_{\hat{\pi}}\sup_{(\Theta,\pi)\in \mathcal{D}_2(t)} \E [\tau_K(\hat{\pi},\pi)]$ follows from $\mathcal{D}_L\subset \mathcal{D}$. The rate $1/p^2$ follows by setting $t=C_2\sigma\sqrt{\log p}$ for some $C_2>0$.
\qed

\section*{SUPPLEMENTARY MATERIALS}

In our online Supplemental Materials, we prove Theorem 3-4, Proposition 1-4, as well as the technical lemmas. Some supplementary simulations, figures and tables are included in the appendix.

\bibliographystyle{chicago}
\bibliography{reference}

\newpage

\title{Supplement to ``Optimal Permutation Recovery in Permuted Monotone Matrix Model"}
\author{Rong Ma$^1$, T. Tony Cai$^2$ and Hongzhe Li$^1$ \\
	Department of Biostatistics, Epidemiology and Informatics$^1$\\
	Department of Statistics$^2$\\
	University of Pennsylvania\\
	Philadelphia, PA 19104}
\date{}
\maketitle
\thispagestyle{empty}

\begin{abstract}
In this Supplementary Material, we prove Theorem 3 \& 4 and Proposition 1 to 4 in the main paper and the technical lemmas. The numerical comparisons of $\hat{\pi}$ and an SVD-based estimator $\tilde{\pi}$ are also included as an Appendix. 
\end{abstract}

\setcounter{section}{0}

\section{Proofs of Theorem 3 and 4}

\paragraph{Proof of Theorem 3.} As in the proof of Theorem 1, it suffices to bound the probability $P\big(  \sum_{k=1}^n\hat{w}_k(X_{ki}-X_{k,i+1})\ge  0 \big)$. Note that $\hat{w}$ is the first eigenvector of
\[
A=X\bigg(I-\frac{1}{p}ee^\top\bigg)\bigg(I-\frac{1}{p}ee^\top\bigg)  X^\top = TT^\top
\]
where $T\in \R^{n\times p}$ and $T$ admits the decomposition
\[
T = \Theta\bigg(I-\frac{1}{p}ee^\top\bigg)+Z\bigg(I-\frac{1}{p}ee^\top\bigg) =\Theta'+E \in \R^{n\times p}
\]
where $E_{ij}\sim N(0, (p-1)\sigma^2/p)$  and $\Theta'$ has SVD
\[
\Theta'=U'D'V'^\top.
\]
Define $w = u'_1 \in \R^n$ as the first eigenvector of $\Theta'\Theta'^\top$. Following the same argument that leads to (13) of the main paper, we have,
up to a change of sign for $\hat{w}$,
\begin{align*}
P\bigg( \sum_{k=1}^n \hat{w}_k (X_{ki}-X_{k,i+1})\ge 0 \bigg) \le P\bigg( w^\top \phi_{i,i+1}\ge- \sqrt{2\delta}\|\phi_{i,i+1}\|_2 \bigg) +P\bigg( |1-(\hat{w}^\top w)^2|>\delta  \bigg).
\end{align*}
The following lemmas parallel Lemma 1 and Lemma 2 in the proof of Theorem 1. In particular, Lemma 8 provides a perturbation bound for the leading eigenvector of approximate rank-one matrices, which could be of independent interest.

\bel \label{lem2.term1}
Under the conditions of Theorem 3, let $\Xi_i=\|\Theta'_{.i}-\Theta'_{.i+1}\|_2$. Then for any $\delta>0$, we have
\begin{align}\label{term1.eq2}
&P\bigg( w^\top \phi_{i,i+1}\ge- \sqrt{2\delta}\|\phi_{i,i+1}\| \bigg)\\
&\quad\le \Phi\bigg( {C\sqrt{\delta}}\bigg[(\sqrt{n}+\sqrt{\log p})^2+\frac{\Xi_i^2}{\sigma^2}+\sqrt{\log p}\frac{\Xi_i}{\sigma} \bigg]^{1/2}+\frac{w^\top(\Theta'_{\cdot i}-\Theta'_{\cdot i+1})}{\sigma} \bigg)+\frac{C}{p^c}
\end{align}
for some universal constant $C>0$.
\eel

\bel  \label{conc.lem3}
Suppose $n\lesssim p$ and $\lambda_1^2(\Theta') \ge \lambda_2^2(\Theta')+C\sigma^2(n+\sqrt{np})$ for some $C>0$, it follows that
\[
P\bigg( |1-(\hat{w}^\top w)^2|\le \frac{C\sigma^2(\lambda_1^2(\Theta')+\sigma^2p)(n+\log p)}{(\lambda_1^2(\Theta')-\lambda_2^2(\Theta'))^2}\bigg)\ge 1-\frac{C}{p^c}.
\]
\eel

Combining Lemma \ref{lem2.term1} and Lemma \ref{conc.lem3}, since
\beq \label{cond.lambda.delta2}
\lambda_1^2(\Theta')-\lambda_2^2(\Theta')>C_0\sigma^2(n+\sqrt{np}) \quad \text{and let}\quad \delta = \frac{C_0\sigma^2(\lambda_1^2(\Theta')+\sigma^2p)(n+\log p)}{(\lambda_1^2(\Theta')-\lambda_2^2(\Theta'))^2}
\eeq
for some $C_0>0$, 
we have
\begin{align*}
&P\bigg( \sum_{k=1}^n \hat{w}_k (X_{ki}-X_{k,i+1})\ge 0 \bigg) \\
&\quad\le \Phi\bigg( {C\sqrt{\delta}}\bigg[ (\sqrt{n}+\sqrt{\log p})^2+\frac{\Xi_i^2}{\sigma^2} +\frac{\Xi_i}{\sigma}\sqrt{\log p} \bigg]^{1/2}+\frac{w^\top(\Theta'_{.i}-\Theta'_{.i+1})}{\sigma} \bigg)+\frac{C}{p^c}
\end{align*}
for some $C,c>0$.
The rest of the analysis is divided into several cases.
\paragraph{Case 1. $\log p \lesssim n$.} In this case, under (\ref{cond.lambda.delta2}), we have
\begin{align*}
&P\bigg( \sum_{k=1}^n \hat{w}_k (X_{ki}-X_{k,i+1})\ge 0 \bigg) \\
&\quad\le \Phi\bigg( {C\sqrt{\delta}}\bigg[ n+\frac{\|\Theta'_{.i}-\Theta'_{.i+1}\|_2^2}{\sigma^2} +\frac{\|\Theta'_{.i}-\Theta'_{.i+1}\|_2}{\sigma}\sqrt{\log p} \bigg]^{1/2}+\frac{w^\top(\Theta'_{.i}-\Theta'_{.i+1})}{\sigma}  \bigg)+\frac{C}{p^c}
\end{align*}
then
\begin{enumerate}
	\item if $\|\Theta'_{.i}-\Theta'_{.i+1}\|_2 \lesssim \sigma\sqrt{n}$, we have
	\begin{align*}
	P\bigg( \sum_{k=1}^n \hat{w}_k (X_{ki}-X_{k,i+1})\ge 0 \bigg) &\le \Phi\bigg( C\sqrt{\delta n}+\frac{w^\top(\Theta'_{.i}-\Theta'_{.i+1})}{\sigma}  \bigg)+\frac{C}{p^c}\le \frac{C'}{p^c},
	\end{align*}
	where the last inequality follows from 
	\[
	\sqrt{\log p}\lesssim \Gamma/\sigma,\quad \quad 
	\lambda_1^2(\Theta') \gtrsim \lambda_2^2(\Theta')+\frac{{\sigma^2n(n+\log p)}}{\Gamma^2/\sigma^2}+\frac{\sigma^2\sqrt{np(n+\log p)}}{\Gamma/\sigma}.
	\]
	\item $\|\Theta'_{.i}-\Theta'_{.i+1}\|_2 \gtrsim \sigma\sqrt{n}$, we have 
	\begin{align*}
	P\bigg( \sum_{k=1}^n \hat{w}_k (X_{ki}-X_{k,i+1})\ge 0 \bigg) &\le \Phi\bigg( {C\sqrt{\delta}}\frac{\|\Theta'_{.i}-\Theta'_{.i+1}\|_2}{\sigma} +\frac{w^\top(\Theta'_{.i}-\Theta'_{.i+1})}{\sigma} \bigg)+\frac{C}{p^c}\le \frac{C'}{p^c},
	\end{align*}
	where the last inequality follows from
	\[
	\Gamma/\sigma\gtrsim \sqrt{n}\gtrsim \sqrt{\log p},\quad\quad  \lambda_1^2(\Theta') \gtrsim  \lambda_2^2(\Theta')+\frac{\sigma^2\Xi^2(n+\log p)}{\Gamma^2}+\frac{\sigma^2\Xi\sqrt{p(n+\log p)}}{\Gamma}.
	\]
\end{enumerate}
This completes the proof of case 1.

\paragraph{Case 2. $\log p \gtrsim n$.} In this case, under (\ref{cond.lambda.delta2}), we have
\begin{align*}
&P\bigg( \sum_{k=1}^n \hat{w}_k (X_{ki}-X_{k,i+1})\ge 0 \bigg) \\
&\quad\le \Phi\bigg( {C\sqrt{\delta}}\bigg[ \log p+\frac{\|\Theta'_{.i}-\Theta'_{.i+1}\|_2^2}{\sigma^2} +\frac{\|\Theta'_{.i}-\Theta'_{.i+1}\|_2}{\sigma}\sqrt{\log p} \bigg]^{1/2}+\frac{w^\top(\Theta'_{.i}-\Theta'_{.i+1})}{\sigma} \bigg)+\frac{C}{p^c}
\end{align*}
In addition, since $\|\Theta'_{.i}-\Theta'_{.i+1}\|_2 \gtrsim \sigma\sqrt{\log p}$, we have
\begin{align*}
P\bigg( \sum_{k=1}^n \hat{u}_k (X_{ki}-X_{k,i+1})\ge 0 \bigg) &\le \Phi\bigg( \frac{C\sqrt{\delta}}{\sigma}\|\Theta'_{.i}-\Theta'_{.i+1}\|_2+\frac{w^\top(\Theta'_{.i}-\Theta'_{.i+1})}{\sigma} \bigg)+\frac{C}{p^c}\le \frac{C'}{p^c}
\end{align*}
where the last inequality follows from
\[
\Gamma/\sigma \gtrsim \sqrt{\log p}\gtrsim \sqrt{n},\quad\quad  \lambda_1^2(\Theta') \gtrsim  \lambda_2^2(\Theta')+\frac{\sigma^2\Xi^2(n+\log p)}{\Gamma^2}+\frac{\sigma^2\Xi\sqrt{p(n+\log p)}}{\Gamma}.
\]
This completes the proof of case 2. As a result, it follows that, up to a change of sign of $\hat{w}$,
\[
P((\frak{r}(\hat{w}^\top X))^{-1}\ne id) \le \sum_{i=1}^{p-1}P\bigg(  \sum_{k=1}^n\hat{w}_k(X_{ki}-X_{k,i+1})\ge  0 \bigg) \le  \frac{C}{p^c}
\]
for some constant $C,c>0$.
\qed

\paragraph{Proof of Theorem 4.}  Similar to the proof of Theorem 2, we write
\begin{align*}
\E [\tau_K(\hat{w}^\top Y,\pi) ] &= \frac{2}{p(p-1)}\sum_{i< j} P\bigg( \sum_{k=1}^n\hat{w}_k(X_{ki}-X_{kj})\ge 0 \bigg)
\end{align*}
where
\begin{align*}
\sum_{i< j} P\bigg( \sum_{k=1}^n\hat{w}_k(X_{ki}-X_{kj})\ge 0 \bigg)
&= \sum_{(i,j): j=i+1} P\bigg(\sum_{k=1}^n\hat{w}_k(X_{ki}-X_{kj})\ge 0 \bigg)\\
&\quad+\sum_{(i,j): j> i+1}P\bigg( \sum_{k=1}^n\hat{w}_k(X_{ki}-X_{kj})\ge 0 \bigg).
\end{align*}
In the following, we first show
\beq \label{eq.ge.1}
P\bigg(\sum_{k=1}^n\hat{w}_k(X_{ki}-X_{k,i+1})\ge 0 \bigg) \le \frac{ce^{-\Gamma^2/2\sigma^2}}{\Gamma/\sigma+\sqrt{\Gamma^2/\sigma^2+8/\pi}}+\frac{C}{p^c}
\eeq
for some constant $c>0$ and therefore
\[
\sum_{(i,j): j=i+1} P\bigg(\sum_{k=1}^n\hat{w}_k(X_{ki}-X_{kj})\ge 0 \bigg) \le \frac{cpe^{-\Gamma^2/2\sigma^2}}{\Gamma/\sigma+\sqrt{\Gamma^2/\sigma^2+8/\pi}}+\frac{C}{p^c}
\]
Then we show that 
\beq \label{eq.ge.2}
\sum_{(i,j): j> i+1}P\bigg( \sum_{k=1}^n\hat{w}_k(X_{ki}-X_{kj})\ge 0 \bigg) \le C\frac{p\sigma}{\Gamma}\min\bigg\{1,e^{-\Gamma^2/2\sigma^2}\log\bigg(1+\frac{2\sigma^2}{\Gamma^2} \bigg) \bigg\}+\frac{C}{p^c}
\eeq
and conclude that
\begin{align*}
\E [\tau_K(\hat{w}^\top Y,\pi) ] &\le \frac{C\sigma}{p\Gamma}\min\bigg\{1,e^{-\Gamma^2/2\sigma^2}\log\bigg(1+\frac{2\sigma^2}{\Gamma^2} \bigg)  \bigg\}+\frac{Ce^{-\Gamma^2/2\sigma^2}}{p(\Gamma/\sigma+\sqrt{8/\pi})}+\frac{C}{p^{c+2}}.
\end{align*}
The bound $\E [\tau_K(\hat{w}^\top Y,\pi) ]\le 1$ is trivial.

{\bf{Proof of (\ref{eq.ge.1}).}}
Following the same argument as in the proof of Theorem 3, we have for $1\le i\le p-1$ and $\delta = \frac{\sigma^2(n+\log p)(\lambda_1^2(\Theta')+\sigma^2p)}{(\lambda_1^2(\Theta')-\lambda_2^2(\Theta'))^2}$, 
\begin{align*}
P\bigg( \sum_{k=1}^n \hat{w}_k (X_{ki}-X_{k,i+1})\ge 0 \bigg) &\le P\bigg( w^\top \phi_{ij}\ge- \sqrt{2\delta}\|\phi_{i,i+1}\|_2 \bigg) +P\bigg( |1-(\hat{w}^\top w)^2|>\delta  \bigg)\\
&\le P\bigg( w^\top \phi_{i,i+1}\ge- \sqrt{2\delta}\|\phi_{i,i+1}\|_2\bigg) + \frac{C}{p^c}
\end{align*}
where the last inequality follows from Lemma \ref{conc.lem3}. For the first term, by Lemma \ref{lem2.term1}, for $\delta<1$,
\begin{align*}
P\bigg( \sum_{k=1}^n \hat{w}_k (X_{ki}-X_{k,i+1})\ge 0 \bigg) &\le \Phi\bigg( {C}\sqrt{\delta}\bigg( (\sqrt{n}+\sqrt{\log p})^2+ \frac{\|\Theta'_{.i}-\Theta'_{.i+1}\|_2^2}{\sigma^2}\\
&\quad +\frac{\|\Theta'_{.i}-\Theta'_{.i+1}\|_2}{\sigma}\sqrt{\log p} \bigg)^{1/2}+\frac{w^\top(\Theta'_{.i}-\Theta'_{.i+1})}{\sigma} \bigg)+\frac{C}{p^c}.
\end{align*}
The rest of the analysis is divided into several cases.
\paragraph{Case 1. $\log p \le n$.} In this case, we have 
\begin{align*}
P\bigg( \sum_{k=1}^n \hat{w}_k (X_{ki}-X_{k,i+1})\ge 0 \bigg) &\le \Phi\bigg( {C}\sqrt{\delta}\bigg(n+ \frac{\|\Theta'_{.i}-\Theta'_{.i+1}\|_2^2}{\sigma^2} +\frac{\|\Theta'_{.i}-\Theta'_{.i+1}\|_2}{\sigma}\sqrt{\log p} \bigg)^{1/2}\\
&\quad+\frac{w^\top(\Theta'_{.i}-\Theta'_{.i+1})}{\sigma} \bigg)+\frac{C}{p^c}.
\end{align*}
In addition,  
\begin{enumerate}
	\item if $\|\Theta'_{.i}-\Theta'_{.i+1}\|_2 \lesssim \sigma\sqrt{n}$, we have
	\begin{align*}
	P\bigg( \sum_{k=1}^n \hat{w}_k (X_{ki}-X_{k,i+1})\ge 0 \bigg) &\le \Phi\bigg( C\sqrt{\delta n}+\frac{w^\top(\Theta'_{.i}-\Theta'_{.i+1})}{\sigma} \bigg)+\frac{C}{p^c}\\
	&\le  \Phi\bigg( \frac{C\Gamma}{\sigma} \bigg)+\frac{C}{p^c}
	\end{align*}
	where the last inequality follows from
	\[
	\lambda_1^2(\Theta') \gtrsim \lambda_2^2(\Theta')+\frac{\sigma^2{n(n+\log p)}}{\Gamma^2/\sigma^2} +\frac{\sigma^2\sqrt{np(n+\log p)}}{\Gamma/\sigma}.
	\]
	Thus, using Formula 7.1.13 from \cite{abramowitz1965handbook},
	\beq \label{as.form}
	\Phi(-x) < \frac{2}{t+\sqrt{t^2+8/\pi}}\phi(x),
	\eeq
	for $1\le i\le p-1$, we have
	\[
	P\bigg( \sum_{k=1}^n \hat{w}_k (X_{ki}-X_{k,i+1})\ge 0 \bigg) \le \frac{C}{\Gamma/\sigma+\sqrt{\Gamma^2/\sigma^2+8/\pi}}\exp(-\Gamma^2/(2\sigma^2))+\frac{C}{p^c}.
	\]
	\item if $\|\Theta'_{.i}-\Theta'_{.i+1}\|_2 \gtrsim \sigma\sqrt{n}$, we have 
	\begin{align*}
	P\bigg( \sum_{k=1}^n \hat{w}_k (X_{ki}-X_{k,i+1})\ge 0 \bigg) &\le \Phi\bigg( \frac{C\sqrt{\delta}}{\sigma}\|\Theta'_{.i}-\Theta'_{.i+1}\|_2+\frac{w^\top(\Theta'_{.i}-\Theta'_{.i+1})}{\sigma} \bigg)+\frac{C}{p^c}\\
	&\le \frac{C}{\Gamma/\sigma+\sqrt{\Gamma^2/\sigma^2+8/\pi}}\exp(-\Gamma^2/(2\sigma^2))+ \frac{C}{p^c},
	\end{align*}
	where the last inequality follows from (\ref{as.form}) and
	\[
	\lambda_1^2(\Theta') \gtrsim  \lambda^2_2(\Theta')+\frac{\sigma^2\Xi^2({n+\log p})}{\Gamma^2} +\frac{\sigma^2\Xi\sqrt{p(n+\log p)}}{\Gamma}.
	\]
\end{enumerate}
This completes the proof of case 1.

\paragraph{Case 2. $\log p > n$.} In this case, we have 
\begin{align*}
P\bigg( \sum_{k=1}^n \hat{w}_k (X_{ki}-X_{k,i+1})\ge 0 \bigg) &\le \Phi\bigg( {C}\sqrt{\delta}\bigg(\log p+ \frac{\|\Theta'_{.i}-\Theta'_{.i+1}\|_2^2}{\sigma^2} +\frac{\|\Theta'_{.i}-\Theta'_{.i+1}\|_2}{\sigma}\sqrt{\log p} \bigg)^{1/2}\\
&\quad+\frac{w^\top(\Theta'_{.i}-\Theta'_{.i+1})}{\sigma} \bigg)+\frac{C}{p^c}.
\end{align*}
In addition, 
\begin{enumerate}
	\item if $\|\Theta'_{.i}-\Theta'_{.i+1}\|_2 \gtrsim \sigma\sqrt{\log p}$, we have
	\begin{align*}
	P\bigg( \sum_{k=1}^n \hat{u}_k (X_{ki}-X_{k,i+1})\ge 0 \bigg) &\le \Phi\bigg( \frac{C\sqrt{\delta}}{\sigma}\|\Theta'_{.i}-\Theta'_{.i+1}\|_2 +\frac{w^\top(\Theta'_{.i}-\Theta'_{.i+1})}{\sigma} \bigg)+\frac{C}{p^c}\\
	&\le \frac{C}{\Gamma/\sigma+\sqrt{\Gamma^2/\sigma^2+8/\pi}}\exp(-\Gamma^2/(2\sigma^2))+ \frac{C}{p^c}.
	\end{align*}
	\item if $\|\Theta'_{.i}-\Theta'_{.i+1}\|_2 \lesssim \sigma\sqrt{\log p}$, we have
	\begin{align*}
	P\bigg( \sum_{k=1}^n \hat{u}_k (X_{ki}-X_{k,i+1})\ge 0 \bigg) &\le \Phi\bigg( {C\sqrt{\delta \log p}}+\frac{w^\top(\Theta'_{.i}-\Theta'_{.i+1})}{\sigma} \bigg)+\frac{C}{p^c}\\
	&\le \frac{\sqrt{2/\pi}}{\Gamma/\sigma+\sqrt{\Gamma^2/\sigma^2+8/\pi}}\exp(-\Gamma^2/(2\sigma^2))+\frac{C}{p^c}
	\end{align*}
	where the last inequality follows from
	\[
	\lambda_1^2(\Theta') \gtrsim \lambda^2_2(\Theta')+\frac{\sigma^2{\log p(n+\log p)}}{\Gamma^2/\sigma^2}+\frac{\sigma^2\sqrt{p\log p (n+\log p)}}{\Gamma/\sigma}.
	\]
\end{enumerate}
This completes the proof of Case 2 and in sum we have proven (\ref{eq.ge.1}).

\paragraph{{Proof of (\ref{eq.ge.2})}.} Following the same construction of $S_1$ and $S_2$ in the proof of (18) of the main paper, we have
\begin{align} \label{eq.ge.2.5}
\sum_{(i,j): j> i+1}P\bigg( \sum_{k=1}^n\hat{w}_k(X_{ki}-X_{kj})\ge 0 \bigg) &= \sum_{(i,j)\in S_1}P\bigg( \sum_{k=1}^n\hat{w}_k(X_{ki}-X_{kj})\ge 0 \bigg)\nonumber\\
&\quad+\sum_{(i,j)\in S_2}P\bigg( \sum_{k=1}^n\hat{w}_k(X_{ki}-X_{kj})\ge 0 \bigg)
\end{align}
For the first term, for any $(i,j)\in S_1$, we have
\begin{align*}
P\bigg( \sum_{k=1}^n \hat{w}_k (X_{ki}-X_{kj})\ge 0 \bigg) &\le P\bigg( w^\top \phi_{ij}\ge- \sqrt{2\delta}\|\phi_{ij}\|_2 \bigg) +P\bigg( |1-(\hat{w}^\top w)^2|>\delta  \bigg)\\
&\le P\bigg( w^\top \phi_{ij}\ge- \sqrt{2\delta}\|\phi_{ij}\|_2\bigg) + \frac{C}{p^c}
\end{align*}
where the last inequality follows from Lemma \ref{conc.lem3}.
To bound $P\bigg( w^\top \phi_{ij}\ge- \sqrt{2\delta}\|\phi_{ij}\|_2\bigg)$, similar argument as in Theorem 3 implies, for $\delta<1$,
\begin{align*}
P\bigg( \sum_{k=1}^n \hat{w}_k (X_{ki}-X_{kj})\ge 0 \bigg)  &\le \Phi\bigg( {C}\sqrt{\delta}\bigg( (\sqrt{n}+\sqrt{\log p})^2+ \frac{\|\Theta'_{.i}-\Theta'_{.j}\|_2^2}{\sigma^2} +\frac{\|\Theta'_{.i}-\Theta'_{.j}\|_2}{\sigma}\sqrt{\log p} \bigg)^{1/2}\\
&\quad+\frac{w^\top(\Theta'_{.i}-\Theta'_{.i+1})}{\sigma} \bigg)+\frac{C}{p^c}.
\end{align*}
The rest of the analysis is divided into two cases.
\paragraph{Case 1. $\log p \lesssim n$.} In this case, we have 
\begin{align*}
P\bigg( \sum_{k=1}^n \hat{w}_k (X_{ki}-X_{kj})\ge 0 \bigg) &\le \Phi\bigg( {C}\sqrt{\delta}\bigg( n+ \frac{\|\Theta'_{.i}-\Theta'_{.j}\|_2^2}{\sigma^2} +\frac{\|\Theta'_{.i}-\Theta'_{.j}\|_2}{\sigma}\sqrt{\log p} \bigg)^{1/2}\\
&\quad+\frac{w^\top(\Theta'_{.i}-\Theta'_{.j})}{\sigma} \bigg)+\frac{C}{p^c}.
\end{align*}
In addition, let $\Gamma_{ij}=\lambda_1(\Theta')|v'_{1i}-v'_{1j}|\ge |i-j|\Gamma, $
\begin{enumerate}
	\item if $\|\Theta'_{.i}-\Theta'_{.j}\|_2 \lesssim \sigma\sqrt{n}$, we have
	\begin{align*}
	P\bigg( \sum_{k=1}^n \hat{w}_k (X_{ki}-X_{kj})\ge 0 \bigg) &\le \Phi\bigg( C\sqrt{\delta n}+\frac{w^\top(\Theta'_{.i}-\Theta'_{.j})}{\sigma} \bigg)+\frac{C}{p^c} \le \frac{C}{p^c},
	\end{align*}
	since
	\begin{align*}
	\lambda_1^2(\Theta') &\gtrsim \lambda^2_2(\Theta')+\frac{\sigma^2{n(n+\log p)}}{\Gamma^2/\sigma^2} +\frac{\sigma^2\sqrt{np (n+\log p)}}{\Gamma/\sigma}\\
	& \ge \lambda^2_2(\Theta')+\frac{\sigma^2{n(n+\log p)}}{\Gamma_{ij}^2/\sigma^2} +\frac{\sigma^2\sqrt{np (n+\log p)}}{\Gamma_{ij}/\sigma},
	\end{align*}
	and $\Gamma_{ij}\ge \sigma \sqrt{C\log p}$ for $(i,j)\in S_1$. 
	\item if $\|\Theta'_{.i}-\Theta'_{.j}\|_2\gtrsim \sigma\sqrt{n}$, then let $\Xi_{ij}=\|\Theta'_{.i}-\Theta'_{.j}\|_2 \le |i-j|\Xi$, we have 
	\begin{align*}
	P\bigg( \sum_{k=1}^n \hat{w}_k (X_{ki}-X_{kj})\ge 0 \bigg) &\le \Phi\bigg( \frac{C\sqrt{\delta}}{\sigma}\|\Theta'_{.i}-\Theta'_{.j}\|_2 +\frac{w^\top(\Theta'_{.i}-\Theta'_{.j})}{\sigma} \bigg)+\frac{C}{p^c}\le \frac{C}{p^c},
	\end{align*}
	where the last inequality follows from
	\begin{align*}
	\lambda_1^2(\Theta') &\gtrsim \lambda^2_2(\Theta')+\frac{\sigma^2\Xi^2(n+\log p)}{\Gamma^2} +\frac{\sigma^2\Xi\sqrt{p (n+\log p)}}{\Gamma}\\
	& \ge \lambda^2_2(\Theta')+\frac{\sigma^2\Xi_{ij}^2(n+\log p)}{\Gamma_{ij}^2} +\frac{\sigma^2\Xi_{ij}\sqrt{p (n+\log p)}}{\Gamma_{ij}},
	\end{align*}
	and $\Gamma_{ij} \ge \sigma\sqrt{C\log p}.$
\end{enumerate}
This completes the proof of case 1.

\paragraph{Case 2. $\log p \gtrsim n$.} In this case, we have 
\begin{align*}
P\bigg( \sum_{k=1}^n \hat{w}_k (X_{ki}-X_{kj})\ge 0 \bigg) &\le \Phi\bigg( {C}\sqrt{\delta}\bigg( \log p+ \frac{\|\Theta'_{.i}-\Theta'_{.j}\|_2^2}{\sigma^2} +\frac{\|\Theta'_{.i}-\Theta'_{.j}\|_2}{\sigma}\sqrt{\log p} \bigg)^{1/2}\\
&\quad+\frac{w^\top(\Theta'_{.i}-\Theta'_{.j})}{\sigma} \bigg)+\frac{C}{p^c}.
\end{align*}
In addition, since $\Gamma_{ij}\ge \sigma\sqrt{C\log p}$, similar arguments yield
\begin{align*}
P\bigg( \sum_{k=1}^n \hat{u}_k (X_{ki}-X_{kj})\ge 0 \bigg) & \le \frac{C}{p^c}.
\end{align*}
This completes the proof of Case 2. Combining Case 1 and Case 2, we have proven
\beq \label{eq.ge.3}
\sum_{(i,j)\in S_1}P\bigg( \sum_{k=1}^n\hat{w}_k(X_{ki}-X_{kj})\ge 0 \bigg) \le \frac{C|S_1|}{p^c} \le \frac{C}{p^{c_0}}.
\eeq
Now for the second term in (\ref{eq.ge.2.5}), since $\Gamma |i-j| \le \Gamma_{ij}$, we have for $(i,j)\in S_2$,
\[
P\bigg( \sum_{k=1}^n\hat{w}_k(X_{ki}-X_{kj})\ge 0 \bigg) \le \frac{ce^{-|i-j|^2\Gamma^2/(2\sigma^2)}}{|i-j|\Gamma/\sigma+\sqrt{|i-j|^2\Gamma^2/\sigma^2+8/\pi}}+\frac{C}{p^c}.
\]
Note that, on the one hand,
\[
\frac{ce^{-|i-j|^2\Gamma^2/(2\sigma^2)}}{|i-j|\Gamma/\sigma+\sqrt{|i-j|^2\Gamma^2/\sigma^2+8/\pi}}\le \frac{c\sigma e^{-|i-j|^2\Gamma^2/2\sigma^2}}{|i-j|\Gamma}.
\]
Following similar argument that leads to (22) of the main paper, we have
\beq \label{eq.ge.4}
\sum_{(i,j)\in S_2}P\bigg( \sum_{k=1}^n\hat{w}_k(X_{ki}-X_{kj})\ge 0 \bigg) \le \frac{C\sigma p}{\Gamma}e^{-\Gamma^2/2\sigma^2}\log\bigg(1+\frac{2\sigma^2}{\Gamma^2} \bigg)+\frac{C}{p^c}.
\eeq
On the other hand, note that
\[
\frac{e^{-|i-j|^2\Gamma^2/(2\sigma^2)}}{|i-j|\Gamma/\sigma+\sqrt{|i-j|^2\Gamma^2/\sigma^2+8/\pi}}\le {ce^{-|i-j|^2\Gamma^2/(2\sigma^2)}},
\]
we have
\beq \label{eq.ge.5}
\sum_{(i,j)\in S_2} P\bigg( \sum_{k=1}^n\hat{w}_k(X_{ki}-X_{kj})\ge 0 \bigg) \le Cp{\sigma}/{\Gamma}+\frac{C'}{p^c}.
\eeq
Combining (\ref{eq.ge.4}) and (\ref{eq.ge.5}), we have
\beq \label{eq.ge.6}
\sum_{(i,j)\in S_2} P\bigg( \sum_{k=1}^n\hat{w}_k(X_{ki}-X_{kj})\ge 0 \bigg) \le C\frac{p\sigma}{\Gamma}\min\bigg\{1,e^{-\Gamma^2/2\sigma^2}\log\bigg(1+\frac{2\sigma^2}{\Gamma^2} \bigg) \bigg\}+\frac{C}{p^c}
\eeq
Combining (\ref{eq.ge.3}) and (\ref{eq.ge.6}), we have
\[
\sum_{(i,j): j> i+1}P\bigg( \sum_{k=1}^n\hat{w}_k(X_{ki}-X_{kj})\ge 0 \bigg) \le C\frac{p\sigma}{\Gamma}\min\bigg\{1,e^{-\Gamma^2/2\sigma^2}\log\bigg(1+\frac{2\sigma^2}{\Gamma^2} \bigg) \bigg\}+\frac{C}{p^c},
\]
which leads to (\ref{eq.ge.2}). 
\qed

\section{Proofs of Proposition 1-4}

\paragraph{Proof of Proposition 1.} Since $\Theta\in \mathcal{D}$, for any nonnegative unit vector $w\in \R^n$, the vector $w^\top\Theta \in \R^p$ is monotonic increasing, so that $\frak{r}(w^\top\Theta)=id$. It then follows that
\[
\frak{r}(w^\top\Theta\Pi)=\frak{r}(w^\top\Theta)\circ \pi^{-1}=\pi^{-1}.
\]\qed

\paragraph{Proof of Proposition 2.} The key observation is, for any $\Gamma$, we have 
\beq \label{cond.lambda}
\Lambda \ge C\Gamma^2 p^3.
\eeq
Then if $\Gamma \gtrsim \sigma\sqrt{n}$ holds, apparently GSS can be implied by (\ref{cond.lambda}) and MSG; if $\Gamma \lesssim \sigma\sqrt{n}$, then the results depend on the relative magnitude of $\sigma^2 n^2/\Gamma^2$ and $p$. Specifically, if $p\gtrsim \sigma^2n^2/\Gamma^2$, then $p\gtrsim (\sigma^3n/\Gamma^3)^{2/5}$ implies GSS; if $p\lesssim \sigma^2n^2/\Gamma^2$, then $p\gtrsim (\sigma^4n^2/\Gamma^4)^{1/3}$ implies GSS.
\qed

\paragraph{Proof of Proposition 3.}
Note that by SVD of $\Theta'$, we have 
\beq \label{v1.def}
v'_1 = \argmax_{\|x\|_2=1} x^\top \Theta^\top \Theta x=\argmax_{\|x\|_2=1} \sum_{i=1}^n \bigg( \sum_{j=1}^p x_j \theta'_{ij}\bigg)^2.
\eeq
To prove that $v'_1$ is monotone, we need the following rearrangement inequality.
\bel [Rearrangement Inequality] \label{rearragement.ineq}
If $a_1\ge a_2\ge ...\ge a_n$ and $b_1\ge b_2\ge ... \ge b_n$, then
\[
a_nb_1+...+a_1b_n \le a_{\sigma(1)}b_1+...+a_{\sigma(n)}b_n \le a_1b_1+...+a_nb_n,
\]
where $\sigma$ is any permutation in $S_n$.
\eel
\bel \label{pos.mat.lem}
For any $\Theta'$ defined above, its first left singular vector $u'_1$ is either nonpositive or nonnegative.
\eel
Now since $u'_1$ is either nonnegative or nonpositive, then we know that $\sum_{j=1}^pv_j\theta'_{ij}$ have the same sign for all $i=1,...,n$. By Lemma \ref{rearragement.ineq}, we have that the components of $v'_1$ are either in increasing order or in decreasing order.
\qed

\paragraph{Proof of Proposition 4.} Note that 
\beq 
v'_1 = \argmax_{\|x\|_2=1} x^\top {\Theta'}^\top \Theta' x=\argmax_{\|x\|_2=1} \sum_{i=1}^n ( x^\top \Theta'_{i.})^2.
\eeq
Let $f(x,\Theta',\lambda) = \sum_{i=1}^n ( x^\top \Theta'_{i.})^2+\lambda(\|x\|_2-1)$. By Lagrange's multiplier method, we have
\[
\frac{\partial f}{\partial x_j} =  \sum_{i=1}^n 2\theta_{ij}( x^\top \Theta'_{i.})+2\lambda x_j=0,
\]
for $j=1,...,p$. It follows that
\[
\frac{\sum_{i=1}^n\theta'_{ij}({v'}_1^\top \Theta'_{i.})}{v'_{1j}} = -\lambda,\quad\quad\text{ for all $j=1,...,p$,}
\]
where
\[
\lambda^2=\sum_{j=1}^p \bigg[\sum_{i=1}^n\theta'_{ij}({v'}_1^\top \Theta'_{i.})\bigg]^2 = \|{v'}_1^\top {\Theta'}^\top {\Theta'}\|_2^2.
\]
Thus, by Lemma \ref{pos.mat.lem}, we have
\begin{align*}
|v'_{1,j+1}-v'_{1,j}| &= \frac{\sum_{i=1}^n|\theta'_{i,j+1}-\theta'_{i,j}||{v'}_1^\top \Theta'_{i.}|}{|\lambda|} \\
&\ge \frac{\sum_{i=1}^n\delta\|\Theta_{i.}\|_2|{v'}_1^\top \Theta'_{i.}|}{|\lambda|}\\
&\ge \delta,
\end{align*}
where the last inequality follows from
\[
|\lambda| = \|{v'}_1^\top {\Theta'}^\top {\Theta'}\|_2=\bigg\|\sum_{i=1}^n (v_1^\top\Theta'_{i.})\Theta'_{i.} \bigg\|_2 \le \sum_{i=1}^n \|({v'}_1^\top \Theta'_{i.})\Theta'_{i.} \|_2 \le \sum_{i=1}^n({v'}_1^\top \Theta'_{i.})\|\Theta'_{i.}\|_2.
\]
This completes the proof.
\qed

\section{Proofs of Technical Lemmas}

\paragraph{Proof of Lemma 1.} Note that
\[
\phi_{i,i+1}^\top w = \sum_{k=1}^nw_k(T_{ki}-T_{k,i+1})\sim N\bigg(\|a\|_2(\eta_i-\eta_{i+1}), \frac{2(p-1)\sigma^2}{p}\bigg),
\]
and
\[
\|\phi_{i,i+1}\|_2^2 = \sum_{k=1}^n (T_{ki}-T_{k,i+1})^2,\quad\text{where}\quad T_{ki}-T_{k,i+1} \sim N(a_k(\eta_i-\eta_{i+1}), 2(p-1)\sigma^2/p). 
\]
We construct the standardized chi-square statistic
\[
Y_k^2 = \bigg[\frac{T_{k,i}-T_{k,i+1}-a_k(\eta_i-\eta_{i+1})}{\sigma\sqrt{2(p-1)/{p}}} \bigg]^2.
\]
Define
\begin{align*}
Q&=\frac{p-1}{p}2\sigma^2\cdot \sum_{k=1}^n(Y_k^2-1),
\end{align*}
we have
\begin{align*}
Q&=\|\phi_{i,i+1}\|_2^2+\|a\|_2^2(\eta_i-\eta_{i+1})^2-2(\eta_i-\eta_{i+1})\sum_{k=1}^n(T_{ki}-T_{k,i+1})a_k-\frac{p-1}{p}\cdot 2n\sigma^2.
\end{align*}
Our next result follows from an exponential inequality for chi-square random variables proved by \cite{laurent2000adaptive}.

\bel \label{conc.lem.3}
Let $(Y_1,...,Y_n)$ be i.i.d. Gaussian variables with mean $0$ and variance $1$. Let $a_1,...,a_n$ be nonnegative. We set $\|a\|_\infty = \sup_{1\le i\le n}|a_i|, \|a\|_2^2 = \sum_{i=1}^n a_i^2$. Let $Z = \sum_{i=1}^na_i(Y_i^2-1)$. Then the following inequalities hold for any positive $x$:
\begin{align}
P(Z \ge 2\|a\|_2\sqrt{x}+2\|a\|_\infty x) &\le \exp(-x),\\
P(Z\le -2\|a\|_2\sqrt{x}) &\le \exp(-x).
\end{align}
\eel
By choosing $x=c\log p$, we have, with probability at least $1-\frac{1}{p^c}$,
\begin{align} \label{conc.eq.chi}
\|\phi_{i,i+1}\|_2^2&\le c\sigma^2\sqrt{n\log p}+c\sigma^2 \log p-\|a\|_2^2(\eta_i-\eta_{i+1})^2+2(\eta_i-\eta_{i+1})\sum_{k=1}^n(T_{ki}-T_{k,i+1})a_k+\frac{p-1}{p}\cdot 2n\sigma^2\nonumber\\
&\le c\sigma^2(\sqrt{n}+\sqrt{\log p})^2-\|a\|_2^2(\eta_i-\eta_{i+1})^2+2(\eta_i-\eta_{i+1})\sum_{k=1}^n(T_{ki}-T_{k,i+1})a_k
\end{align}
Note that the term
\[
(\eta_i-\eta_{i+1})\sum_{k=1}^n(T_{ki}-T_{k,i+1})a_k \sim N\bigg( \|a\|_2^2(\eta_i-\eta_{i+1})^2, 2\|a\|_2^2(\eta_i-\eta_{i+1})^2\cdot \frac{p-1}{p}\sigma^2\bigg)
\]
It follows from the tail bound of standard Gaussian distribution $\Phi(-x)\le \frac{1}{x}\phi(x)$ that
\beq \label{conc.eq.norm}
P\bigg( (\eta_i-\eta_{i+1})\sum_{k=1}^n(T_{ki}-T_{k,i+1})a_k\ge \|a\|_2^2(\eta_i-\eta_{i+1})^2+2\sigma\sqrt{c\log p}\cdot \|a\|_2|\eta_i-\eta_{i+1}| \bigg) \le \frac{1}{p^c}.
\eeq
Combining (\ref{conc.eq.chi}) and (\ref{conc.eq.norm}), we have
\[
P\bigg(\|\phi_{i,i+1}\|_2^2\le c'\sigma^2(\sqrt{n}+\sqrt{\log p})^2+\|a\|_2^2(\eta_i-\eta_{i+1})^2+c'\sigma\sqrt{\log p}\|a\|_2|\eta_i-\eta_{i+1}| \bigg) \ge 1-\frac{C}{p^c}
\]
for some constant $C,c,c'>0.$
Thus, 
\begin{align*}
&P\bigg( w^\top \phi_{i,i+1}\ge- \sqrt{2\delta}\|\phi_{i,i+1}\| \bigg)\\
&\quad\le \Phi\bigg( \frac{C\sqrt{\delta}}{\sigma}\bigg[\sigma^2(\sqrt{n}+\sqrt{\log p})^2+\|a\|_2^2(\eta_i-\eta_{i+1})^2 +\|a\|_2|\eta_i-\eta_{i+1}|\sigma\sqrt{\log p} \bigg]^{1/2}+\frac{\|a\|_2(\eta_i-\eta_{i+1})}{\sigma} \bigg)+\frac{C}{p^c}
\end{align*}
\qed

\paragraph{Proof of Lemma 2.} The proof follows from the following deterministic perturbation bound and concentration inequalities adapted from \cite{cai2018rate}.

\bel \label{pert.lem}
Suppose $A\in \R^{n\times p}$, $\tilde{V}=\begin{bmatrix}
V & V_{\perp}\end{bmatrix} \in \mathbb{O}_{n}$ are left singular vectors of $A$, $V\in  \mathbb{O}_{n,r}, V_{\perp}\in  \mathbb{O}_{n,n-r}$ correspond to the first $r$ and last $(n-r)$ singular vectors respectively. $\tilde{W} =\begin{bmatrix}
W & W_{\perp}\end{bmatrix}\in  \mathbb{O}_{n}$ is any orthogonal matrix with $W\in  \mathbb{O}_{n,r}$, $W_{\perp}\in  \mathbb{O}_{n,n-r}$. Given that $\lambda_{r}(W^\top A)>\lambda_{r+1}(A)$, we have
\[
\|\sin \Theta(V,W)\| \le \frac{\lambda_{r}(W^\top A)\| \mathbb{P}_{A^\top W}A^\top W_{\perp}\|^2}{\lambda_r^2(W^\top A)-\lambda_{r+1}^2(A)} \land 1.
\]
\eel

\bel\label{conc.lem.2} 
Suppose $X\in \R^{n\times p}$ is a rank-$r$ matrix with left singular space as $V\in  \mathbb{O}_{n,r}, Z\in \R^{n\times p},$ $Z$ is a i.i.d. sub-Gaussian random matrix with sub-Gaussian parameter $\tau$. $Y=X+Z$. Then there exists constants $C,c$ such that for any $x>0$,
\beq \label{conc.1}
P\big(\lambda_{r}^2(V^\top Y)\le (\lambda_r^2(X)+\tau^2p)(1-x)\big)\le C\exp\big\{ C r-c\tau^{-2}(\lambda_r^2(X)+\tau^2p) (x^2 \land x) \big\},
\eeq
\beq \label{conc.2}
P\big( \lambda_{r+1}^2(Y) \ge \tau^2p(1+x)\big) \le C\exp\big\{ Cn-cp\cdot (x^2 \land x)\big\}
\eeq
Moreover, there exists $C_0,C,c$ such that whenever $\lambda_r^2(X)\ge C_0\tau^2n,$ for any $x>0$ we have
\beq \label{conc.3}
P(\|\mathbb{P}_{V^\top Y}V_{\perp}^\top Y\| \le x) \ge 1- C\exp \big\{ Cn-c\tau^{-2}\min(x^2,x\sqrt{\lambda_r^2(X)+\tau^2p})\big\}-C\exp\big\{ -c\tau^{-2}(\lambda_r^2(X)+\tau^2p) \big\}.
\eeq
\eel

Note that $T=\Theta'+E$ where $\Theta'=\|a\| \cdot w\eta'^{\top}$ is rank 1, and $\hat{w}$ is the first left singular vector of $T$, by Lemma \ref{pert.lem}, we have
\beq
|1-(\hat{w}^\top w)^2|=\|\sin \Theta(w,\hat{w})\|^2 \le  \frac{\lambda_1^2({w}^\top T)\| \mathbb{P}_{{w}^\top T}{w}_{\perp}^\top T\|^2}{(\lambda_1^2({w}^\top T)-\lambda_2^2(T))^2} \land 1.
\eeq
To bound the quantity $\lambda_{1}^2(w^\top T)$, using Lemma \ref{conc.lem.2}, by choosing $x=\frac{\lambda^2}{3(\lambda^2+\sigma^2p)}$ in (\ref{conc.1}), since $\lambda^2 \ge C\sigma^2(n+\sqrt{pn})$ for some sufficiently large $C$, or $\sigma^2(\lambda^2+\sigma^2 p)(\sigma^2n+\log p)/\lambda^4 \le 1$, we have
\[
c\sigma^{-2}\min\bigg\{ \frac{\lambda^4}{\lambda^2+\sigma^2p},  \lambda^2 \bigg\} -C \ge \frac{c_0\sigma^{-2}\lambda^4}{\lambda^2+\sigma^2p}-C\ge c_1\log p.
\]
Thus
\beq \label{conc.4}
P\big(\lambda_{1}^2(w^\top T)\le \frac{2\lambda^2}{3}+\sigma^2p\big)\le C\exp\bigg\{C -\frac{c_0\sigma^{-2}\lambda^4}{\lambda^2+p} \bigg\}\le C\exp\{ -c_1\log p\}.
\eeq
To bound $\lambda_2^2(T)$, if $p\le n+\log p$, by choosing $x={\frac{n+\log p}{p}}$ in (\ref{conc.2}), we have
\[
c \min\bigg\{ \frac{(n+\log p)^2}{p}, {n+\log p} \bigg\} -Cn\ge c\log p
\]
for some $C,c>0$, so that
\beq \label{conc.5}
P\big( \lambda_{2}^2(T) \ge \sigma^2(p+n+\log p)\big) \le C\exp\big\{ Cn-cp\cdot x^2 \land x\big\}\le C\exp \{-c\log p\} = \frac{C}{p^c}.
\eeq
If $p>n+\log p$, by choosing $x=\sqrt{\frac{n+\log p}{p}}$ in (\ref{conc.2}), we have
\[
c \min\bigg\{ \sqrt{p(n+\log p)}, {n+\log p} \bigg\} -Cn\ge c\log p
\]
for some $C,c>0$, so that
\beq \label{conc.5.1}
P\big( \lambda_{2}^2(T) \ge \sigma^2(p+\sqrt{p(n+\log p)})\big) \le C\exp\big\{ Cn-cp\cdot x^2 \land x\big\}\le C\exp \{-c\log p\} = \frac{C}{p^c}.
\eeq
Lastly, to bound $\|\mathbb{P}_{w^\top T}w_{\perp}^\top T \|$, choosing $x= \sigma\sqrt{n+\log p}$ in (\ref{conc.3}), since $\lambda^2 \ge \sigma^2n$, we have
\[
\sigma^2(n+\log p) \le c(\lambda^2+\sigma^2p)
\]
for sufficiently large $c>0$ and therefore
\[
c\sigma^{-2} \min\bigg\{ x^2,x\sqrt{\lambda^2+\sigma^2p} \bigg\} -Cn \ge c\log p.
\]
So that
\beq \label{conc.6}
P(\|\mathbb{P}_{w^\top T}w_{\perp}^\top T \| \le \sigma\sqrt{n+\log p}) \ge 1- \frac{C}{p^c}-C\exp\big\{ -c\sigma^{-2}(\lambda^2+\sigma^2p) \big\}.
\eeq
Combining (\ref{conc.4}) (\ref{conc.5}) (\ref{conc.5.1}) and (\ref{conc.6}), by the fact that $x/(x-y)^2$ is decreasing in $x$ and increasing in $y$ for $x>y>0$, we have
\[
P\bigg(  \frac{\lambda^2({w}^\top T)\| \mathbb{P}_{{w}^\top T}{w}_{\perp}^\top T\|^2}{(\lambda_r^2({w}^\top T)-\lambda_2^2(T))^2} \le C'\frac{\sigma^2(\lambda^2+\sigma^2p)(n+\log p)}{\lambda^4} \bigg) \ge 1-Cp^{-c}.
\]
\qed

\paragraph{Proof of Lemma \ref{lem2.term1}.} Note that
\[
\phi_{i,i+1}^\top w = \sum_{k=1}^nw_k(T_{ki}-T_{k,i+1})\sim N\bigg(w^\top(\Theta'_{\cdot i}-\Theta'_{\cdot i+1}), \frac{2(p-1)\sigma^2}{p}\bigg),
\]
and
\[
\|\phi_{i,i+1}\|_2^2 = \sum_{k=1}^n (T_{ki}-T_{k,i+1})^2,\quad\text{where}\quad T_{ki}-T_{k,i+1} \sim N(\theta'_{ki}-\theta'_{k,i+1}, 2(p-1)\sigma^2/p). 
\]
We construct the standardized chi-square statistic
\[
Y_k^2 = \bigg[\frac{T_{k,i}-T_{k,i+1}-(\theta'_{ki}-\theta'_{k,i+1})}{\sigma\sqrt{2(p-1)/{p}}} \bigg]^2.
\]
Define
\begin{align*}
Q&=\frac{p-1}{p}2\sigma^2\cdot \sum_{k=1}^n(Y_k^2-1),
\end{align*}
we have
\begin{align*}
Q&=\|\phi_{i,i+1}\|_2^2+\|\Theta'_{\cdot i}-\Theta'_{\cdot i+1}\|_2^2-2\sum_{k=1}^n(T_{ki}-T_{k,i+1})(\theta'_{ki}-\theta'_{k,i+1})-\frac{p-1}{p}\cdot 2n\sigma^2.
\end{align*}
By Lemma \ref{conc.lem.3}, choosing $x=c\log p$, we have, with probability at least $1-\frac{1}{p^c}$,
\begin{align} \label{conc.eq.chi.2}
\|\phi_{i,i+1}\|_2^2&\le c\sigma^2\sqrt{n\log p}+c\sigma^2 \log p-\|\Theta'_{\cdot i}-\Theta'_{\cdot i+1}\|_2^2+2\sum_{k=1}^n(T_{ki}-T_{k,i+1})(\theta'_{ki}-\theta'_{k,i+1})+\frac{p-1}{p}\cdot 2n\sigma^2\nonumber\\
&\le c\sigma^2(\sqrt{n}+\sqrt{\log p})^2-\|\Theta'_{\cdot i}-\Theta'_{\cdot i+1}\|_2^2+2\sum_{k=1}^n(T_{ki}-T_{k,i+1})(\theta'_{ki}-\theta'_{k,i+1}).
\end{align}
Note that the term
\[
\sum_{k=1}^n(T_{ki}-T_{k,i+1})(\theta'_{ki}-\theta'_{k,i+1}) \sim N\bigg( \|\Theta'_{\cdot i}-\Theta'_{\cdot i+1}\|_2^2, 2\|\Theta'_{\cdot i}-\Theta'_{\cdot i+1}\|_2^2\cdot \frac{p-1}{p}\sigma^2\bigg)
\]
It follows from the tail bound of standard Gaussian distribution $\Phi(-x)\le \frac{1}{x}\phi(x)$ that
\beq \label{conc.eq.norm.2}
P\bigg( \sum_{k=1}^n(T_{ki}-T_{k,i+1})(\theta'_{ki}-\theta'_{k,i+1})\ge \|\Theta'_{\cdot i}-\Theta'_{\cdot i+1}\|_2^2+2\sigma\sqrt{c\log p}\cdot \|\Theta'_{\cdot i}-\Theta'_{\cdot i+1}\|_2 \bigg) \le \frac{1}{p^c}.
\eeq
Combining (\ref{conc.eq.chi.2}) and (\ref{conc.eq.norm.2}), we have
\[
P\bigg(\|\phi_{i,i+1}\|_2^2\le c'\sigma^2(\sqrt{n}+\sqrt{\log p})^2+\|\Theta'_{\cdot i}-\Theta'_{\cdot i+1}\|_2^2+c'\sigma\sqrt{\log p}\|\Theta'_{\cdot i}-\Theta'_{\cdot i+1}\|_2 \bigg) \ge 1-\frac{C}{p^c}
\]
for some constant $C,c,c'>0.$
Thus, 
\begin{align*}
&P\bigg( w^\top \phi_{i,i+1}\ge- \sqrt{2\delta}\|\phi_{i,i+1}\| \bigg)\\
&\quad\le \Phi\bigg( {C\sqrt{\delta}}\bigg[(\sqrt{n}+\sqrt{\log p})^2+\frac{\|\Theta'_{\cdot i}-\Theta'_{\cdot i+1}\|_2^2}{\sigma^2}+c'\sqrt{\log p}\frac{\|\Theta'_{\cdot i}-\Theta'_{\cdot i+1}\|_2}{\sigma} \bigg]^{1/2}+\frac{w^\top(\Theta'_{\cdot i}-\Theta'_{\cdot i+1})}{\sigma} \bigg)+\frac{C}{p^c}
\end{align*}
\qed

\paragraph{Proof of Lemma \ref{conc.lem3}.} The proof of Lemma \ref{conc.lem3} depends on Lemma \ref{pert.lem} as well as the following concentration inequalities.

\bel\label{conc.lem.4} 
Suppose $X\in \R^{n\times p}$ is a rank-$r$ matrix with first left singular vector $v_1\in \R^n$,  $Z$ is a i.i.d. sub-Gaussian random matrix with sub-Gaussian parameter $\tau$. $Y=X+Z$. Then there exists constants $C,c$ such that for any $x>0$, $\lambda_1^2(X) \ge 3\sqrt{\log (n\lor p)}$
\beq \label{conc4.1}
P\bigg( \lambda_{1}^2(v_1^\top Y)  \le  \lambda_1^2(X)+\tau^2(p-\sqrt{p\log (p \lor n)}-\log (p\lor n) ) \bigg)  \le \frac{C}{p^c},
\eeq
\beq \label{conc4.2}
P(\lambda_2^2(Y)> (\lambda^2_2(X)+\tau^2p)(1+ t)) \le C\exp(Cn-c\tau^{-2}(\lambda_n^2+ \tau^2p)t^2\land t).
\eeq
Moreover, there exists $C_0,C,c$ such that for any $x>0$ we have
\beq \label{conc4.3}
P(\|\mathbb{P}_{Y^\top v_1}Y^\top {v_1}_{\perp}\| < t) \ge 1-C\exp (Cn-c\tau^{-2}\min(t^2,t\sqrt{\lambda_1^2(X)+\tau^2p} ))-C\exp\big\{ -c\tau^{-2}(\lambda_1^2(X)+\tau^2p) \big\}.
\eeq
\eel

Following a similar argument as the proof of Lemma 2, by taking $t=\sigma\sqrt{\frac{n+\log p}{\lambda_n^2+\sigma^2p}}$ in (\ref{conc4.2}), since $n\lesssim p$, we have
\[
P(\lambda_2^2(Y)>\lambda_2^2(\Theta')+\tau^2p)\le \frac{C}{p^c}.
\]
Similarly, by taking $t=\sigma\sqrt{n+\log p}$, we have
\[
P(\|\mathbb{P}_{Y^\top w}Y^\top {w}_{\perp}\|^2 < \sigma^2(n+\log p)) \ge 1-\frac{C}{p^c}.
\]
Then Lemma \ref{pert.lem} implies, if $\lambda_1^2(\Theta') - \lambda_2^2(\Theta')\ge C\sigma^2(n+\sqrt{np})$ for some $C>0$,
\[
P\bigg( |1-(\hat{w}^\top w)^2|\le \frac{C\sigma^2(\lambda_1^2(\Theta')+\sigma^2p)(n+\log p)}{(\lambda_1^2(\Theta')-\lambda_2^2(\Theta'))^2}\bigg)\ge 1-\frac{C}{p^c}.
\]
\qed

\paragraph{Proof of Lemma \ref{pos.mat.lem}.} By SVD of $\Theta'=(\theta'_{ij})\in \R^{n\times p}$, its first left singular vector
\beq 
u'_1 = \argmax_{\|x\|_2=1} x^\top \Theta'{\Theta'}^\top x=\argmax_{\|x\|_2=1} \sum_{j=1}^p \bigg( \sum_{i=1}^n x_i \theta'_{ij}\bigg)^2.
\eeq
In the following, we show that, for any unit vector $x\in \R^n$,
\beq \label{ineq.2}
\sum_{j=1}^p \bigg( \sum_{i=1}^n x_i \theta'_{ij}\bigg)^2\le \sum_{j=1}^p \bigg( \sum_{i=1}^n |x_i| \theta'_{ij}\bigg)^2,
\eeq
from which we conclude that $u'_1$ is either nonpositive or nonnegative. Toward this end, note that
\[
\sum_{j=1}^p \bigg( \sum_{i=1}^n x_i \theta'_{ij}\bigg)^2=\sum_{j=1}^p\sum_{i=1}^nx^2_i{\theta'}^2_{ij}+\sum_{i\ne k}x_ix_k\bigg(\sum_{j=1}^p{\theta'}_{ij}{\theta'}_{kj}\bigg).
\]
Then the inequality (\ref{ineq.2}) follows from
\[
\sum_{i\ne k}x_ix_k\bigg(\sum_{j=1}^p{\theta'}_{ij}{\theta'}_{kj}\bigg)\le \sum_{i\ne k}|x_ix_k|\bigg(\sum_{j=1}^p{\theta'}_{ij}{\theta'}_{kj}\bigg),
\]
which is true as long as $\sum_{j=1}^p{\theta'}_{ij}{\theta'}_{kj}\ge 0$ for all pairs $i\ne k$. We conclude this proof by showing that

\indent {\bf Fact.} For any nondecreasing vectors $a,b\in \R^n$, such that $\sum_{i=1}^n a_i=0$. Then it follows that $a^\top b\ge 0$.\\
To see this, note that since $a$ and $b$ are both nondecreasing, there exist a constant $\delta$ such that the components of $b+\delta\cdot \bf{1}$ has the same sign as $a$. Hence the claim follows from $0\le a^\top(b+\delta\cdot {\bf{1}}) = a^\top b+\delta a^\top {\bf{1}}=a^\top b.$
\qed

\paragraph{Proof of Lemma \ref{conc.lem.4}.} Note that
\[
\lambda_{1}^2(v_1^\top Y) = \|v_1^\top Y\|_2^2 = \lambda_1^2(X)+2\lambda_1(X)v_1^\top Z u_1+\|v_1^\top Z\|_2^2.
\]
The linear term $v_1^\top Z u_1 = \sum_{i,j}v_{1i}u_{1j}Z_{ij}$ is subgaussian with parameter $c\| v_1u_1^\top\|_F=c$ for some constant $c$ only depending on $\tau$. Therefore, by concentration inequality for subgaussian random variables 
\beq \label{conc4.1.1}
P ( v_1^\top Z u_1 \ge t ) \le \exp\bigg\{-\frac{t^2}{2c\tau^2} \bigg\}.
\eeq
On the other hand, the quadratic term $\|v_1^\top Z\|_2^2 = \sum_{i=1}^p (v_1^\top Z_i)^2$ where $Z_i$ for $i=1,...,p$ are the $i$-th column of $Z$. Each of $(v_1^\top Z_i)^2$ is subexponential with mean
\[
\E(v_1^\top Z_i)^2 = \E Z_{ij}^2.
\]
Then concentration inequality for subexponential random variables yields
\beq \label{conc4.1.2}
P\bigg( \bigg|\sum_{i=1}^p (v_1^\top Z_i)^2- p \E Z_{ij}^2\bigg| \ge t \bigg) \le 2\exp\bigg\{-c\min\bigg( \frac{t^2}{\tau^4p} , \frac{t}{\tau^2} \bigg)  \bigg\}
\eeq
for some constant $c>0$. By taking $t=\tau\sqrt{\log (p \lor n)}$ in (\ref{conc4.1.1}), we have
\[
P ( v_1^\top Z u_1 \ge \tau\sqrt{\log (p \lor n)} ) \le \frac{C}{p^c}.
\]
Taking $t=\tau^2\sqrt{p\log (p \lor n)}$ if $\sqrt{\log (p \lor n)}\le \sqrt{p}$ and $t=\tau^2\log (p\lor n)$ otherwise in (\ref{conc4.1.2}), we have
\[
P\bigg( \sum_{i=1}^p (v_1^\top Z_i)^2 \le  \tau^2p- \tau^2\sqrt{p\log (p \lor n)}-\tau^2\log (p\lor n) \bigg)  \le \frac{C}{p^c}.
\]
Combining these, we have, with probability at least $1-O(p^{-c})$,
\begin{align*}
\lambda_{1}^2(v_1^\top Y) &=\lambda_1^2(X)+2\tau\lambda_1(X)v_1^\top Z u_1+\|v_1^\top Z\|_2^2\\
&\ge  \lambda_1^2(X)-2\tau\lambda_1(X)\sqrt{\log (p \lor n)} +\tau^2p- \tau^2\sqrt{p\log (p \lor n)}-\tau^2\log (p\lor n).
\end{align*}
If in addition {$\lambda_1(X)\ge 3\tau\sqrt{\log (n\lor p)}$}, then
\[
P\bigg( \lambda_{1}^2(v_1^\top Y)  \le  \lambda_1^2(X)+\tau^2p-\tau^2\sqrt{p\log (p \lor n)}-\tau^2\log (p\lor n) \bigg)  \le \frac{C}{p^c}.
\]
For $ \lambda_{2}^2(Y) $, note that
\[
\lambda_{2}(Y) = \min_{\text{rank}(B)\le 1}\|Y-B\|\le\| Y-[v_1\quad {\bf{0}}]^\top\cdot Y  \| = \| {v_1}_{\perp}^\top Y\|
\]
for $ {v_1}_{\perp} \in \R^{n\times (n-1)}$.
It suffices to obtain an upper bound for $\lambda_1({v_1}_{\perp}^\top Y)$ with high probability. Next,
\begin{align}\label{v1Y}
\|{v_1}_{\perp}^\top Y\|_2^2 &= \|{v_1}_{\perp}^\top Y Y^\top {v_1}_{\perp} \| \nonumber \\
&\le \|\E{v_1}_{\perp}^\top Y Y^\top {v_1}_{\perp} \| +\|{v_1}_{\perp}^\top Y Y^\top {v_1}_{\perp}-\E{v_1}_{\perp}^\top Y Y^\top {v_1}_{\perp} \|\nonumber \\
&=  \lambda_2^2(X)+\tau^2p+\|{v_1}_{\perp}^\top Y Y^\top {v_1}_{\perp}-\E{v_1}_{\perp}^\top Y Y^\top {v_1}_{\perp} \|
\end{align}
Define the normalization matrix $M\in \R^{(n-1)\times (n-1)}$ as
\[
M = \begin{bmatrix}
(\lambda_2^2(X)+\tau^2p)^{-1/2} & & & & \\
& \ddots &  &&\\
&   & (\lambda_r^2(X)+\tau^2p)^{-1/2} && \\
&&& \ddots &\\
&&&& (\lambda_n^2(X)+\tau^2p)^{-1/2}
\end{bmatrix},
\]
we have
\[
M^\top {v_1}_{\perp}^\top \E Y Y^\top {v_1}_{\perp} M=I_{n-1}.
\]
Let $Q={v_1}_{\perp}^\top Y Y^\top {v_1}_{\perp}-\E{v_1}_{\perp}^\top Y Y^\top {v_1}_{\perp}$, we have
\[
\|Q\| =\| (M^{-1})^\top M^\top Q M M^{-1}\|\le \|M^\top Q M \| \|M^{-1}\|^2.
\]
By construction we have
\[
\|M^{-1}\| = (\lambda_2^2(X)+\tau^2p)^{1/2},
\]
then
\[
\|Q\| \le (\lambda_2^2(X)+\tau^2p)\| M^\top {v_1}_{\perp}^\top  Y Y^\top {v_1}_{\perp} M-M^\top {v_1}_{\perp}^\top \E Y Y^\top {v_1}_{\perp} M\|.
\]
Now it suffices to obtain a concentration inequality for $\| M^\top {v_1}_{\perp}^\top  Y Y^\top {v_1}_{\perp} M-I_{n-1}\|$. The main idea is to use the $\epsilon$-net argument to split the spectral norm deviation to the deviations of single random variables, which can be further controlled by the Hanson-Wright inequality. Specifically, for any unit vector $u\in \R^{n-1}$, we have
\begin{align*}
&u^\top M^\top {v_1}_{\perp}^\top  Y Y^\top {v_1}_{\perp} M u- u^\top I_{n-1} u\\
&= (u^\top M^\top {v_1}_{\perp}^\top  X X^\top {v_1}_{\perp} M u - \E u^\top M^\top {v_1}_{\perp}^\top  X X^\top {v_1}_{\perp} M u)\\
&\quad +2(u^\top M^\top {v_1}_{\perp}^\top  X Z^\top {v_1}_{\perp} M u - \E u^\top M^\top {v_1}_{\perp}^\top  X Z^\top {v_1}_{\perp} M u)\\
&\quad + (u^\top M^\top {v_1}_{\perp}^\top  Z Z^\top {v_1}_{\perp} M u - \E u^\top M^\top {v_1}_{\perp}^\top ZZ^\top {v_1}_{\perp} M u)\\
&=2 (X^\top  {v_1}_{\perp} M u )^\top Z^\top ({v_1}_{\perp} M u)+({v_1}_{\perp} Mu)^\top(ZZ^\top - \tau^2 pI_n) ({v_1}_{\perp} Mu).
\end{align*}
In the following, we shall bound the two terms separately.

To bound the second term, for any fixed unit vector $u\in \R^{n-1}$, we vectorize $Z\in \R^{n\times p}$ into $\text{vec}(Z)\in\R^{np}$ as
\[
\text{vec}(Z)=(Z_{11},Z_{21},...,Z_{n1},Z_{12},...,Z_{n2},...,Z_{1p},Z_{np})^\top.
\]
We also define the block diagonal matrix
\[
D = \begin{bmatrix}
({v_1}_{\perp} Mu)({v_1}_{\perp} Mu)^\top & & \\
& \ddots &  \\
&   &  ({v_1}_{\perp} Mu)({v_1}_{\perp} Mu)^\top
\end{bmatrix} \in \R^{np\times np}.
\]
It then follows that 
\[
({v_1}_{\perp} Mu)^\top(ZZ^\top - \tau^2 pI_n) ({v_1}_{\perp} Mu) = \text{vec}(Z)^\top D \cdot \text{vec}(Z)-\E \text{vec}(Z)^\top D\cdot \text{vec}(Z).
\]
Besides,
\[
\| D\| = \| ({v_1}_{\perp} Mu)({v_1}_{\perp} Mu)^\top\| = \|Mu\|_2^2\le \|M\|^2=(\lambda^2_n(X)+\tau^2p)^{-1},
\]
\[
\|D\|_F^2=p\| ({v_1}_{\perp} Mu)({v_1}_{\perp} Mu)^\top\|_F^2\le p\| Mu\|_2^4\le p(\lambda^2_n(X)+\tau^2p)^{-2}.
\]
By Hansen-Wright inequality \citep{rudelson2013hanson},
\begin{align} \label{conc.term2}
&P( |({v_1}_{\perp} Mu)^\top(ZZ^\top - \tau^2 pI_n) ({v_1}_{\perp} Mu)|>t  ) \nonumber \\
&= P(| \text{vec}(Z)^\top D \cdot \text{vec}(Z)-\E \text{vec}(Z)^\top D\cdot \text{vec}(Z) | >t) \nonumber \\
&\le 2\exp\bigg\{ -c\min \bigg( \frac{ t^2(\lambda^2_n(X)+\tau^2p)^2}{\tau^4p}, \frac{t(\lambda^2_n(X)+\tau^2p)}{\tau^2} \bigg)  \bigg\}
\end{align}
for some $c>0$.

Next, we bound the first term
\begin{align*}
(X^\top  {v_1}_{\perp} M u )^\top Z^\top ({v_1}_{\perp} M u)&=\text{tr}(Z^\top ({v_1}_{\perp} M u)(X^\top  {v_1}_{\perp} M u )^\top)\\
&=\text{vec}(Z)^\top \cdot \text{vec}(({v_1}_{\perp} M u)(X^\top  {v_1}_{\perp} M u )^\top).
\end{align*}
Since
\[
X^\top  {v_1}_{\perp} M=U\begin{bmatrix}
0 & ...& 0 &0\\
\lambda_2(X)(\lambda_2^2(X)+\tau^2p)^{-1/2} && &0 \\
& \ddots & & \vdots \\
&& \lambda_n(X)(\lambda_n^2(X)+\tau^2p)^{-1/2} & 0
\end{bmatrix},
\]
we know $\|X^\top  {v_1}_{\perp} M\|\le 1$ and
\begin{align*}
\| \text{vec}(({v_1}_{\perp} M u)(X^\top  {v_1}_{\perp} M u )^\top)\|_2^2 &= \|({v_1}_{\perp} M u)(X^\top  {v_1}_{\perp} M u )^\top\|_F^2\\
&=\|{v_1}_{\perp} M u\|_2^2\cdot \|X^\top  {v_1}_{\perp} M u\|_2^2 \\
&\le \|M\|^2 \le (\lambda_n^2+\tau^2 p)^{-1}.
\end{align*}
By the concentration inequality for i.i.d. subgaussian random variables, we have
\beq \label{conc.term1} 
P(|(X^\top  {v_1}_{\perp} M u )^\top Z^\top ({v_1}_{\perp} M u)|>t)\le C\exp\bigg( -c\frac{t^2 (\lambda_n^2+\tau^2 p)}{\tau^2}\bigg)
\eeq
for some constant $C,c>0$. Combining (\ref{conc.term1}) and (\ref{conc.term2}), we have, for any fixed unitary $u\in \R^{n-1}$,
\[
P(|u^\top M^\top {v_1}_{\perp}^\top  Y Y^\top {v_1}_{\perp} M u- u^\top I_{n-1} u|>t) \le C\exp(-c\tau^{-2}(\lambda_n^2+\tau^2 p)t^2\land t)
\]
for all $t>0$. Next, we use the following lemma proved by \cite{cai2018rate} concerning the $\epsilon$-net argument for unit ball.

\bel \label{enet.lem}
For any $p\ge 1$, denote $\mathbb{B}^p = \{x\in \R^p: \|x\|_2\le 1\}$ as the $p$-dimensional unit ball in the Euclidean space. Suppose $K\in \R^{p_1\times p_2}$ is a random matrix. Then we have for $t>0$,
\[
P(\|K\|\ge 3t)\le 7^{p_1+p_2}\cdot \max_{u\in \mathbb{B}^{p_1},v\in \mathbb{B}^{p_2}} P(|u^\top K v|\ge t).
\]
\eel

It then follows that
\beq \label{conc.M}
P(\| M^\top {v_1}_{\perp}^\top  Y Y^\top {v_1}_{\perp} M -  I_{n-1} \|>t) \le C\exp(Cn-c\tau^{-2}(\lambda_n^2+\tau^2 p)t^2\land t).
\eeq
Recall (\ref{v1Y}), we have
\[
P(\lambda_2^2(Y)> \lambda^2_2(X)+\tau^2p+ (\lambda_2^2(X)+\tau^2 p) t) \le C\exp(Cn-\tau^{-2}(\lambda_n^2+\tau^2 p)t^2\land t).
\]
Finally, we consider $\| \mathbb{P}_{v_1^\top Y} v_{1\perp}^\top Y\|$. Define the constant
\[
m = (\lambda_1^2(X)+\tau^2p)^{-1/2}.
\]
It then follows that
\[
m^2 v_1^\top \E YY^\top v_1 = 1.
\]
Since
\begin{align*}
\| \mathbb{P}_{Y^\top v_1} Y^\top v_{1\perp} \| &= \| \mathbb{P}_{m Y^\top v_1} Y^\top v_{1\perp}\|\\
&=\|m Y^\top v_1((mY^\top v_1)^\top( mY^\top v_1))^{-1} (mY^\top v_1)^\top Y^\top v_{1\perp} \|\\
&\le \|m Y^\top v_1\|_2^{-1}  \|(m Y^\top v_1)^\top Y^\top v_{1\perp}\|.
\end{align*}
In the following we analyze $\|m Y^\top v_1\|_2$ and $\|  (m Y^\top v_1)^\top Y^\top v_{1\perp}\|$ separately. 

Since
\[
\|m Y^\top v_1\|^2 = m^2 |v_1^\top YY^\top v_1| \ge 1 - |m^2v_1^\top YY^\top v_1-m^2 v_1^\top \E YY^\top v_1|,
\]
following the same argument that leads to (\ref{conc.M}), we have
\[
P(|m^2v_1^\top YY^\top v_1| >1- t) \ge 1-C\exp(C-c\tau^{-2}(\lambda_1^2(X)+\tau^2p)t^2\land t).
\]
Now set $t=1/2$, we can choose $C_0$ large enough, such that $\lambda_1^2(X) \ge \tau^2C_0$ and therefore $C-c\tau^{-2}(\lambda_1^2(X)+\tau^2p)\le -c'\tau^{-2}(\lambda_1^2(X)+\tau^2p)$ for some $c,c'>0$. In this case,
\[
P(|m^2v_1^\top YY^\top v_1| >1/2) \ge 1-C\exp(-c\tau^{-2}(\lambda_1^2(X)+\tau^2p)).
\]
For $\|  (m Y^\top v_1)^\top Y^\top v_{1\perp}\|$, note that $mv_1^\top XX^\top v_{1\perp}=0$, we have
\begin{align*}
(m Y^\top v_1)^\top Y^\top v_{1\perp} &= mv_1^\top (X+Z)(X+Z)^\top v_{1\perp}\\
&= mv_1^\top XZ^\top v_{1\perp}+mv_1^\top ZX^\top v_{1\perp}+mv_1^\top ZZ^\top v_{1\perp}.
\end{align*}
Following similar idea of the proof of (\ref{conc.M}), we can show for any unit vector $u\in \R^{n-1}$,
\[
P(|mv_1^\top XZ^\top v_{1\perp} u| >t )\le C\exp\bigg( \frac{-ct^2}{\tau^2\| (v_{1\perp} u) (mv_1^\top X) \|_F^2} \bigg) \le C\exp(-c\tau^{-2}t^2),
\]
\[
P(|mv_1^\top ZZ^\top v_{1\perp} u| >t )=P(|mv_1^\top (ZZ^\top-\tau^2 p I_n) v_{1\perp} u| >t )\le C\exp(-c\tau^{-2}\min(t^2,\sqrt{\lambda_1^2(X)+\tau^2p} t)).
\]
By the $\epsilon$-net argument again (Lemma \ref{enet.lem}), we have
\[
P(\|  (m Y^\top v_1)^\top Y^\top v_{1\perp}\|>t)\le C\exp (Cn-c\tau^{-2}\min(t^2,t\sqrt{\lambda_1^2(X)+\tau^2p})).
\]
\qed

\appendix
\section{Comparison with an SVD-based Esitmator}

\begin{figure}[h!]
	\centering
	\includegraphics[angle=0,width=16cm]{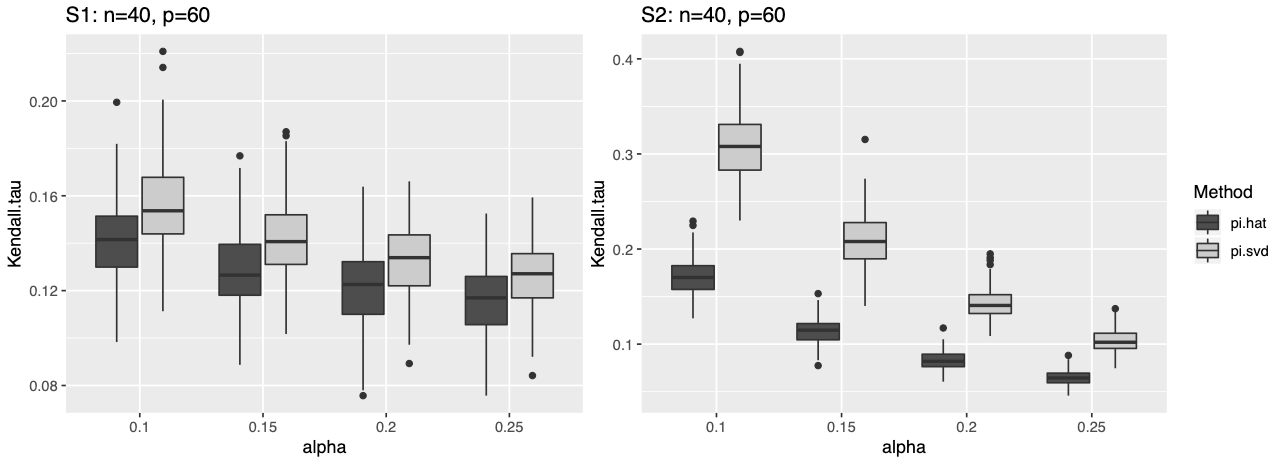}
	\includegraphics[angle=0,width=16cm]{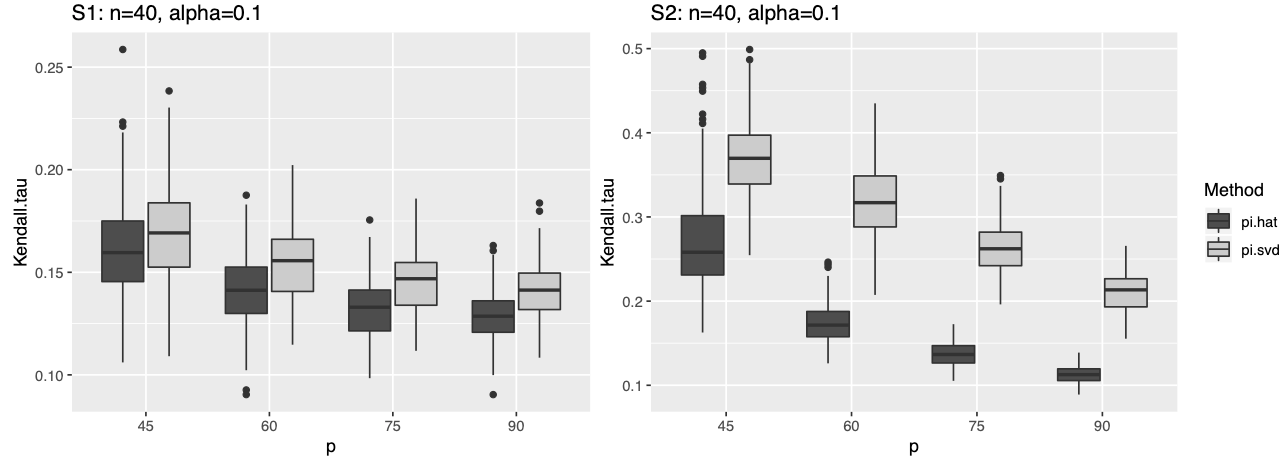}
	\includegraphics[angle=0,width=16cm]{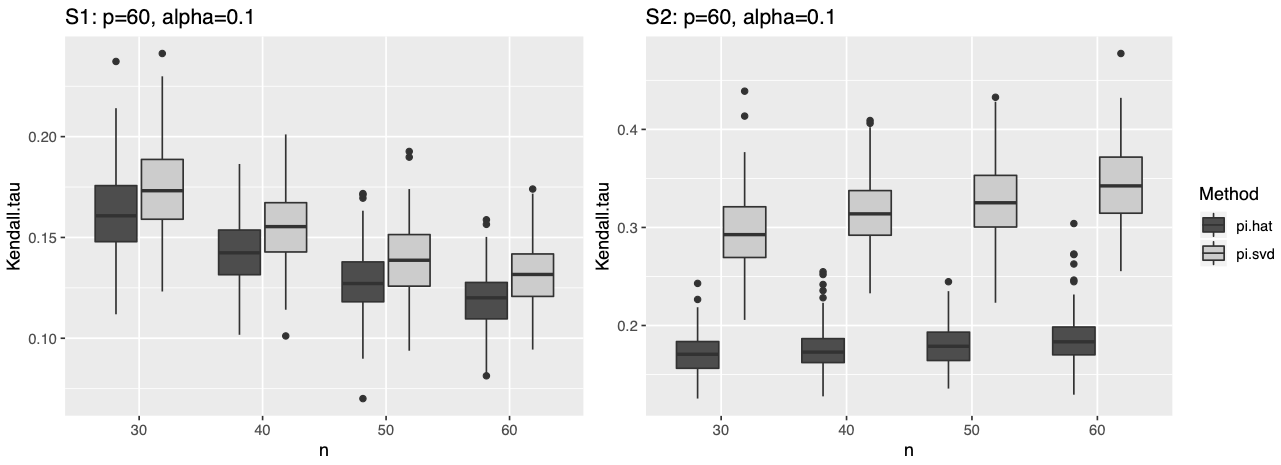}
	\caption{Boxplots of the empirical normalized Kendall's distance between the estimated permutation and true permutation under models $S_1(\alpha,p,n)$ and $S_2(\alpha,p,n)$.  $\hat{\pi}$: proposed estimator; $\pi_{svd}$: estimator based on SVD.} 
	\label{simu.1}
\end{figure}

In this section, we compare the empirical performance of our proposed estimator $\hat{\pi}$ to that of an alternative estimator $\tilde{\pi}=(\frak{r}(\hat{v}_1))^{-1}$ where $\hat{v}_1$ is the first right singular vector of $Y$. This estimator is closely related to $\hat{\pi}$ except that  it does not centralize the rows in $Y$ before estimating  its singular subspaces.  However,  this normalization step is essential in order for the resulting estimator to be  invariant to the unknown intercepts of the growth models. 
The signal matrix $\Theta=(\theta_{ij})\in \R^{n\times p}$ is generated under the following two regimes:
\begin{itemize}
	\item $S_1(\alpha,n,p)$: For any $1\le j\le p$, $\theta_{ij}=\log(1+\alpha j+\beta_i)$ for $1\le i\le n/2$ where as $\theta_{ij}=0$ for $n/2<i\le n$, $\beta_i\sim \text{Unif}(1,3)$ for all $1\le i\le n$;
	\item $S_2(\alpha,n,p)$: For any $1\le j\le p$, $\theta_{ij}=\alpha j+\beta_i$ for $i=1$ where as $\theta_{ij}=0$ for $2\le i\le n$, $\beta_i\sim \text{Unif}(1,3)$ for all $1\le i\le n$;
\end{itemize}
In each setting, we evaluate the empirical performance of each method over a range of $n$, $p$ and $\alpha$. Each setting is repeated for 200 times. The empirical normalized Kendall's tau is reported using boxplots, as shown in Figure 1 of our supplementary material. From Figure \ref{simu.1}, our proposed estimator $\hat{\pi}$ performs better than $\tilde{\pi}$ in all the settings, especially in $S_2(\alpha,n,p)$ where the signals are concentrated at one row.

\section{Supplementary Figures and Tables}

\begin{figure}[h!]
	\centering\begin{tabular}{c}
		\includegraphics[height=0.48\textheight, angle=0,width=0.60\textwidth]{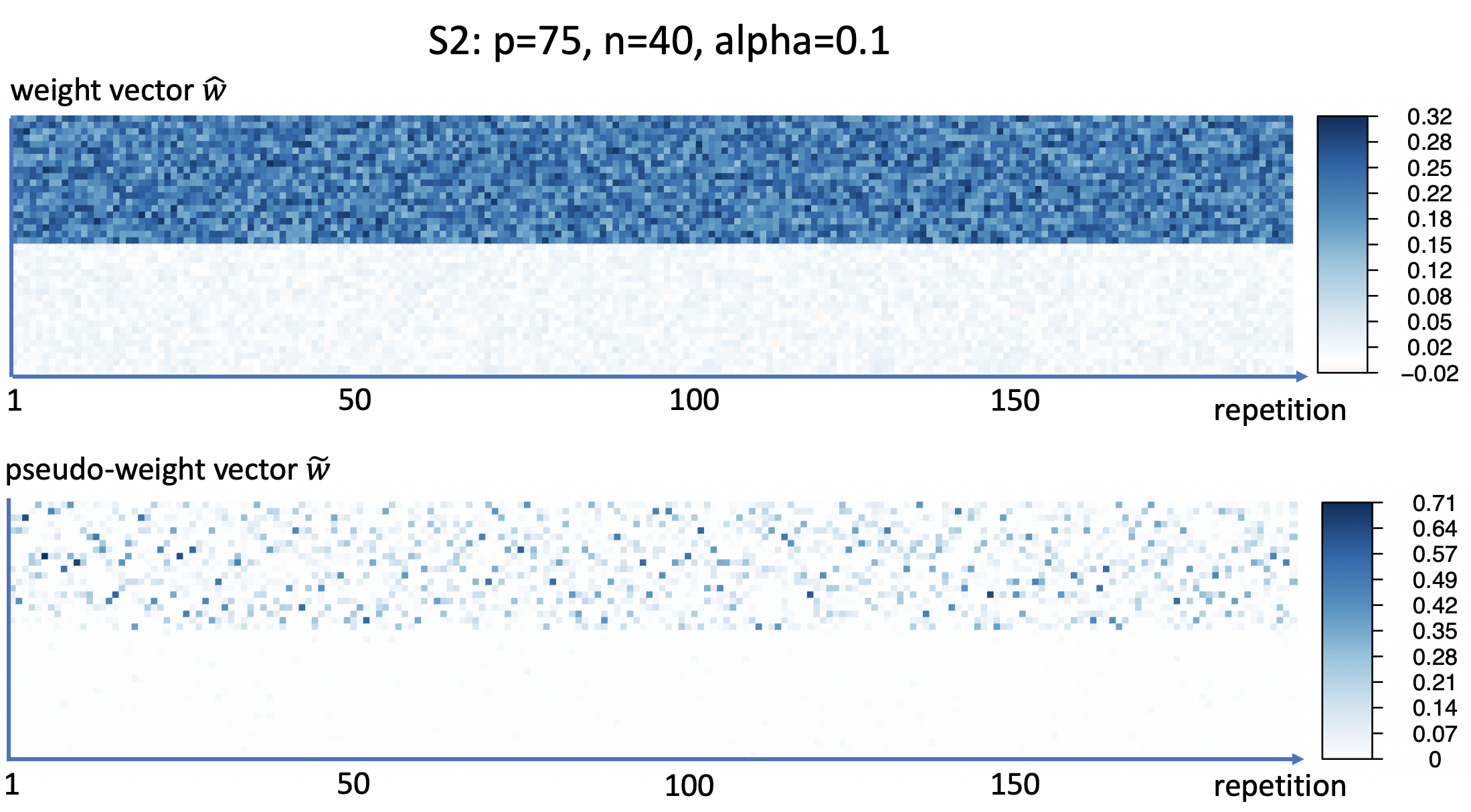}\\
		\includegraphics[height=0.48\textheight, angle=0,width=0.60\textwidth]{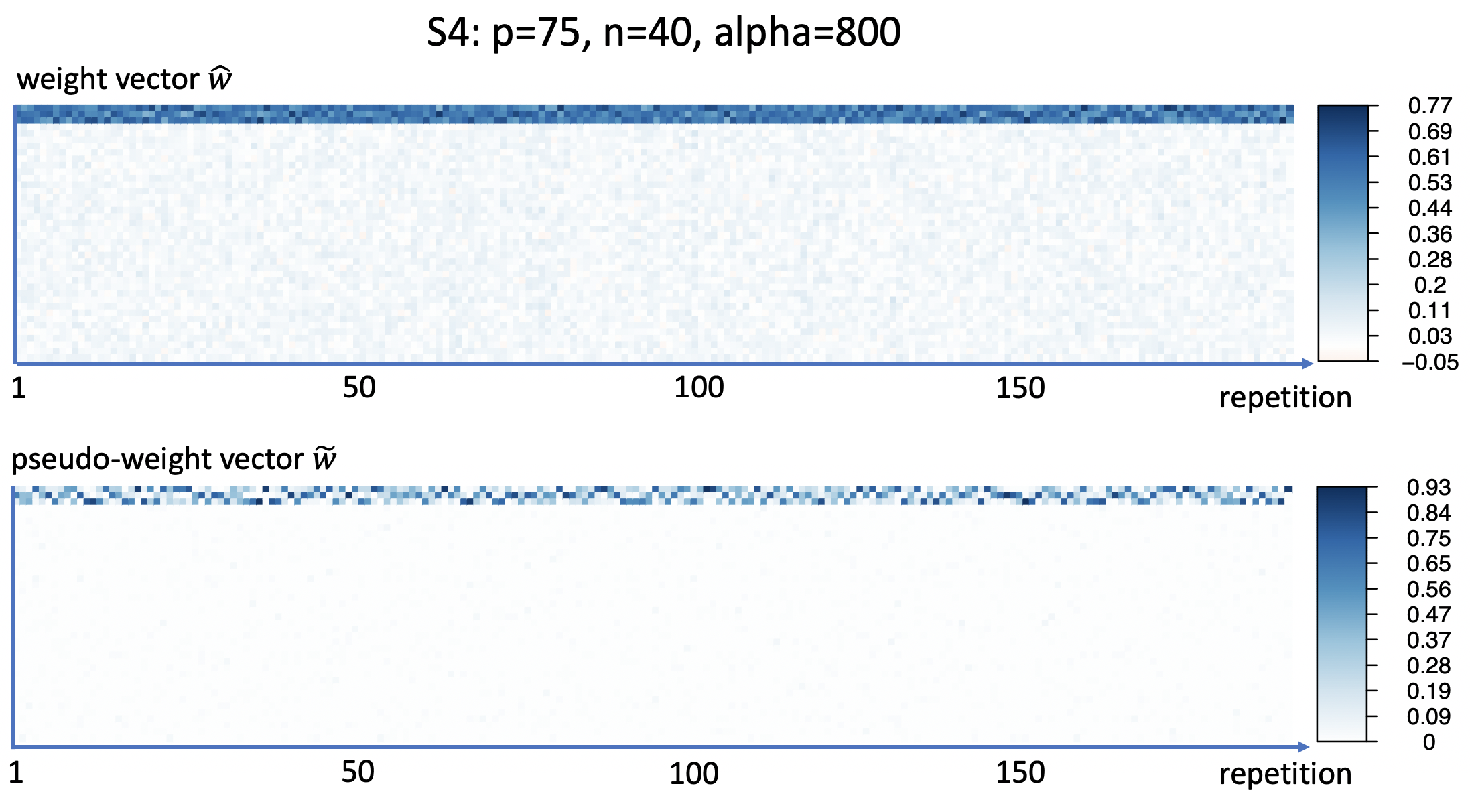}
	\end{tabular}
	\caption{The graphical representation of the weight vectors $\hat{w}$ for our proposed estimator $\hat{\pi}$, and the pseudo-weight vector $\tilde{w}$ for the estimator $\pi_{max}$ based on 200 simulations under four different models.  Each column represents an $n$ dimensional weight vector, and there are 200 columns in each plot.}
	\label{S_w1}
\end{figure}

In Figure \ref{S_w1}, the graphical representation of the weight vectors $\hat{w}$ for our proposed estimator $\hat{\pi}$, and the pseudo-weight vector $\tilde{w}$ for the estimator $\pi_{max}$ based on 200 simulations under four different models in Section 6.1 of our main paper is given.

\begin{figure}[h!]
	\centering
	\includegraphics[angle=0, height=0.75\textheight, width=0.95\textwidth]{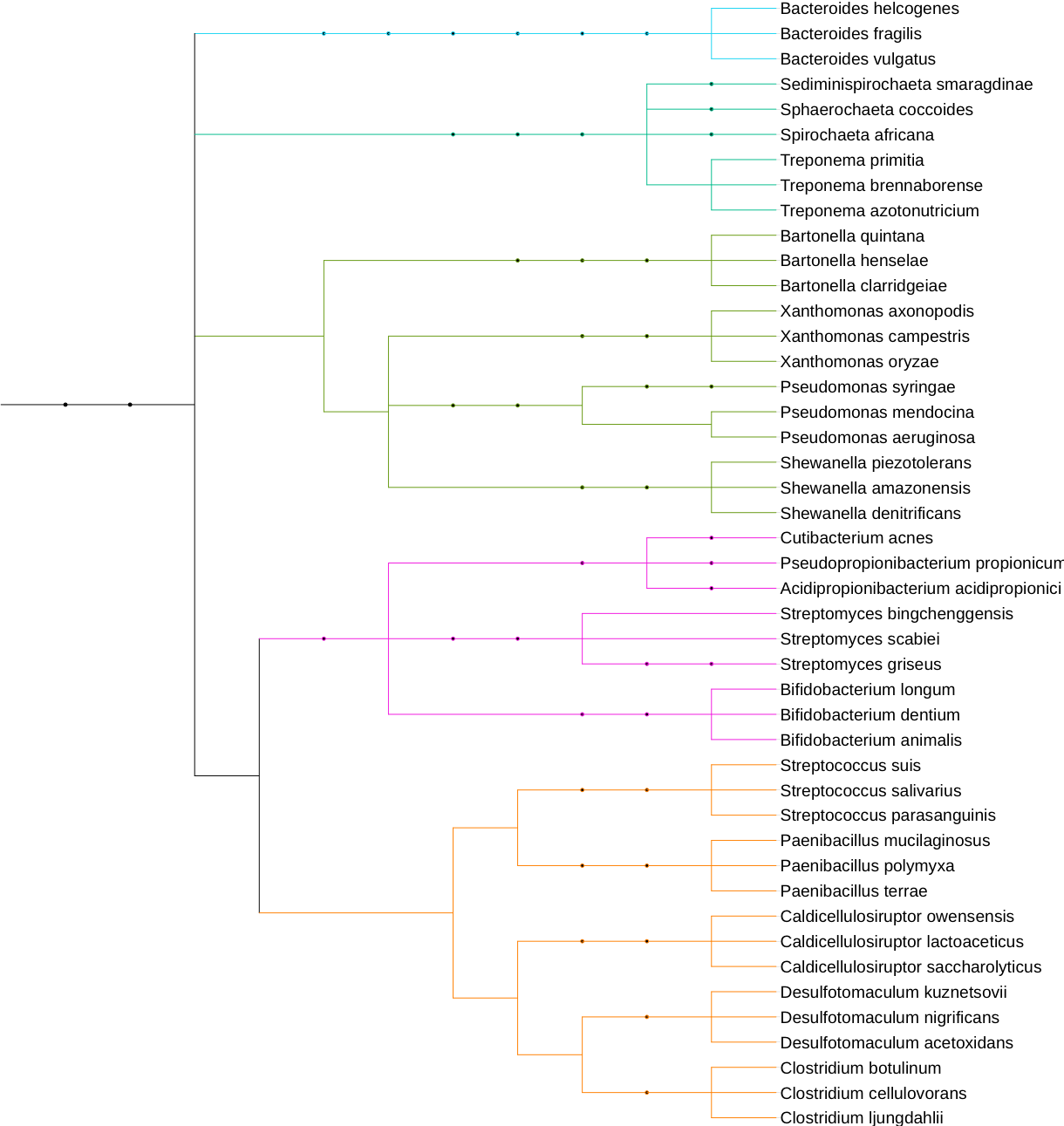}
	\caption{Taxonomic tree of 45 closely related species used in generating the shotgun metagenomic data used on in \cite{gao2018quantifying}. }
	\label{synth.tree}
\end{figure}

In Figure \ref{synth.tree}, the Taxonomic tree of 45 closely related species used in generating the shotgun metagenomic data used on in s \cite{gao2018quantifying} as well as Section 6.2 of our main paper is given.

Table \ref{tax.ann} lists the p-values of 8 contig clusters from the Wilcoxon rank sum test of the ePTRs between the responser and non-responser groups, and the taxonomic annotations with lineage scores indicating the quality of each taxonomic classification (see Section 6.3 of our main paper).

\begin{table}[h!]
	\centering
	\renewcommand\thetable{S.1}
	\caption{The p-values of 8 contig clusters from the Wilcoxon rank sum test of the ePTRs between the responser and non-responser groups, and the taxonomic annotations with lineage scores indicating the quality of each taxonomic classification.}
	\vspace{0.3cm}
	\begin{tabular*}{ 1.056 \textwidth}{lll}
		\hline 
		Contig Clusters ($n_1,n_2$) & P-values &  Taxonomic Annotations with Lineage Scores  \\  
		\hline  
		bin.004 (8,8) & 0.9592 &    Firmicutes (phylum): 0.88;	Clostridia (class): 0.78;\\
		&&	Clostridiales (order): 0.78;	Lachnospiraceae (family): 0.40\\
		bin.007 (9,11) & 0.4119  &Firmicutes (phylum): 0.96;	Clostridia (class): 0.92;\\
		&&	Clostridiales (order): 0.92; \\
		bin.016 (5,8) & 0.3543&  Firmicutes (phylum): 0.86;	Clostridia (class): 0.74;\\
		&&	Clostridiales (order): 0.74;\\
		bin.017 (7,7) & 0.2086 & Firmicutes (phylum): 0.92;	Erysipelotrichia (class): 0.47;\\
		&&Erysipelotrichales (order): 0.47;\\
		&&	Erysipelotrichaceae (family): 0.47; \\
		bin.026 (7,9)& 0.0418 & Firmicutes (phylum): 0.90;	Clostridia (class): 0.88;\\
		&&	Clostridiales (order): 0.88; \\
		bin.041 (6,6) & 0.2402 & Firmicutes (phylum): 0.95;	Clostridia (class): 0.91;\\
		&&	Clostridiales (order): 0.91;	Lachnospiraceae (family): 0.49;\\
		&&	Roseburia (genus): 0.45;\\
		bin.058 (8,15) & 0.5063 & Bacteroidetes (phylum): 0.89;	Bacteroidia (class): 0.88;\\
		&&	Bacteroidales (order): 0.88;	Bacteroidaceae (family): 0.84;\\
		&&	Bacteroides (genus): 0.84;\\
		bin.065 (5,8)& 0.0653 & Bacteroidetes (phylum): 0.88;	Bacteroidia (class): 0.88;\\
		&&	Bacteroidales (order): 0.88;	Bacteroidaceae (family): 0.85;\\
		&&	Bacteroides (genus): 0.85;\\
		\hline
	\end{tabular*}
	\label{tax.ann}
\end{table}

\label{lastpage}

\end{document}